\documentclass[11pt,twoside,a4paper]{article}

\usepackage[utf8]{inputenc}
\usepackage[T1]{fontenc}
\usepackage{amsmath,amssymb,dsfont}
\numberwithin{equation}{section}
\usepackage{microtype}
\usepackage{lmodern}
\usepackage{graphicx,tikz,pgfplots}
\usepackage{tikz-3dplot}
\usetikzlibrary{arrows,shapes,fadings,intersections,decorations.pathreplacing,positioning}
\graphicspath{{images/}}
\pgfplotsset{compat=newest}
\usepackage[hyperref,amsmath,thmmarks]{ntheorem}
\usepackage{aliascnt}
\usepackage[a4paper,centering,bindingoffset=0cm,marginpar=2cm,margin=2.5cm]{geometry}
\usepackage[pagestyles]{titlesec}
\usepackage[font=footnotesize,format=plain,labelfont=sc,textfont=sl,width=0.75\textwidth,labelsep=period]{caption}
\usepackage{comment}
\usepackage{siunitx}
\usepackage{mathtools}

\usetikzlibrary{quotes}
\usetikzlibrary{angles}
\usetikzlibrary{babel}

\usetikzlibrary{decorations.pathreplacing}
\usetikzlibrary{decorations.markings}

\usepackage{mathabx}
\tikzset{
    >=stealth',
    punkt/.style={
           rectangle,
           rounded corners,
           draw=black, very thick,
           text width=6.5em,
           minimum height=2em,
           text centered},
    pil/.style={
           ->,
           thick,
           shorten <=2pt,
           shorten >=2pt,}
}
\tikzstyle{block} = [rectangle, rounded corners, minimum width= 3cm, minimum height=1cm, text centered, draw=black, fill=blue!20,]

\tikzstyle{invisbleblock} = [minimum width= 3em, minimum height=1cm, text centered, ]

\tikzstyle{decision} = [diamond, minimum width=3cm, minimum height=1cm, text centered, draw=black, fill=orange!20]

\tikzstyle{arrow} = [thick,->,>=stealth]
\tikzstyle{line} = [thick,-]

\usepackage[boxed]{algorithm2e}

\usepackage{enumerate}

\usepackage{subcaption}
\usepackage{multirow}

\usepackage[backend=biber,maxnames=10,backref=true,hyperref=true,giveninits=true,safeinputenc]{biblatex}
\ifdefined\biberExtract
\else
\fi

\DefineBibliographyStrings{english}{%
	backrefpage = {cited on page},
	backrefpages = {cited on pages},
}

\title{Diffraction Tomography, Fourier Reconstruction, and \\Full Waveform Inversion}
\author{Florian Faucher$^{1,2}$\\
	{\footnotesize\href{mailto:florian.faucher@univie.ac.at}{florian.faucher@univie.ac.at}}
	\and Clemens Kirisits$^1$\\
	{\footnotesize\href{mailto:clemens.kirisits@univie.ac.at}{clemens.kirisits@univie.ac.at}}  
	\and Michael Quellmalz$^3$\\
	{\footnotesize\href{mailto:quellmalz@math.tu-berlin.de}{quellmalz@math.tu-berlin.de}}
	\and Otmar Scherzer$^{1,4}$\\
	{\footnotesize\href{mailto:otmar.scherzer@univie.ac.at}{otmar.scherzer@univie.ac.at}}
	\and Eric Setterqvist$^4$\\
	{\footnotesize\href{mailto:eric.setterqvist@ricam.oeaw.ac.at}{eric.setterqvist@ricam.oeaw.ac.at}}
}

\date{\today}

\DeclareFieldFormat[report]{title}{``#1''}
\DeclareFieldFormat[book]{title}{``#1''}
\AtEveryBibitem{\clearfield{url}}
\AtEveryBibitem{\clearfield{note}}

\titleformat{\section}[block]{\large\sc\filcenter}{\thesection.}{0.5ex}{}[]
\titleformat{\subsection}[runin]{\bf}{\thesubsection.}{0.5ex}{}[.]
\titleformat{\subsubsection}[runin]{\bf}{\thesubsubsection.}{0.5ex}{}[.] 

\usepackage[pdftex,colorlinks=true,linkcolor=blue,citecolor=green!50!black,urlcolor=blue,bookmarks=true,bookmarksnumbered=true,pdfusetitle]{hyperref}
\hypersetup
{
    pdfauthor={F. Faucher, C. Kirisits, M. Quellmalz, O. Scherzer, E. Setterqvist},
}

\newpagestyle{headers}
{
	\headrule
	\sethead[\footnotesize\thepage][\footnotesize\sc F.~Faucher, C.~Kirisits, M.~Quellmalz, O.~Scherzer, E.~Setterqvist][]{}{\footnotesize\sc Fourier reconstruction for diffraction tomography 
  }{\footnotesize\thepage}
	\setfoot{}{}{}
}
\newaliascnt{proposition}{lemma}

\aliascntresetthe{proposition}

\newaliascnt{corollary}{lemma}

\aliascntresetthe{corollary}

\newaliascnt{theorem}{lemma}
\newtheorem{theorem}[theorem]{Theorem}
\aliascntresetthe{theorem}

\theorembodyfont{\normalfont}
\newaliascnt{definition}{lemma}

\aliascntresetthe{definition}

\newaliascnt{assumption}{lemma}

\aliascntresetthe{assumption}

\newaliascnt{notation}{lemma}

\aliascntresetthe{notation}

\newaliascnt{example}{lemma}

\aliascntresetthe{example}

\newaliascnt{experiment}{lemma}

\aliascntresetthe{experiment}

\newaliascnt{remark}{lemma}
\newtheorem{remark}[remark]{Remark}
\aliascntresetthe{remark}

\theoremstyle{nonumberplain}
\theoremseparator{:}
\theoremheaderfont{\normalfont\itshape}

\theoremsymbol{\ensuremath{\square}}

\newcommand{\strlegend}{}

\newcommand{\misfit}  {\mathcal{J}}
\newcommand{\restrict}{\mathcal{R}}
\newcommand{\dataA}{} \newcommand{\legendA}{}
\newcommand{\dataB}{} \newcommand{\legendB}{}
\newcommand{\dataC}{} \newcommand{\legendC}{}
\newcommand{\dataD}{} \newcommand{\legendD}{}
\newcommand{\dataE}{} \newcommand{\legendE}{}
\newcommand{\dataF}{} \newcommand{\legendF}{}

\newcommand{\R}{\mathds{R}}

\newcommand{\abs}[1]{\left|#1\right|}
\newcommand{\norm}[1]{\left\|#1\right\|}

\newcommand{\e}{\mathrm e}
\let\ii\i
\renewcommand{\i}{\mathrm i}

\newcommand{\dd}{\, \mathrm{d} }

\newcommand{\bs}{{\bf s}}
\newcommand{\bx}{{\bf x}}
\newcommand{\by}{{\bf y}}
\newcommand{\bk}{{\bf k}}

\newcommand{\be}{{\bf e}}

\newcommand{\bo}{{\bf 0}}

\newcommand{\ui}{u^{\mathrm{inc}}}
\newcommand{\uborn}{u^{\mathrm{Born}}}

\newcommand{\urytov}{u^{\mathrm{Rytov}}}

\newcommand{\utot}{u^{\mathrm{tot}}}
\newcommand{\phitot}{\varphi^{\mathrm{tot}}}

\newcommand{\us}{u^{\mathrm{sca}}}
\newcommand{\phir}{\varphi^{\mathrm{Rytov}}}
\newcommand{\phii}{\varphi^{\mathrm{inc}}}

\newcommand{\phis}{\varphi^{\mathrm{sca}}}

\newcommand{\rs}{r_{\mathrm{s}}}
\newcommand{\rM}{r_{\mathrm{M}}}
\newcommand{\lM}{l_{\mathrm{M}}}

\newcommand{\supp}{\operatorname{supp}}

\newcommand{\ktran}{\mathcal{F}}

\newlength{\modelwidth} \newlength{\modelheight}
\newcommand{\modelfile}{}
\newlength{\plotwidth} \newlength{\plotheight}
\newcommand{\datafile}{}
\newcommand{\datafileB}{}

\begin{document}

\maketitle
\thispagestyle{empty}
\begin{center}
	\parbox[t]{11em}{\footnotesize
		\hspace*{-1ex}$^1$Faculty of Mathematics\\
		University of Vienna\\
		Oskar-Morgenstern-Platz 1\\
		A-1090 Vienna, Austria}
	\hspace{1em}
	\parbox[t]{11em}{\footnotesize
		\hspace*{-1ex}$^2$ Project-team Makutu \\
		    Inria Bordeaux Sud-Ouest,\\
		    E2S--UPPA, UMR CNRS 5142, 
            France}
	\hspace{1em}
	\parbox[t]{11em}{\footnotesize
		\hspace*{-1ex}$^3$Institute of Mathematics\\
		Technical University Berlin\\
		Straße des 17. Juni 136 \\
		D-10623 Berlin, Germany}
	\\[1ex]
	\parbox[t]{15em}{\footnotesize
		\hspace*{-1ex}$^4$Johann Radon Institute for Computa-\\
		\hspace*{1em}tional and Applied Mathematics \\
		\hspace*{1em}(RICAM)\\
		Altenbergerstraße 69\\
		A-4040 Linz, Austria} 

\end{center}

\begin{abstract}
In this paper, we study the mathematical imaging problem 
of diffraction tomography (DT), which is an inverse scattering 
technique used to find material properties of an object by 
illuminating it with probing waves and recording the scattered waves. 
Conventional DT relies on the Fourier diffraction theorem, which is 
applicable under the condition of weak scattering. 
However, if the object has high contrasts or is too large compared 
to the wavelength, it tends to produce multiple scattering, which
complicates the reconstruction. 
We give a survey on diffraction tomography and compare 
the reconstruction of low and high contrast objects. 
We also implement and compare the reconstruction using the full 
waveform inversion method which, contrary to the Born and Rytov 
approximations, works with the total field and is more robust to multiple scattering.
\end{abstract}

\section{Introduction}

\emph{Diffraction tomography} (DT) is a technique for reconstructing the scattering potential of an object from measurements of waves scattered by that object.
DT can be understood as an alternative to, or extension of, classical computerized tomography. In computerized tomography a crucial assumption is that the radiation, X-rays for instance, essentially propagates along straight lines through the object. The attenuated rays are recorded and can be related to material properties $f$ of the object by means of the Radon, or X-ray, transform. A central result for the inversion of this relation is the Fourier slice theorem. Roughly speaking, it says that the Fourier transformed measurements are equal to the Fourier transform of $f$ evaluated along slices through the origin, \cite{Nat86}.

The straight ray assumption of computerized tomography can be considered valid as long as the wavelength of the incident field is much smaller than the size of the relevant details in the object. As soon as the wavelength is similar to or greater than those details, for instance in situations where X-rays are replaced by visible light, diffraction effects are no longer negligible. 
As an example of a medical application, an optical diffraction experiment in \cite{sung2009optical} utilized a red laser of wavelength \SI{633}{nm} to illuminate human cells of diameter around \SI{10}{\micro\meter}, which include smaller subcellular organelles.
One way to achieve better reconstruction quality in such cases is to drop the straight ray assumption and adopt a propagation model based on the wave equation instead.

The theoretical groundwork for DT was laid more than half a century ago \cite{Wol69}. The central result derived there, sometimes called the \emph{Fourier diffraction theorem}, says that the Fourier transformed measurements of the scattered wave are equal to the Fourier transform of the scattering potential evaluated along a hemisphere. This result relies on a series of assumptions: (i) the object is immersed in a homogeneous background, (ii) the incident field is a monochromatic plane wave, (iii) the scattered wave is measured on a plane in $\R^3$, and (iv) the first Born approximation of the scattered field is valid.

On the one hand the Born approximation greatly simplifies the relationship between scattered wave and scattering potential. On the other hand, however, it generally requires the object to be weakly scattering, thus limiting the applicability of the Fourier diffraction theorem. An alternative is to assume validity of the first Rytov approximation instead \cite{IwaNag75}. While mathematically this amounts to essentially the same reconstruction problem, the underlying physical assumptions are not identical to those of the Born approximation, leading to a different range of applicability in general \cite{CheSta98,SlaKakLar84}. Nevertheless, the restriction to weakly scattering objects remains.

\emph{Full waveform inversion} (FWI) is a different approach that can overcome 
some of the limitations of the first-order methods, typically at the cost of being 
computationally more demanding. 
It relies on the iterative minimization of a cost functional which penalizes the 
misfit between measurements and forward simulations of the total field, 
cf. \cite{Bamberger1979, Lailly1983, Pratt1998, Tarantola1984, Virieux2009}.
Here, the forward model consists of the solution of the full wave equation, 
without simplification of first-order approximations. 
It results in a nonlinear minimization problem to be solved, typically with
Newton-type methods \cite{Virieux2009,Nocedal2006}.

In practical experiments, there are sometimes only measurements of the intensity, i.e., the absolute value of the complex-valued wave, available.
Different phase retrieval methods were investigated, e.g., in \cite{MalDev93,GbuWol02,HorChuOuZheYan16,BeiQue22}.
For this paper, we assume that both the phase and amplitude information is present, which can be achieved by interferometry, cf.\ \cite{WedSta95}.

\subsection*{Contribution and outline}
In this paper, we present a numerical comparison of three reconstruction approaches for diffraction tomography on simulated data, based on (i) the Born approximation, (ii) the Rytov approximation, and (iii) FWI. The setting we use for this comparison is 2D transmission imaging in a homogeneous background with (approximate) plane wave irradiation. The object is assumed to make a full turn during the experiment, providing measurements for a uniform set of incidence angles. In addition we investigate how providing additional data by varying the wavelength affects the reconstruction. The scattering potentials considered here are test phantoms of varying sizes, shapes, and contrasts. Moreover, for data generation purposes we compare several forward models.

For numerical reconstruction under the Born and Rytov approximations a well-known method is the backpropagation algorithm \cite{Dev82}, which is widely used in practice, cf.\ \cite{ODTbrain} and also \cite{FaSmLiSa17}.
Our algorithms rely on the nonuniform discrete Fourier transform (NDFT), which was used in 3D Fourier diffraction tomography yielding better results than discrete backpropagation \cite{KirQueRitSchSet21}.
Our FWI-based reconstruction uses an iterative Newton-type method on 
an $L^2$ distance between data and simulations. Here, the discretization 
of the partial differential equations associated with the wave propagation
uses the hybridizable discontinuous Galerkin method (HDG), \cite{Cockburn2009,Faucher2020adjoint}.
It is implemented, together with the inverse procedure, in the open-source parallel software \texttt{hawen}\footnote{\url{https://ffaucher.gitlab.io/hawen-website/} \label{footnote:fwi:hawen}},
\cite{Hawen2020}.

The outline of this paper is as follows. The conceptual experiment is detailed in \autoref{sec:exp}. Forward models are presented in \autoref{sec:data}, and their numerical performance is compared in \autoref{sec:data_numerics}. The Fourier diffraction theorem is formulated and discussed in \autoref{sec:coverage}. Further, \autoref{section:reconstruction-methods} covers the reconstruction algorithms used for the numerical experiments, which are presented in \autoref{sec:numerical-experiments}. A concluding discussion of our findings is given in \autoref{sec:conclusion}.

\section{Experimental setup} \label{sec:exp}

\begin{figure}[ht!]
	\begin{center}
		\begin{tikzpicture}[scale=0.5]
		
		\draw[->] (-5,0)--(5,0) node[right] {$x_1$};
		\draw[->] (0,-5.5)--(0,4.5) node[right] {$x_2$};
		
		\fill[red,opacity=0.3] (0,0) circle (2cm and 2cm);
		\fill[red,opacity=0.8] (-.8,.8) circle (.4);
		\fill[red,opacity=0.3,rotate=45] (0,0) circle (1.2cm and .3cm);
		\fill[red,opacity=1] (0.5,-0.5) circle (.15);
		\fill[red,opacity=1] (1,-0.5) circle (.15);
		\fill[red,opacity=1] (0.5,-1) circle (.15);
		\fill[red,opacity=1] (1.0,-1) circle (.15);
		\node at (3,-2) {\footnotesize{object}};
		\draw[->]  (0.43, 2.46) arc (80:10:2.5);
		
		\draw[dashed] (-5,-3.5) -- (5,-3.5);
		\draw[dashed] (-5,-4) -- (5,-4);
		\draw[dashed] (-5,-4.5) -- (5,-4.5);
		\node at (-8,-4) {\footnotesize{incident field}};
		\draw[->] (5.5,-4.5) -- (5.5,-3.5);
		\node at (6,-4) {$\be_2$};
		
		\draw[line width = 3pt] (5,3.5) -- (-5,3.5) node [left] {\footnotesize{measurement line}};
		\node at (7,3.5) {$x_2 = \rM$};
		\end{tikzpicture}
	\end{center}
\caption{Experimental setup.}
	\label{fig:trans}
\end{figure}
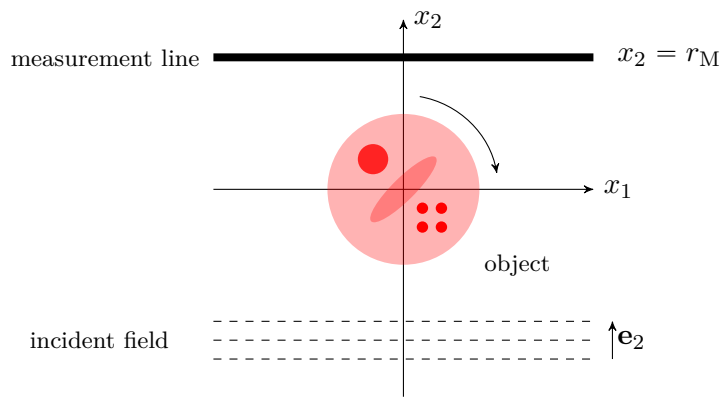

We consider the tomographic reconstruction of a two-dimensional object taking into account diffraction of the incident field. The object is assumed to be embedded in a homogeneous background and illuminated or insonified by a monochromatic plane wave.
In fact, for the computational experiments, we implement and compare 
         several approaches to approximate the incident plane wave, 
         see \autoref{sec:data} and \autoref{sec:data_numerics}.
We restrict ourselves to transmission imaging, where the incident field propagates in direction $\be_2=(0,1)^\top$ and the resulting field is measured on the line $x_2 = \rM$. The distance between the measurement line and the origin, $\rM>0$, is sufficiently large so that it does not intersect the object. 
From the measurements, we aim to reconstruct the object's scattering properties. In order to improve the reconstruction quality, we generate additional data by rotating the object or changing the incident field's wavelength. 
See \autoref{fig:trans} for an illustration of the experimental setup. 

We now introduce the physical quantities needed subsequently. Let $\lambda>0$ denote the \emph{wavelength} of the incident wave and $k_0={2\pi}/{\lambda}$ the background \emph{wave number}. Furthermore, let $n(\bx)$ denote the \emph{refractive index} at position $\bx \in \R^2$ and $n_0$ the constant refractive index of the background. 
From these quantities, we define the wave number
\begin{equation*} 
   k(\bx)\coloneqq k_0 \frac{n(\bx)}{n_0}.
\end{equation*}
Furthermore, the wave number $k$ can be equivalently expressed in terms of the
    angular frequency $\omega$ and the wave speed $c$ such that
\begin{equation} \label{eq:refractive:omega-c}
   k(\bx) \,=\, \dfrac{\omega}{c(\bx)} \, \qquad \text{and} \qquad
       k_0 \,=\, \dfrac{\omega}{c_0} \,,
\end{equation}
where $c_0$ is the constant wave speed in the background.
The \emph{scattering potential} $f$ is obtained by 
subtracting the background wave number $k_0$:
\begin{equation}\label{eq:f}
f(\bx)\coloneqq k^2(\bx)-k_0^2
= k_0^2 \left( \frac{n(\bx)^2}{n_0^2}-1\right) \, .
\end{equation}
Note that, for all practical purposes, $f$ can be assumed to be 
bounded and compactly supported in the disk 
$\mathcal B_{\rM} = \{\bx\in\R^2: \norm{\bx}<\rM\}$.

In our subsequent reconstructions with Born and Rytov approximations,
$f$ is the quantity to be reconstructed from the measured data and $k_0$ is
known. On the other hand, with FWI, we reconstruct $c$, see \autoref{rk:fwi}.
These two quantities can be related to each other via
\begin{equation} \label{eq:fwi:c_and_f}
c(\bx) \,=\, \sqrt{\dfrac{\omega^2}{k_0^2 \,+\, f(\bx)}} \, .
\end{equation}

\section{Forward Models}
\label{sec:data}
In this section, we propose several forward models for the 
experiment presented above. For all of them, the starting point 
is the system of equations
\begin{equation} \label{eq:start}
\left. 
\begin{gathered}
\big(-\Delta \,-\, k(\bx)^2 \big) \, \utot(\bx) \,=\, g(\bx), \\
\big(-\Delta \,-\, k_0^2 \big) \, \ui(\bx) \,=\, g(\bx), \\
\utot(\bx)\,=\,\ui(\bx)\,\,+\,\,\us(\bx),
\end{gathered}
\right\} \quad \bx \in \R^2.
\end{equation}
Here, $\ui$ is the given \emph{incident field}, the \emph{total field} $\utot$ is what is recorded on the measurement line $\{\bx\in\R^2: x_2=\rM\}$, and the difference between the two constitutes the \emph{scattered field} $\us$. We describe different sources $g$ in the following subsections. The scattered field $\us$ is assumed to satisfy the Sommerfeld radiation condition which requires that
\begin{equation*}
	\lim_{\norm{\bx} \to \infty} \sqrt{\norm{ \bx} } \left( \frac{\partial\us}{\partial \norm{ \bx} } - \i k_0 \us \right) = 0
\end{equation*}
uniformly for all directions $\bx/\norm{\bx}$. It guarantees that $\us$ is an outgoing wave. Further details concerning derivation and analytical properties of problems like \autoref{eq:start} can be found, for instance, in \cite{ColKre13}.

The models considered below are based on the following specifications of \autoref{eq:start}. Their numbers agree with the corresponding subsection numbers, where the models will be explained in more detail.
\begin{description}
	\item[1.\ Plane wave] No source ($g=0$) and $\ui$ is an ideal plane wave, see \autoref{eq:planewave}. 
	\item[1.1 Born model] No source, $\ui$ is an ideal plane wave and $\us$ satisfies the Born approximation.
	\item[1.2 Rytov model] No source, $\ui$ is an ideal plane wave and $\us$ satisfies the Rytov approximation.
\end{description}
In addition we propose the following two models with sources.
\begin{description}
	\item[2.1\ Point source] $g$ represents a point source located far from the object.
	\item[2.2\ Line source] $g$ represents simultaneous point sources positioned along 
	          a straight line. We refer to this configuration as a `line source'.
\end{description}
\autoref{sec:data_numerics} contains a numerical comparison of these forward models.

The proposed selection of equations is motivated in part by the availability of methods for their numerical inversion. While the Born and Rytov models can be inverted using nonuniform Fourier methods, the point and line source models are well-suited for FWI.

\subsection{Incident plane wave}\label{sec:planewave}
\emph{Monochromatic plane waves} are basic solutions $u$ of the homogeneous Helmholtz equation
\begin{equation*}
	\big(-\Delta \,-\, k_0^2 \big) \, u \,=\, 0.
\end{equation*}
They take the form $u(\bx) = \e^{\i k_0 \bx \cdot \bs}$, where the unit vector $\bs$ specifies the direction of propagation of $u$. Plane waves are widely studied in imaging applications and theory and we refer to \cite{ColKre13,Dev12,KakSla01} for further information.

In the first model we consider the incident field is a monochromatic plane wave propagating in direction $\be_2$
\begin{equation}\label{eq:planewave}
\ui(\bx) = \e^{\i k_0 x_2}.
\end{equation}
In this case, we obtain from \autoref{eq:start} the following equation for the scattered field
\begin{equation} \label{eq:uscat-rhs-uinc}
\left(-\Delta \,-\, k(\bx)^2\right) \, \us(\bx) \,=\, f(\bx) \, \e^{\i k_0 x_2} \, .
\end{equation}

\subsubsection{The Born approximation} \label{sub:Born}
\autoref{eq:uscat-rhs-uinc} can be written as
\begin{equation*}
	\left(-\Delta \,-\, k_0^2\right) \, \us(\bx) \,=\, f(\bx) \, \left( \e^{\i k_0 x_2} + \us(\bx) \right).
\end{equation*}
If the scattered field $\us$ is negligible compared to the incident field $\e^{\i k_0 x_2}$, 
we can ignore $\us$ on the right-hand side and obtain
\begin{equation} \label{eq:uborn-rhs-uinc}
(-\Delta \,-\, k_0^2) \, \uborn(\bx) \,=\, f(\bx) \, \e^{\i k_0 x_2} \,,
\end{equation}
where $\uborn$ is the \emph{(first-order) Born approximation} to the scattered field.
Supplementing this equation with the Sommerfeld radiation condition we have a unique solution 
corresponding to an outgoing wave \cite{ColKre13}. It can be written as a convolution
\begin{equation}\label{eq:conv-2d}
  \uborn(\bx) =
  \int_{\R^2} G(\bx-\by) \, f(\by)\,\e^{\i k_0 y_2} \dd\by,
\end{equation}
where $G$ is the outgoing Green's function for the Helmholtz equation. 
In $\R^2$, it is given by
\begin{equation} \label{eq:G2}
G(\bx)
= \frac{\i}{4} H^{(1)}_0(k_0 \norm{\bx})
,\qquad \bx\in\R^2\setminus\{\bo\},
\end{equation}
where $H^{(1)}_0$ denotes the zeroth order Hankel function of the first kind, see \cite[Sect.\ 3.4]{ColKre13}. We note that, in spite of a singularity at the origin, $G$ is locally integrable in $\R^2$.

The second-order Born approximation can be obtained by replacing the plane wave $\e^{\i k_0 y_2}$ in \autoref{eq:conv-2d} by the sum $\e^{\i k_0 y_2} + \uborn(\by)$. Iterating this procedure yields Born approximations of arbitrary order. For more details we refer to \cite[Sect.\ 6.2.1]{KakSla01} and \cite{Dev12}.

\subsubsection{The Rytov approximation} \label{sub:Rytov}
In this subsection, we derive an alternative approximation for the scattered field.
Introducing the complex phases $\phii$, $\phitot$ and $\phis$ according to
\begin{equation} \label{eq:exp}
\utot= \e^{\phitot}, \quad \ui = \e^{\phii}, \quad  \phitot= \phii + \phis,
\end{equation}
one can derive from \autoref{eq:start}, with $g=0$, the following relation
\begin{equation} \label{eq:totphase}
(-\Delta -k_0^2)(\ui \phis)=\left(f + \left(\nabla\phis\right)^2\right)\ui,
\end{equation}
where $\left(\nabla\phis\right)^2 = \left( \partial \phis / \partial x_1 \right)^2 + \left( \partial \phis / \partial x_2 \right)^2$.
The details of this derivation can be found, for instance, in \cite[Sect.\ 6.2.2]{KakSla01}.
Neglecting $\left(\nabla\phis\right)^2$ in \autoref{eq:totphase} we obtain
\begin{equation} \label{eq:Rytov}
(-\Delta -k_0^2)(\ui \phir)=f\ui,
\end{equation}
where $\phir$ is the \emph{Rytov approximation} to $\phis$. Note that we still assume $\ui$ to be a monochromatic plane wave, as given in \autoref{eq:planewave}. Thus the product $\ui \phir$ solves the same equation as $\uborn$.
If we define the Rytov approximation to the scattered field, $\urytov$, in analogy to \autoref{eq:exp} via
\begin{equation*}
	\urytov = \e^{\phir + \phii} - \ui,
\end{equation*}
and replace $\phir$ by $\frac{\uborn}{\ui}$, we obtain a relation between the two approximate scattered fields that can be expressed as
\begin{equation}\label{eq:Born-Rytov}
\uborn
= \ui\, \log \left( \frac{\urytov}{\ui}+1 \right).
\end{equation}
The relation between Born and Rytov in \autoref{eq:Born-Rytov}
is not unique because of the multiple branches of the complex logarithm.
In practical computations, this is addressed by a phase unwrapping as we will see in \autoref{eq:unwrap}.

\begin{remark} \label{rem:Born-Rytov}
  There have been many investigations about the validity of the Born and Rytov approximations, see, e.g., \cite{CheSta98,SlaKakLar84} or \cite[chap.\ 6]{KakSla01}.
  The Born approximation is reasonable only for a relatively (to the wavelength) small object.
  In particular, for a homogeneous cylinder of radius $a$, the Born approximation is valid if $a (n-n_0) < \lambda/4$, where $\lambda$ is the wavelength of the incident wave and $n$ is the constant refractive index inside the object.
  In contrast, the Rytov approximation only requires that
  $ n-n_0 > (\nabla \phis)^2 / k_0^2$,
  i.e., the phase change of the scattered phase $\phis$, see \autoref{eq:exp}, is small over one wavelength, but it has no direct requirements on the object size and is therefore applicable for a larger class of objects.
  The latter is also observed in numerical simulations in \cite{CheSta98}.
  However, for objects that are small and have a low contrast $n-n_0$, the Born and Rytov approximation produce approximately the same results.
\end{remark}

\subsection{Modeling the total field using line and point sources}
As an alternative to ideal incident plane waves we consider models with one or several point sources. 
If arranged the right way, the resulting incident field can resemble a monochromatic plane wave. We refer to \autoref{sec:data_numerics} for a numerical comparison of the different models presented here.

\subsubsection{Point source far from object}
In this case the right-hand side in \autoref{eq:start} is a Dirac delta function so that we obtain
\begin{equation} \label{eq:fwi:helmholtz-point}
\left\lbrace \begin{aligned}
\big(-\Delta \,-\, k(\bx)^2 \big) \, \utot_{\mathrm{P}}(\bx) \,=\, \delta(\bx -\bx_0)\, , \\
\big(-\Delta \,-\, k_0^2 \big) \, \ui_{\mathrm{P}}(\bx) \,=\, \delta(\bx -\bx_0) \, .
\end{aligned} \right. \end{equation}
If the position of the source is given by $\bx_0 = -r_0\be_2$ with $r_0>0$ sufficiently large, then, after appropriate rescaling, $\ui_{\mathrm{P}}$ approximates a plane wave with wave number $k_0$ and propagation direction $\be_2$ in a neighborhood of $\bo$. 

\subsubsection{Line source}
\label{sec:line-src}
Alternatively, we let $g$ be a sum of Dirac functions and consider
\begin{equation} \label{eq:fwi:helmholtz-full}
\left\lbrace \begin{aligned}
\big(-\Delta \,-\, k(\bx)^2 \big) \, \utot_{\mathrm{L}}(\bx) \,=\, \sum_{j=1}^{N_{\mathrm{sim}}} 
\delta(\bx-\bx_j)\, , \\
\big(-\Delta \,-\, k_0^2 \big) \, \ui_{\mathrm{L}}(\bx) \,=\, \sum_{j=1}^{N_{\mathrm{sim}}} 
\delta(\bx-\bx_j) \, ,
\end{aligned} \right. \end{equation}
where the number $N_{\mathrm{sim}}$ of simultaneous point sources should be sufficiently large. Moreover, the positions $\bx_j$ should be arranged uniformly along a line perpendicular to the propagation direction $\be_2$ of the plane wave. 
This is illustrated in \autoref{subsection:modeling:line-src}.

\section{Numerical Comparison of Forward Models}
\label{sec:data_numerics}
In this section, we numerically compare 
the forward models presented above.
For the discretization of partial differential 
equations, several approaches exist, we mention for 
instance the finite differences that approximates 
the problem on a nodal grid (e.g., \cite{Virieux1984}), 
or methods that use the variational formulation, such as
finite elements \cite{Monk2003} or discontinuous Galerkin 
methods, \cite{Hesthaven2007}. 
In our work, we use the \emph{hybridizable discontinuous Galerkin 
method} (HDG) for the discretization and refer to 
\cite{Cockburn2009,Kirby2012,Faucher2020adjoint} for more details.
The implementation precisely follows the steps prescribed 
in \cite{Faucher2020adjoint}, and it is carried out in the 
open-source parallel software \texttt{hawen},
see \cite{Hawen2020} and \autoref{footnote:fwi:hawen}.
While the propagation is assumed on infinite space,
the numerical simulations are performed on a finite 
discretization domain $\Omega\subset\R^2$, with absorbing boundary 
conditions \cite{Engquist1977} implemented to simulate
free-space.
It corresponds to the following Robin-type condition applied 
on the boundary $\Gamma$ of the discretization domain $\Omega$:
\begin{equation} \label{eq:abc}
-\mathrm{i} \, k(\bx) \, u(\bx) \,+\ \partial_n u(\bx) \,=\,0, \qquad
\text{for $x$ on $\Gamma$.}
\end{equation}
where $\partial_n u$ denotes the normal derivative of $u$. The test sample used below is a homogeneous medium encompassing a circular object of radius $4.5$ around the origin with contrast $f=1$.
This corresponds to the characteristic function
\begin{equation}\label{eq:f0}
  \mathbf 1^{\mathrm{disk}}_a(\bx) := 
  \begin{cases}
    1, & \bx\in\mathcal B_a,\\
    0, & \bx\in\R^2\setminus\mathcal B_a,
  \end{cases}
\end{equation}
of the disk $\mathcal B_a$ with radius $a>0$.
The incident plane wave has wave number $k_0 = 2\pi$.

\subsection{Modeling the scattered field assuming incident plane waves}

We consider the solutions $\us$ of \autoref{eq:uscat-rhs-uinc} and $\uborn$ of \autoref{eq:uborn-rhs-uinc}, 
both satisfying boundary condition \autoref{eq:abc}, and simulated following the HDG discretization
indicated above.
As an alternative for computing the Born approximation $\uborn$, we discretize the convolution \autoref{eq:conv-2d} 
with the Green's function $G$ given in \autoref{eq:G2}.
In particular, applying an $N\times N$ quadrature on the uniform grid $\mathcal R_N = \{ -\rM, -\rM + 2\rM/N, \dots,\rM-2\rM/N \}^2$ to \autoref{eq:conv-2d}, we obtain
\begin{equation} \label{eq:uborn-conv}
\uborn(\bx) 
\approx 
u^{\mathrm{Born}}_{\mathrm{conv},N} (\bx) \coloneqq
\left(\frac{2\rM}{N}\right)^2
\sum_{\by\in\mathcal R_N} G(\bx-\by) \, f(\by)\,\e^{\i k_0 y_2}
,\qquad\bx\in\R^2.
\end{equation}

In \autoref{fig:fwi:modeling-2d_Uinc-rhs}, we illustrate the solutions
obtained with the different formulations. 
We observe that the transmission waves contain the most energy, that is, 
waves that ``cross'' the object. On the other hand, 
the solution has very low amplitude elsewhere, including the reflected
waves.
In \autoref{fig:fwi:modeling-2d_Uinc-rhs_line}, we see that the Born
approximation leads to an incorrect amplitude of the solution, in 
particular the imaginary part. 
In addition, the imaginary part of $u^{\mathrm{Born}}_{\mathrm{conv},200}$ 
does not match the one of $u^{\mathrm{Born}}$.

\setlength{\modelwidth}{4.50cm}
\setlength{\modelheight}{\modelwidth}
\setlength{\plotwidth}{13cm}
\setlength{\plotheight}{4cm}
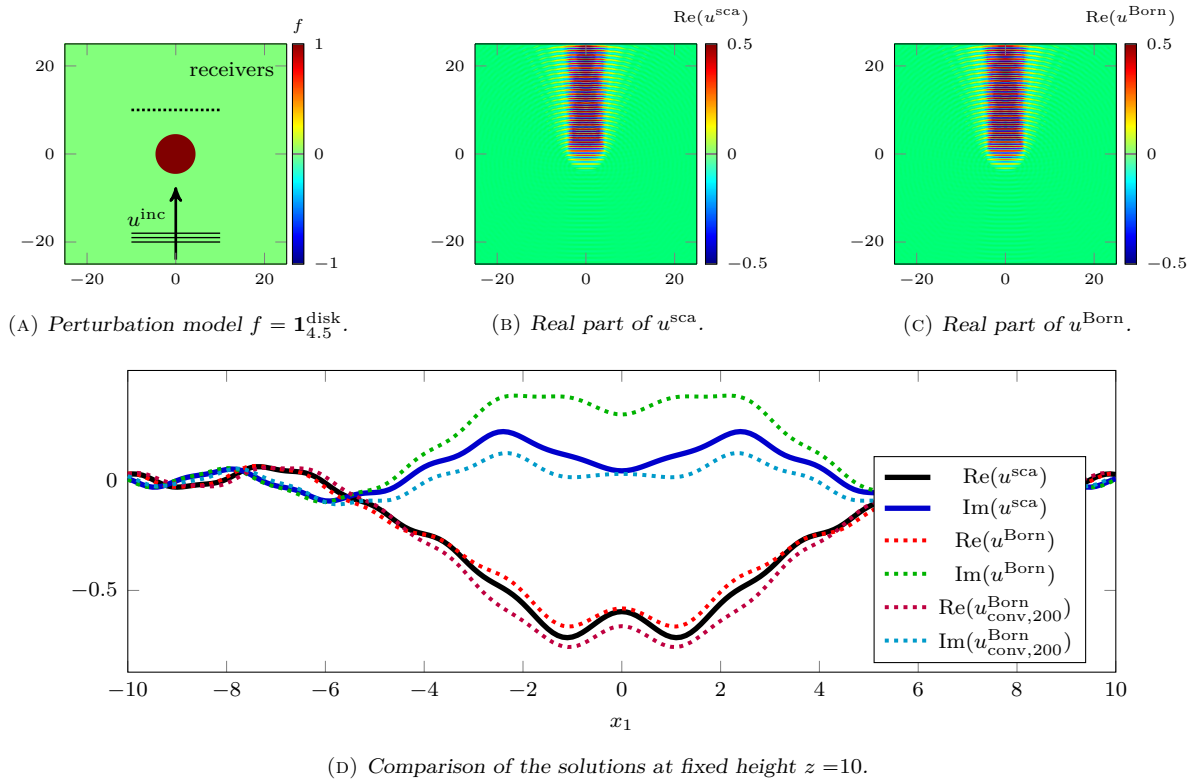
\begin{figure}[ht!] \centering
	\pgfmathsetmacro{\xminloc}{-25}\pgfmathsetmacro{\xmaxloc}{25}
	\pgfmathsetmacro{\zminloc}{-25}\pgfmathsetmacro{\zmaxloc}{25}
	\renewcommand{\modelfile}{images/fwi/models/domain100x100_f_R4.50_-1-1}
	\pgfmathsetmacro{\cmin} {-1} \pgfmathsetmacro{\cmax} { 1}
	\begin{subfigure}[t]{.31\textwidth} \centering
		\begin{tikzpicture}

\pgfmathsetmacro{\xmingb} {-50} \pgfmathsetmacro{\xmaxgb}{50}
\pgfmathsetmacro{\zmingb} {-50} \pgfmathsetmacro{\zmaxgb}{50}

\begin{axis}[
  width=\modelwidth, height=\modelheight,
  axis on top, separate axis lines,
  xmin=\xminloc, xmax=\xmaxloc, 
  ymin=\zminloc, ymax=\zmaxloc, 
  x label style={xshift=-0.0cm, yshift= 0.00cm}, 
  y label style={xshift= 0.0cm, yshift=-0.50cm},
  colormap/jet,colorbar,colorbar style={title={\tiny $f$},title style={yshift=-2mm, xshift=0mm},
  width=.15cm, xshift=-0.6em, 
  point meta min=\cmin,point meta max=\cmax,ytick={\cmin,0,\cmax},
  },  
  label style={font=\tiny},
  tick label style={font=\tiny},
  legend style={font=\tiny\selectfont},
]
\addplot [forget plot] graphics [xmin=\xmingb,xmax=\xmaxgb,ymin=\zmingb,ymax=\zmaxgb] {{\modelfile}.png};
\draw[line width=1.00, black, ->]  ( 0,-24) --  ( 0,-7.5) ;
\draw[line width=0.50, black]     (-10,-20) -- ( 10,-20)  ;
\draw[line width=0.50, black]     (-10,-19) -- ( 10,-19)  ;
\draw[line width=0.50, black]     (-10,-18) -- ( 10,-18)  ;
\node[xshift=9.5em,yshift=10.50em] (tt) {\scriptsize receivers} ;

\draw[line width=1.00, black, densely dotted]  (-10,10) --  (10,10) ;
\node[xshift=6.6em,yshift=5.40em] (tt) {\scriptsize $\ui$} ;

\end{axis}
\end{tikzpicture}%
		\caption{Perturbation model $f=\mathbf 1_{4.5}^{\mathrm{disk}}$.}
		\label{fig:fwi:modeling-2d_R4.5-f1}
	\end{subfigure}\hfill
	\pgfmathsetmacro{\cmin} {-0.50} \pgfmathsetmacro{\cmax} {0.50}
	\renewcommand{\modelfile}{images/fwi/modeling_2d_Uscat/Uscat_R4.50-f1_50x50_rhsUinc_real_scale0.50}
	\renewcommand{\strlegend}{$\mathrm{Re}(\us)$}
	\begin{subfigure}[t]{.31\textwidth} \centering
		\begin{tikzpicture}

\pgfmathsetmacro{\xmingb} {-25} \pgfmathsetmacro{\xmaxgb}{25}
\pgfmathsetmacro{\zmingb} {-25} \pgfmathsetmacro{\zmaxgb}{25}


\begin{axis}[
  width=\modelwidth, height=\modelheight,
  axis on top, separate axis lines,
  xmin=\xminloc, xmax=\xmaxloc, 
  ymin=\zminloc, ymax=\zmaxloc, 
  x label style={xshift=-0.0cm, yshift= 0.00cm}, 
  y label style={xshift= 0.0cm, yshift=-0.50cm},
  colormap/jet,colorbar,colorbar style={title={\tiny \strlegend},
  title style={yshift=-1mm, xshift=0mm},
  width=.15cm, xshift=-0.5em, 
  point meta min=\cmin,point meta max=\cmax,
  ytick={\cmin,0,\cmax},   
  },  
  label style={font=\tiny},
  tick label style={font=\tiny},
  legend style={font=\tiny\selectfont},
]
\addplot [forget plot] graphics [xmin=\xmingb,xmax=\xmaxgb,ymin=\zmingb,ymax=\zmaxgb] {{\modelfile}.png};
\end{axis}
\end{tikzpicture}%
		\caption{Real part of $\us$.}
	\end{subfigure}\hfill
	\pgfmathsetmacro{\cmin} {-0.50} \pgfmathsetmacro{\cmax} {0.50}
	\renewcommand{\modelfile}{images/fwi/modeling_2d_Uscat/Uborn_R4.50-f1_50x50_rhsUinc_real_scale0.50}
	\renewcommand{\strlegend}{$\mathrm{Re}(\uborn)$}
	\begin{subfigure}[t]{.31\textwidth} \centering
		\begin{tikzpicture}

\pgfmathsetmacro{\xmingb} {-25} \pgfmathsetmacro{\xmaxgb}{25}
\pgfmathsetmacro{\zmingb} {-25} \pgfmathsetmacro{\zmaxgb}{25}


\begin{axis}[
  width=\modelwidth, height=\modelheight,
  axis on top, separate axis lines,
  xmin=\xminloc, xmax=\xmaxloc, 
  ymin=\zminloc, ymax=\zmaxloc, 
  x label style={xshift=-0.0cm, yshift= 0.00cm}, 
  y label style={xshift= 0.0cm, yshift=-0.50cm},
  colormap/jet,colorbar,colorbar style={title={\tiny \strlegend},
  title style={yshift=-1mm, xshift=0mm},
  width=.15cm, xshift=-0.5em, 
  point meta min=\cmin,point meta max=\cmax,
  ytick={\cmin,0,\cmax},   
  },  
  label style={font=\tiny},
  tick label style={font=\tiny},
  legend style={font=\tiny\selectfont},
]
\addplot [forget plot] graphics [xmin=\xmingb,xmax=\xmaxgb,ymin=\zmingb,ymax=\zmaxgb] {{\modelfile}.png};
\end{axis}
\end{tikzpicture}%
		\caption{Real part of $\uborn$.}
	\end{subfigure} \\[1em]
	
	\renewcommand{\datafile}{images/fwi/modeling_2d_Uscat/data_rcv/data_rcv.txt}
	\renewcommand{\dataA}{uscatuinc_real} \renewcommand{\legendA}{$\mathrm{Re}(\us)$}
	\renewcommand{\dataB}{uscatuinc_imag} \renewcommand{\legendB}{$\mathrm{Im}(\us)$}
	\renewcommand{\dataC}{uborn_real}     \renewcommand{\legendC}{$\mathrm{Re}(\uborn)$}
	\renewcommand{\dataD}{uborn_imag}     \renewcommand{\legendD}{$\mathrm{Im}(\uborn)$}
	\pgfmathsetmacro{\scaleR}{1}  \pgfmathsetmacro{\scaleI}{1}
	\begin{subfigure}[t]{\textwidth} \centering
		
		\begin{tikzpicture}
		\begin{axis}[
		enlargelimits=false, 
		xlabel={$x_1$},
		enlarge y limits=true, 
		enlarge x limits=false,
		yminorticks=true,
		height=\plotheight,width=\plotwidth,
		scale only axis,
		ylabel style = {yshift =0mm, xshift=0mm},
		xlabel style = {yshift =0mm, xshift=0mm},
		clip mode=individual,
		legend columns=1,
		label style={font=\scriptsize},
		tick label style={font=\scriptsize},
		legend style={font=\scriptsize\selectfont},
		legend pos={south east}, 
		]  

		\addplot[color=black,line width=2]
		table[x expr = \thisrow{xrcv}, 
		y expr = \thisrow{\dataA}, 
		]
		{\datafile}; \addlegendentry{\legendA}    
		
		\addplot[color=blue!80!black,line width=2] 
		table[x expr = \thisrow{xrcv}, 
		y expr = \thisrow{\dataB}]
		{\datafile}; \addlegendentry{\legendB}    
		
		\addplot[color=red,line width=1.50,dotted]
		table[x expr = \thisrow{xrcv}, 
		y expr = \scaleR*\thisrow{\dataC}]
		{\datafile}; \addlegendentry{\legendC}    
		
		\addplot[color=green!70!black,line width=1.50,dotted]
		table[x expr = \thisrow{xrcv}, 
		y expr = \scaleI*\thisrow{\dataD}]
		{\datafile}; \addlegendentry{\legendD}
		
		\addplot[color=purple,line width=1.5,dotted]
		table[x index=0,y index=1]{images/fwi/modeling_2d_Uscat/data_rcv/f1-1000mHz-u.txt}; \addlegendentry{$\mathrm{Re}(\uborn_{\mathrm{conv},200})$}    
		
		\addplot[color=cyan!80!black,line width=1.5,dotted] 
		table[x index=0,y index=2]{images/fwi/modeling_2d_Uscat/data_rcv/f1-1000mHz-u.txt}; \addlegendentry{$\mathrm{Im}(\uborn_{\mathrm{conv},200})$} 
		
		%
		%
		
		\end{axis}
		\end{tikzpicture}

		\caption{Comparison of the solutions at fixed height $z=$\num{10}.}
		\label{fig:fwi:modeling-2d_Uinc-rhs_line}
	\end{subfigure}
	
	\caption{Comparison of the scattered wave $\us$ and the Born approximation $\uborn$.
		The computations are performed on a
		domain $[-25,25]\times[-25,25]$
		with boundary conditions given in 
		\autoref{eq:abc}.
		Further, we display $\uborn_{\mathrm{conv},200}$ which is the Born approximation obtained by the convolution \autoref{eq:uborn-conv}.}
	\label{fig:fwi:modeling-2d_Uinc-rhs}
\end{figure}

\subsection{Modeling the total field using line and point sources}
\label{subsection:modeling:line-src}

\begin{figure}[ht!] \centering
	\pgfmathsetmacro{\xminloc}{-25}\pgfmathsetmacro{\xmaxloc}{25}
	\pgfmathsetmacro{\zminloc}{-25}\pgfmathsetmacro{\zmaxloc}{25}
	\renewcommand{\modelfile}{images/fwi/models/domain100x100_f_R4.50_-1-1}
	\pgfmathsetmacro{\cmin} {-1} \pgfmathsetmacro{\cmax} { 1}
	\begin{subfigure}[t]{.31\textwidth}
		\begin{tikzpicture}

\pgfmathsetmacro{\xmingb} {-50} \pgfmathsetmacro{\xmaxgb}{50}
\pgfmathsetmacro{\zmingb} {-50} \pgfmathsetmacro{\zmaxgb}{50}

\begin{axis}[
  width=\modelwidth, height=\modelheight,
  axis on top, separate axis lines,
  xmin=\xminloc, xmax=\xmaxloc, 
  ymin=\zminloc, ymax=\zmaxloc, 
  x label style={xshift=-0.0cm, yshift= 0.00cm}, 
  y label style={xshift= 0.0cm, yshift=-0.50cm},
  colormap/jet,colorbar,colorbar style={title={\tiny $f$},title style={yshift=-2mm, xshift=0mm},
  width=.15cm, xshift=-0.6em, 
  point meta min=\cmin,point meta max=\cmax,ytick={\cmin,0,\cmax},
  },  
  label style={font=\tiny},
  tick label style={font=\tiny},
  legend style={font=\tiny\selectfont},
]
\addplot [forget plot] graphics [xmin=\xmingb,xmax=\xmaxgb,ymin=\zmingb,ymax=\zmaxgb] {{\modelfile}.png};
\draw[line width=0.50, black]     (-22,-15) -- ( 22,-15)  ;
\node[xshift=6.6em,yshift=6.25em] (tt) {\scriptsize line-source} ;

\end{axis}
\end{tikzpicture}%
		\caption{Computational domain $[-25,25]^2$
			for line source with perturbation model $f_{4.5}^{\mathrm{disk}}$.}
	\end{subfigure}\hfill
	\pgfmathsetmacro{\cmin} {-0.01} \pgfmathsetmacro{\cmax} {0.01}
	\renewcommand{\modelfile}{images/fwi/modeling_2d_Uscat/Ugb_R4.50-f1_rhsline_real_scale0.01}
	\renewcommand{\strlegend}{$\mathrm{Re}(\utot_{\mathrm{L}})$}
	\begin{subfigure}[t]{.31\textwidth}
		\begin{tikzpicture}

\pgfmathsetmacro{\xmingb} {-25} \pgfmathsetmacro{\xmaxgb}{25}
\pgfmathsetmacro{\zmingb} {-25} \pgfmathsetmacro{\zmaxgb}{25}


\begin{axis}[
  width=\modelwidth, height=\modelheight,
  axis on top, separate axis lines,
  xmin=\xminloc, xmax=\xmaxloc, 
  ymin=\zminloc, ymax=\zmaxloc, 
  x label style={xshift=-0.0cm, yshift= 0.00cm}, 
  y label style={xshift= 0.0cm, yshift=-0.50cm},
  colormap/jet,colorbar,colorbar style={title={\tiny \strlegend},
  title style={yshift=-1mm, xshift=0mm},
  width=.15cm, xshift=-0.5em, 
  point meta min=\cmin,point meta max=\cmax,
  ytick={\cmin,0,\cmax},   
  },  
  label style={font=\tiny},
  tick label style={font=\tiny},
  legend style={font=\tiny\selectfont},
]
\addplot [forget plot] graphics [xmin=\xmingb,xmax=\xmaxgb,ymin=\zmingb,ymax=\zmaxgb] {{\modelfile}.png};
\end{axis}
\end{tikzpicture}%
		\caption{Real part of $\utot_{\mathrm{L}}$.}
	\end{subfigure}\hfill
	\pgfmathsetmacro{\cmin} {-0.005} \pgfmathsetmacro{\cmax} {0.005}
	\renewcommand{\modelfile}{images/fwi/modeling_2d_Uscat/Uscat_R4.50-f1_rhsline_real_scale5e-3}
	\renewcommand{\strlegend}{$\mathrm{Re}(\us_{\mathrm{L}})$}
	\begin{subfigure}[t]{.31\textwidth}
		\begin{tikzpicture}

\pgfmathsetmacro{\xmingb} {-25} \pgfmathsetmacro{\xmaxgb}{25}
\pgfmathsetmacro{\zmingb} {-25} \pgfmathsetmacro{\zmaxgb}{25}


\begin{axis}[
  width=\modelwidth, height=\modelheight,
  axis on top, separate axis lines,
  xmin=\xminloc, xmax=\xmaxloc, 
  ymin=\zminloc, ymax=\zmaxloc, 
  x label style={xshift=-0.0cm, yshift= 0.00cm}, 
  y label style={xshift= 0.0cm, yshift=-0.50cm},
  colormap/jet,colorbar,colorbar style={title={\tiny \strlegend},
  title style={yshift=-1mm, xshift=0mm},
  width=.15cm, xshift=-0.5em, 
  point meta min=\cmin,point meta max=\cmax,
  ytick={\cmin,0,\cmax},   
  },  
  label style={font=\tiny},
  tick label style={font=\tiny},
  legend style={font=\tiny\selectfont},
]
\addplot [forget plot] graphics [xmin=\xmingb,xmax=\xmaxgb,ymin=\zmingb,ymax=\zmaxgb] {{\modelfile}.png};
\end{axis}
\end{tikzpicture}%
		\caption{Real part of $\us_{\mathrm{L}}$.}
	\end{subfigure} \\[-1em]
	
	\pgfmathsetmacro{\xminloc}{-500}\pgfmathsetmacro{\xmaxloc}{500}
	\pgfmathsetmacro{\zminloc}{-500}\pgfmathsetmacro{\zmaxloc}{500}
	\renewcommand{\modelfile}{images/fwi/models/domain1000x1000_f_R4.50_-1-1}
	
	\pgfmathsetmacro{\cmin} {-1} \pgfmathsetmacro{\cmax} { 1}
	\begin{subfigure}[t]{.30\textwidth}
		\begin{tikzpicture}

\pgfmathsetmacro{\xmingb} {-500} \pgfmathsetmacro{\xmaxgb}{500}
\pgfmathsetmacro{\zmingb} {-500} \pgfmathsetmacro{\zmaxgb}{500}

\begin{axis}[
  width=\modelwidth, height=\modelheight,
  axis on top, separate axis lines,
  xmin=\xminloc, xmax=\xmaxloc, 
  ymin=\zminloc, ymax=\zmaxloc, 
  x label style={xshift=-0.0cm, yshift= 0.00cm}, 
  y label style={xshift= 0.0cm, yshift=-0.50cm},
  colormap/jet,colorbar,colorbar style={title={\tiny $f$},title style={yshift=-2mm, xshift=0mm},
  width=.15cm, xshift=-0.6em, 
  point meta min=\cmin,point meta max=\cmax,ytick={\cmin,0,\cmax},
  },  
  label style={font=\tiny},
  tick label style={font=\tiny},
  legend style={font=\tiny\selectfont},
]
\addplot [forget plot] graphics [xmin=\xmingb,xmax=\xmaxgb,ymin=\zmingb,ymax=\zmaxgb] {{\modelfile}.png};
\draw[line width=1.00, black]     (-25,-480) -- ( 25,-480)  ;
\draw[line width=1.00, black]     (  0,-500) -- ( 0,-460)  ;
\node[xshift=6.00em,yshift=0.50em] (tt) {\scriptsize pt-source} ;

\draw[line width=0.40, black, ] (-50,-400) to[out=60, in=-240]  (50,-400)  ;
\draw[line width=0.40, black, ] (-80,-350) to[out=60, in=-240]  (80,-350)  ;
\draw[line width=0.40, black, ](-110,-300) to[out=60, in=-240] (110,-300)  ;
\draw[line width=0.40, black, ](-140,-250) to[out=60, in=-240] (140,-250)  ;
\draw[line width=0.40, black, ](-170,-200) to[out=60, in=-240] (170,-200)  ;
\end{axis}
\end{tikzpicture}%
		\caption{Computational domain $[-500,500]^2$
			for point source with perturbation model 
			$f_{4.5}^{\mathrm{disk}}$.}
	\end{subfigure}\hfill
	\pgfmathsetmacro{\xminloc}{-25}\pgfmathsetmacro{\xmaxloc}{25}
	\pgfmathsetmacro{\zminloc}{-25}\pgfmathsetmacro{\zmaxloc}{25}
	\pgfmathsetmacro{\cmin} {-0.03} \pgfmathsetmacro{\cmax} {0.03}
	\renewcommand{\modelfile}{images/fwi/modeling_2d_Uscat/Ugb_R4.50-f1_ptsrc_zoom-center-50x50_real_scale0.03}
	\renewcommand{\strlegend}{$\mathrm{Re}(\utot_{\mathrm{P}})$}
	\begin{subfigure}[t]{.30\textwidth}
		\begin{tikzpicture}

\pgfmathsetmacro{\xmingb} {-25} \pgfmathsetmacro{\xmaxgb}{25}
\pgfmathsetmacro{\zmingb} {-25} \pgfmathsetmacro{\zmaxgb}{25}


\begin{axis}[
  width=\modelwidth, height=\modelheight,
  axis on top, separate axis lines,
  xmin=\xminloc, xmax=\xmaxloc, 
  ymin=\zminloc, ymax=\zmaxloc, 
  x label style={xshift=-0.0cm, yshift= 0.00cm}, 
  y label style={xshift= 0.0cm, yshift=-0.50cm},
  colormap/jet,colorbar,colorbar style={title={\tiny \strlegend},
  title style={yshift=-1mm, xshift=0mm},
  width=.15cm, xshift=-0.5em, 
  point meta min=\cmin,point meta max=\cmax,
  ytick={\cmin,0,\cmax},   
  },  
  label style={font=\tiny},
  tick label style={font=\tiny},
  legend style={font=\tiny\selectfont},
]
\addplot [forget plot] graphics [xmin=\xmingb,xmax=\xmaxgb,ymin=\zmingb,ymax=\zmaxgb] {{\modelfile}.png};
\end{axis}
\end{tikzpicture}%
		\caption{Real part of $\utot_{\mathrm{P}}$ near origin.}
	\end{subfigure}\hfill
	\pgfmathsetmacro{\cmin} {-0.01} \pgfmathsetmacro{\cmax} {0.01}
	\renewcommand{\modelfile}{images/fwi/modeling_2d_Uscat/Uscat_R4.50-f1_ptsrc_zoom-center-50x50_real_scale0.01}
	\renewcommand{\strlegend}{$\mathrm{Re}(\us_{\mathrm{P}})$}
	\begin{subfigure}[t]{.30\textwidth}
		\begin{tikzpicture}

\pgfmathsetmacro{\xmingb} {-25} \pgfmathsetmacro{\xmaxgb}{25}
\pgfmathsetmacro{\zmingb} {-25} \pgfmathsetmacro{\zmaxgb}{25}


\begin{axis}[
  width=\modelwidth, height=\modelheight,
  axis on top, separate axis lines,
  xmin=\xminloc, xmax=\xmaxloc, 
  ymin=\zminloc, ymax=\zmaxloc, 
  x label style={xshift=-0.0cm, yshift= 0.00cm}, 
  y label style={xshift= 0.0cm, yshift=-0.50cm},
  colormap/jet,colorbar,colorbar style={title={\tiny \strlegend},
  title style={yshift=-1mm, xshift=0mm},
  width=.15cm, xshift=-0.5em, 
  point meta min=\cmin,point meta max=\cmax,
  ytick={\cmin,0,\cmax},   
  },  
  label style={font=\tiny},
  tick label style={font=\tiny},
  legend style={font=\tiny\selectfont},
]
\addplot [forget plot] graphics [xmin=\xmingb,xmax=\xmaxgb,ymin=\zmingb,ymax=\zmaxgb] {{\modelfile}.png};
\end{axis}
\end{tikzpicture}%
		\caption{Real part of $\us_{\mathrm{P}}$ near origin.}
	\end{subfigure} \\[1em]

	\renewcommand{\datafile}{images/fwi/modeling_2d_Uscat/data_rcv/data_rcv.txt}
	\renewcommand{\dataA}{uscatuinc_real} \renewcommand{\legendA}{$\mathrm{Re}(\us)$}
	\renewcommand{\dataB}{uscatuinc_imag} \renewcommand{\legendB}{$\mathrm{Im}(\us)$}
	\renewcommand{\dataC}{uscatline_real} \renewcommand{\legendC}{$\mathrm{Re}(\us_{\mathrm{L}})$}
	\renewcommand{\dataD}{uscatline_imag} \renewcommand{\legendD}{$\mathrm{Im}(\us_{\mathrm{L}})$}
	\renewcommand{\datafileB}{images/fwi/modeling_2d_Uscat/data_rcv/data_rcv_ptsrc.txt}
	\renewcommand{\dataE}{uscat_real} \renewcommand{\legendE}{$\mathrm{Re}(\us_{\mathrm{P}})$}
	\renewcommand{\dataF}{uscat_imag} \renewcommand{\legendF}{$\mathrm{Im}(\us_{\mathrm{P}})$}
	\pgfmathsetmacro{\scaleR}{-1.103603846935472e+02}  
	\pgfmathsetmacro{\scaleI}{-1.506037323506464e+02}
	\pgfmathsetmacro{\scaleRb}{4.991040370570111e+01}  
	\pgfmathsetmacro{\scaleIb}{2.626829205657547e+01}
	\begin{subfigure}[t]{\textwidth} \centering
\begin{tikzpicture}
\begin{axis}[
             enlargelimits=false, 
             xlabel={$x_1$},
             enlarge y limits=true, 
             enlarge x limits=false,
             yminorticks=true,
             height=\plotheight,width=\plotwidth,
             scale only axis,
             ylabel style = {yshift =0mm, xshift=0mm},
             xlabel style = {yshift =0mm, xshift=0mm},
             clip mode=individual,
             legend columns=2,
             label style={font=\scriptsize},
             tick label style={font=\scriptsize},
             legend style={font=\scriptsize\selectfont},
             legend pos={south east}, 
             ]  

     \addplot[color=black,line width=2]
              table[x expr = \thisrow{xrcv}, 
                    y expr = \thisrow{\dataA}, 
                   ]
              {\datafile}; \addlegendentry{\legendA}    

     \addplot[color=blue!80!black,line width=2] 
              table[x expr = \thisrow{xrcv}, 
                    y expr = \thisrow{\dataB}]
              {\datafile}; \addlegendentry{\legendB}    

     \addplot[color=red,line width=1.50,dotted]
              table[x expr = \thisrow{xrcv}, 
                    y expr = \scaleR*\thisrow{\dataC}]
              {\datafile}; \addlegendentry{\legendC}    

     \addplot[color=green!70!black,line width=1.50,dotted]
              table[x expr = \thisrow{xrcv}, 
                    y expr = \scaleI*\thisrow{\dataD}]
              {\datafile}; \addlegendentry{\legendD}

     \addplot[color=cyan,line width=1.50,densely dashed]
              table[x expr = \thisrow{xrcv}, 
                    y expr = \scaleRb*\thisrow{\dataE}]
              {\datafileB}; \addlegendentry{\legendE}

     \addplot[color=orange,line width=1.50,densely dashed]
              table[x expr = \thisrow{xrcv}, 
                    y expr = \scaleIb*\thisrow{\dataF}]
              {\datafileB}; \addlegendentry{\legendF}
              
\end{axis}
\end{tikzpicture}
		\caption{Comparison of the solutions at fixed height $x_2=$\num{10}.}
		\label{fig:fwi:modeling-2d_Uline-rhs_line}
	\end{subfigure}
	
	\caption{Total field $\utot_{\mathrm{L}}$ and scattered field $\us_{\mathrm{L}}$ (line source) and 
		total field $\utot_{\mathrm{P}}$ and scattered field $\us_{\mathrm{P}}$ (point source).
		The line source is composed of $N_{\mathrm{sim}} = 441$
		points at fixed height $x_2=-15$.
		The computational domain for the point source is 
		very large such that the perturbation is barely 
		visible and the source is positioned in $\bx_0 = (0,-480)^\top$.
		}
	\label{fig:fwi:modeling-2d_Uline-rhs}
\end{figure}
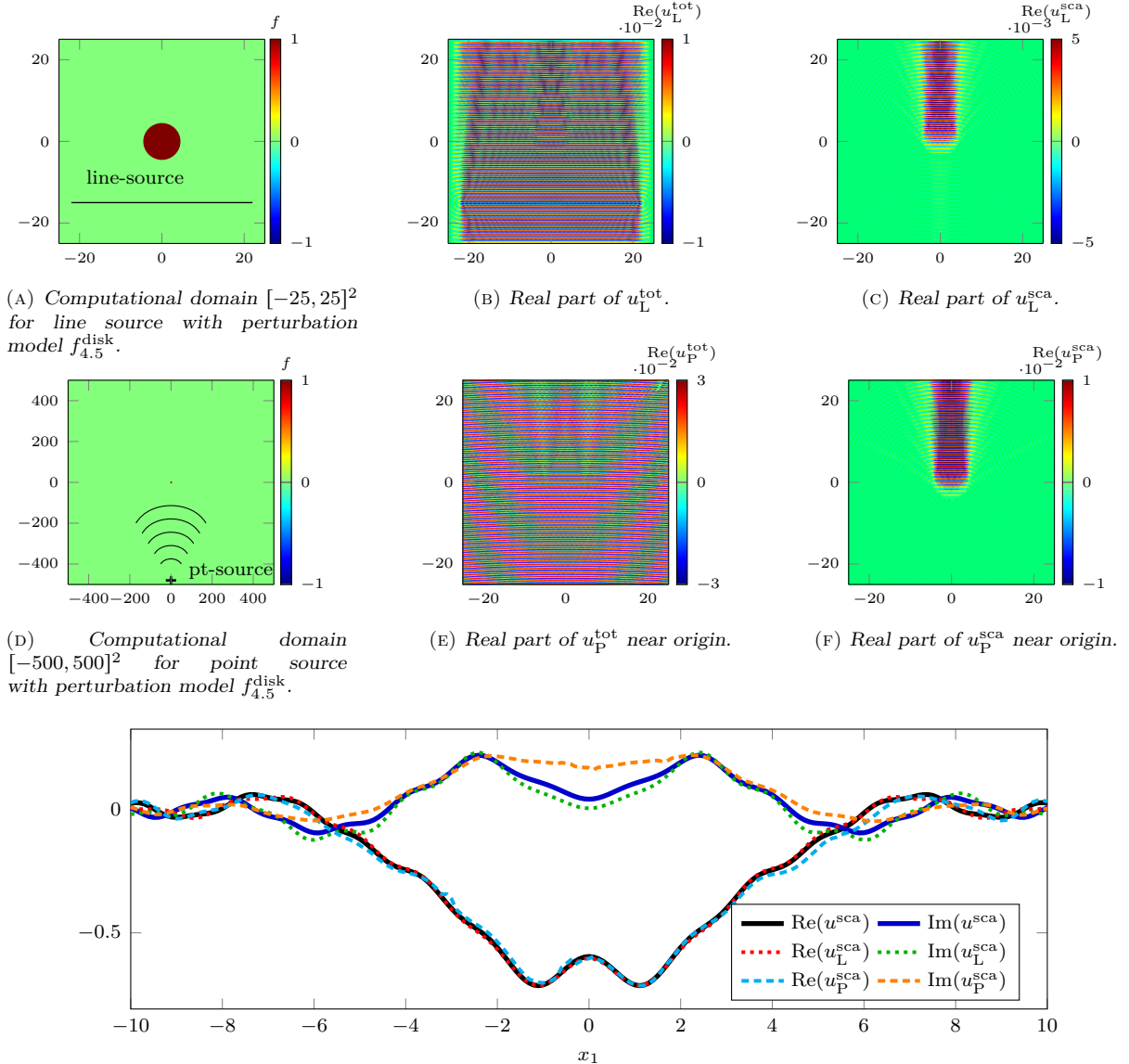

The objective is to evaluate how considering line and point sources
differs from using incident plane waves, and if the data 
obtained with both approaches are comparable.
To compare the scattered fields obtained from \autoref{eq:fwi:helmholtz-point} and \autoref{eq:fwi:helmholtz-full} with the solution $\us$ of 
\autoref{eq:uscat-rhs-uinc}, one needs to rescale according to
\begin{subequations}\begin{align*} 
	\us_{\mathrm{P}}		&=	\alpha_{\mathrm{P}} \left( \utot_{\mathrm{P}} - \ui_{\mathrm{P}} \right), \\
	\us_{\mathrm{L}}		&=	\alpha_{\mathrm{L}} \left( \utot_{\mathrm{L}} - \ui_{\mathrm{L}} \right),
\end{align*}\end{subequations}
where $\alpha_{\mathrm{P}}$ and $\alpha_{\mathrm{L}}$ are constants depending only on the number and positions of the point sources $\bx_j$.
We illustrate in \autoref{fig:fwi:modeling-2d_Uline-rhs}, where 
the line source is positioned at fixed height $x_2=-15$ and composed
of \num{441} points between $x_1=-22$ and $x_1=22$. 
For the case of a point source, we have to consider a very wide
domain, namely $[-500,500]\times[-500,500]$, and the point source is
positioned in $(x_1=0,x_2=-480)$.
In \autoref{fig:fwi:modeling-2d_Uline-rhs_line}, we plot the 
corresponding solutions on a line at height $x_2=10$, i.e., 
for measurements of transmission waves. 

We see that the simulation using the line source is very 
close to the original solution $\us$, in fact, it is a more 
accurate representation than the Born approximation pictured 
in \autoref{fig:fwi:modeling-2d_Uinc-rhs}.
The simulation using a point source positioned far away is 
also accurate, except for the middle area of the imaginary 
part.
Furthermore, the major drawback of using a single point source 
is that it necessitates a very big domain, hence largely 
increasing the computational cost.

 \section{Fourier diffraction theorem} \label{sec:coverage}
In this section we discuss the inverse problem of recovering the scattering potential from measurements of the scattered wave under the Born or Rytov approximations. 
Before stating the fundamental result in this context, see \autoref{thm:fdt} below, we have to introduce further notation.

We denote by $\ktran$ the two-dimensional Fourier transform and by $\ktran_1$ the partial Fourier transform with respect to the first coordinate,
\begin{align*}
	\mathcal{F} \phi(\bk)	&=  (2\pi)^{-1} \int_{\R^2} \phi(\bx) \e^{-\i \bx \cdot \bk} \dd \bx, \qquad \bk \in \R^2, \\
	\mathcal{F}_{1} \phi(k_1,x_2)	&=  (2\pi)^{-\frac{1}{2}} \int_{\R} \phi(\bx) \e^{-\i k_1 x_1} \dd x_1, \qquad (k_1,x_2)^\top \in \R^2.
\end{align*}
For $k_1\in[-k_0,k_0]$, we define
\begin{align*}
	\kappa(k_1) \coloneqq  \sqrt{k_0^2-k_1^2}.
\end{align*}
We can now formulate the Fourier diffraction theorem, see for instance \cite[Sect.\ 6.3]{KakSla01}, \cite[Thm.\ 3.1]{NatWue01} or \cite{Wol69}.
\begin{theorem}\label{thm:fdt}
Let $f$ be bounded with $\supp f\subset \mathcal{B}_{\rM}$. For $k_1 \in (-k_0,k_0)$, we have
	\begin{equation} \label{eq:recon}
		\ktran_{1}\uborn (k_1, \rM) = \sqrt{2 \pi} \frac{\i \e^{\i \kappa \rM}}{2 \kappa} \ktran f\left(k_1, \kappa-k_0 \right).
	\end{equation}
\end{theorem}

According to the Fourier diffraction theorem, \autoref{thm:fdt}, the measurements of the scattered wave $\uborn$ can be related to the scattering potential $f$ on a semicircle in k-space. Below we discuss how this so-called \emph{k-space coverage} of the experiment is affected by (i) rotating the object and (ii) varying the wave number $k_0$ of the incident field $\ui$.
\subsection{Rotating the object} \label{sec:rot_obj}
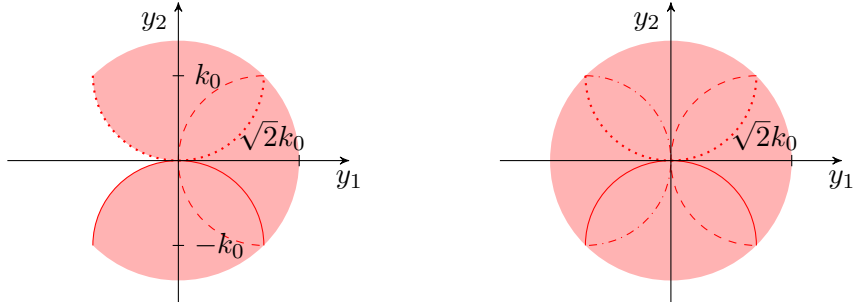
\begin{figure}[h!]
  \centering
	\begin{tikzpicture}[scale=0.75]
    		\fill[red,opacity=.3] (-1.5,1.5) arc (135:-135:2.12) -- (-1.5,-1.5) arc (180:90:1.5) -- (0,0) arc (270:180:1.5);
    		
    		\draw[thick, red, dotted] (-1.5,1.5) arc (180:360:1.5);
    		\draw[red] (-1.5,-1.5) arc (180:0:1.5);
    		\draw[red, dashed] (1.5,1.5) arc (90:270:1.5);
    		
    		\draw[->]   (-3,0) -- (3,0) node[below,black] {$y_1$};
      	\draw[->]   (0,-2.5) -- (0,2.8) node[below left,black] {$y_2$};
      	
      	\draw (2.12,-0.1) -- (2.12,0.1);
      	\node[above] at (1.65,0) {$\sqrt{2}k_0$};
      	\draw (-0.1,1.5) -- (0.1,1.5);
	    \node[right] at (0.1,1.5) {$k_0$};
	    \draw (-0.1,-1.5) -- (0.1,-1.5);
	    \node[right] at (0.1,-1.5) {$-k_0$};
    \end{tikzpicture}
    \qquad\qquad
    \begin{tikzpicture}[scale=0.75]
    		\fill[red,opacity=.3] (0,0) circle (2.12);
    		
    		\draw[thick, red, dotted] (-1.5,1.5) arc (180:360:1.5);
    		\draw[red] (-1.5,-1.5) arc (180:0:1.5);
    		\draw[red, dashed] (1.5,1.5) arc (90:270:1.5);
    		\draw[red, dashdotted] (-1.5,1.5) arc (90:-90:1.5);
    		
		\draw[->]   (-3,0) -- (3,0) node[below,black] {$y_1$};
      	\draw[->]   (0,-2.5) -- (0,2.8) node[below left,black] {$y_2$};

      	\draw (2.12,-0.1) -- (2.12,0.1);
		\node[above] at (1.65,0) {$\sqrt{2}k_0$};
    \end{tikzpicture}
    \caption{
      k-space coverage for a rotating object.
      Left: half turn, $A = [0,\pi]$. Right: full turn, $A = [0,2\pi]$. The k-space coverage (light red) is a union of infinitely many semicircles, each corresponding to a different orientation of the object. Some of the semicircles are depicted in red: solid arc ($\alpha = 0$), dashed ($\alpha = \pi/2$), dotted ($\alpha = \pi$), dash-dotted ($\alpha = 3 \pi/2$).
      \label{fig:2d-cover-half}}
  \end{figure}

Suppose the object rotates around the origin during the experiment. Then the resulting orientation-dependent scattering potential can be written as
\begin{equation*}
	 f^\alpha (\bx) = f(R_\alpha \bx), \quad \bx \in \R^2,
\end{equation*}
where $\alpha$ ranges over a (continuous or discrete) set of angles $A\subset [0,2\pi]$ and
\begin{equation*} 
	R_\alpha = \begin{pmatrix} \cos \alpha & -\sin \alpha \\ \sin \alpha & \phantom{-}\cos \alpha \end{pmatrix}
\end{equation*}
is a rotation matrix. If we let $u^\alpha$ be the Born approximation to the wave scattered by $f^\alpha$, the collected measurements are given by
\begin{equation*}
	u^\alpha(x_1,\rM), \quad x_1 \in \R, \quad \alpha \in A.
\end{equation*}
Exploiting the rotation property of the Fourier transform, $\ktran (f \circ R_\alpha) = (\ktran f) \circ R_\alpha$, we obtain
from \autoref{eq:recon}
\begin{equation*} 
	\ktran_{1}u^\alpha(k_1,\rM)= \sqrt{2\pi} \frac{\i \e^{\i \kappa \rM}}{2 \kappa} \ktran f \left(R_\alpha (k_1,\kappa-k_0)^\top \right).
\end{equation*}
Thus the k-space coverage, that is, the set of all spatial frequencies $\by \in \R^2$ at which $\ktran f$ is accessible via the Fourier diffraction theorem, is given by

\begin{equation*}
	\mathcal{Y} = \left\{ \by = R_\alpha (k_1,\kappa-k_0)^\top \in \R^2 : \abs{k_1} < k_0, \, \alpha \in A \right\}.
\end{equation*}
It consists of rotated versions (around the origin) of the semicircle $(k_1,\kappa-k_0)^\top$, $\abs{k_1}<k_0$, see \autoref{fig:2d-cover-half}.

\subsection{Varying wave number} \label{rem:2d-cover-frequencies}

Now suppose the object is illuminated or insonified by plane waves with wave numbers ranging over a set $K\subset (0,+\infty)$. Recall the definition of the scattering potential $f_{k_0} = k_0^2 (n^2/n_0^2 -1)$ from \autoref{eq:f}, but note that we have now added a subscript to indicate the dependence of $f$ on $k_0$. If the variation of the object's refractive index $n$ with $k_0 \in K$ is negligible, we can write
\begin{equation} \label{eq:f-multifreq}
  f_{k_0}(\bx) = k_0^2 f_1(\bx), \quad \bx \in \R^2.
\end{equation}
If no confusion arises, we may write $f = f_1$. Denoting by $u_{k_0}$ the Born approximation to the wave scattered by $f_{k_0}$, the resulting collection of measurements is
\begin{equation*}
	u_{k_0}(x_1,\rM), \quad x_1 \in \R, \quad k_0 \in K.
\end{equation*}
Then, according to \autoref{eq:recon}, we have
\begin{equation*}
	\ktran_{1}u_{k_0}(k_1,\rM)= \sqrt{2\pi} \frac{\i \e^{\i \kappa \rM}}{2 \kappa} k_0^2 \ktran f_1 \left( k_1, \kappa - k_0 \right).
\end{equation*}
Notice that now $\kappa$ also varies with $k_0$. The resulting k-space coverage
\begin{equation*}
	\mathcal{Y} = \left\{ \by = \left( k_1, \kappa - k_0 \right)^\top \in \R^2 : k_0 \in K, \, \abs{k_1} < k_0 \right\}
\end{equation*}
is a union of semicircles scaled and shifted (in direction of $-\be_2$) according to $k_0 \in K$.

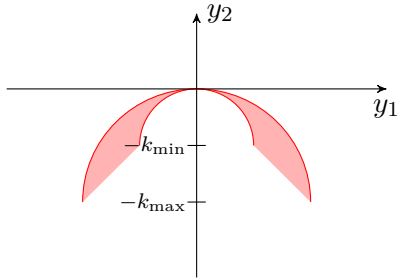
\begin{figure}[h!]\centering
	\begin{tikzpicture}[scale=1] 
	\fill[red,opacity=.3] (-1.5,-1.5) arc (180:0:1.5) -- (.75,-.75) arc (0:180:.75);
	
	\draw[red] (-1.5,-1.5) arc (180:0:1.5);
	\draw[red] (-.75,-.75) arc (180:0:.75);
	
	\draw[->]   (-2.5,0) -- (2.5,0) node[below,black] {$y_1$};
	\draw[->]   (0,-2.5) -- (0,1) node[right,black] {$y_2$};
	
	\draw[black] (-0.12,-.75) -- (0.12,-.75);
	\node[left] at (0,-0.75) {\scriptsize{$-k_{\mathrm{min}}$}};
	\draw[black] (-0.12,-1.5) -- (0.12,-1.5);
	\node[left] at (0,-1.5) {\scriptsize{$-k_{\mathrm{max}}$}};
	\end{tikzpicture}
	\caption{ k-space coverage for $k_0$ covering the interval $[k_{\mathrm{min}},k_{\mathrm{max}}]$.
		\label{fig:2d-cover-frequencies} }
\end{figure}

Consider, for example,
$K = [k_{\mathrm{min}},k_{\mathrm{max}}]$. Then the k-space coverage consists of all points $\left(y_1, y_2 \right)^\top\in\R^2$ such that
$
\abs{y_1} \le k_{\mathrm{max}}$ 
and
$$
\sqrt{k_{\mathrm{max}}^2-y_1^2} - k_{\mathrm{max}}
\ge y_2
\ge 
\begin{cases} -\abs{y_1},& \abs{y_1}\ge k_{\mathrm{min}},\\
\sqrt{k_{\mathrm{min}}^2-y_1^2} - k_{\mathrm{min}}, & \text{otherwise},
\end{cases}
$$
see \autoref{fig:2d-cover-frequencies}. Note that, in contrast to the scenarios of \autoref{sec:rot_obj}, there are large missing parts near the origin.

\subsection{Rotating the object with multiple wave numbers} \label{sec:rot_k0}
We combine the two previous observations by picking a finite set of wave numbers $K \subset(0,\infty)$ and performing a full rotation of the object for each $k_0\in K$. Let $u_{k_0}^\alpha$ be the Born approximation to the wave scattered by $f^\alpha_{k_0} = k_0^2 f_1 \circ R_\alpha $. Then the full set of measurements is given by
\begin{equation} \label{eq:u_a_k}
	u_{k_0}^\alpha(x_1,\rM), \quad x_1 \in \R, \quad \alpha \in [0,2\pi], \quad k_0\in K,
\end{equation}
and the Fourier diffraction theorem yields
\begin{equation} \label{eq:fdt_final}
	\ktran_{1}u_{k_0}^\alpha (k_1,\rM)= \sqrt{2\pi} \frac{\i \e^{\i \kappa \rM}}{2 \kappa} k_0^2 \ktran f_1 \left(R_\alpha (k_1,\kappa - k_0)^\top \right).
\end{equation}
We deduce from \autoref{sec:rot_obj} and \autoref{rem:2d-cover-frequencies} that the resulting k-space coverage $\mathcal Y$ is the union of disks with radii $\sqrt2 k_0$, all centered at the origin. Hence, $\mathcal Y$ is just the largest disk, that is, the one corresponding to the largest wave number $\max K$. However, smaller disks in k-space are covered more often, which might improve reliability of the reconstruction for noisy data.

\section{Reconstruction methods}
\label{section:reconstruction-methods}
In the following, we assume data generated by line-sources according to the 
setup described in \autoref{subsection:modeling:line-src}.
We simulate the total fields solutions to \autoref{eq:fwi:helmholtz-full}, 
which are the synthetic data used for the reconstruction.

\subsection{Reconstruction using Full Waveform Inversion}
\label{sec:fwi}

For the identification of the physical properties of the medium, 
the \emph{Full Waveform Inversion} (FWI) relies on an iterative minimization 
of a misfit functional which evaluates a distance between numerical 
simulation and measurements of the total field.
The Full Waveform Inversion method arises in the context of 
seismic inversion for sub-surface Earth imaging, 
cf.\ \cite{Bamberger1979, Lailly1983, Tarantola1984, Pratt1998, Virieux2009}, 
where the measured seismograms are compared to simulated waves.

With FWI, we invert with respect to the 
wave speed $c$, from which the wave number 
is defined according to \autoref{eq:refractive:omega-c}.
It further connects with the model perturbation $f$ according to 
\autoref{eq:fwi:c_and_f}. 
In our experiment, $c_0$ is used as an initial guess 
(i.e., we start from constant background) and then 
$c$ is inverted rather than $f$, as discussed in \autoref{rk:fwi}.
Given some measurements $\mathbf{d}$ of the 
total field, the quantitative reconstruction of 
the wave speed $c$ is performed following the 
minimization of the misfit functional $\misfit$ such that
\begin{equation} \label{eq:fwi:main-minimization}
  \min_c \, \misfit(c) \, , \qquad
  \misfit  \,=\, \mathrm{dist}\Big(\restrict \utot , \, \mathbf{d}\Big) \, ,
    \qquad \text{where $u$ solves \autoref{eq:fwi:helmholtz-full}} \, .
\end{equation}
Here, $\mathrm{dist}(\cdot)$ is a distance function to evaluate
the difference between the measurements and the simulations, 
and $\restrict$ is a linear operator to restrict the solution 
to the positions of the receivers. 
For simplicity, we do not encode a regularization term 
in~\autoref{eq:fwi:main-minimization} and refer the readers
to, e.g.,~\cite{Faucher2020EV,Kaltenbacher2018} and the 
references therein.

Several formulations of the distance function have been 
studied for FWI (in particular for seismic applications), 
such as a logarithmic criterion, \cite{Shin2007a}, the 
use of the signal's phase or amplitude, \cite{Shin2007b,Shin2007c}, 
the use of the envelope of the signal, \cite{Fichtner2008}, 
criteria based upon cross-correlation, 
\cite{LuoSchuster1991,vanLeeuwen2010,Faucher2019FRgWIGeo,Faucher2020DAS},
or optimal transport distance, \cite{Metivier2016}.
Here, we rely on a least-squares approach where the 
misfit functional is defined as the $L^2$ distance 
between the data and simulations: 
\begin{equation} \label{eq:fwi:misfit-l2}
    \misfit(c) \,:=\,  \dfrac{1}{2} \, 
    \sum_{\omega\in c_0 K} \, \sum_{\alpha\in A} \Vert \, \restrict 
    \utot(c, \omega,\, \alpha)
    \,-\, {\mathbf{d}(\omega,\alpha)} \, \Vert_{L^2(-\lM,\lM)}^2   \, , 
\end{equation}
where $\mathbf{d}(\omega,\alpha)$ refers to the measurement data of the total field at the measurement plane
with respect to the object rotated with angle $\alpha$, 
and $\utot(c, \omega,\, \alpha)$ is the solution of \autoref{eq:fwi:helmholtz-full} with $k(\bx) = \omega / c(R_\alpha^\top\bx)$.
The last term $R_\alpha^\top \bx$ encodes the rotation of the object.
We note that a rotation of the object is equivalent of the rotation of both the measurement line and the direction of the incident field.
We have encoded
a sum over the frequencies $\omega$, which are
chosen in accordance with the frequency content
available in measurements. 
In the computational experiments,
we further investigate uni- and multi-frequency
reconstructions.

The minimization of the misfit functional \autoref{eq:fwi:misfit-l2}
follows an iterative Newton-type method as depicted in \autoref{algo:FWI}.
Due to the computational cost, we use first-order information 
and avoid the Hessian computation, \cite{Virieux2009}: Namely, we 
rely on the nonlinear conjugate gradient method for the model update, 
cf.~\cite{Nocedal2006,Faucher2017}.
Furthermore, to avoid the formation of the dense Jacobian, the 
gradient of the misfit functional is computed using the adjoint-state 
method, cf.~\cite{Pratt1998,Plessix2006,Faucher2019IP,Faucher2020adjoint}.
In \autoref{algo:FWI}, we further implement a progression 
in the frequency content, which is common to mitigate the 
ill-posedness of the nonlinear inverse problem, \cite{Bunks1995}.
We further invert each frequency independently, from low 
to high, as advocated by \cite{Faucher2019IP,Faucher2020basins}.
For the implementation details using the HDG discretization,
we refer to \cite{Faucher2020adjoint}.

\begin{algorithm}[ht!]
 \KwIn{Initial wave speed model $c_0$.}
 Initiate global iteration number $\ell:= 1$;

 \For{\emph{frequency} $\omega \in c_0 K$}{
   \For{\emph{iteration} $j = 1, \, \ldots, \, n_\text{iter}$}{ 
     Compute the solution to the wave equation using 
     current wave speed model $c_\ell$ and frequency 
     $\omega$, that is, the solution to \autoref{eq:fwi:helmholtz-full}
     with $k(\bx)=\omega/c_\ell(\bx)$\;
    Evaluate the misfit functional $\misfit$ in~\autoref{eq:fwi:misfit-l2}\; 
    Compute the gradient of the misfit functional using the adjoint-state method\;
    Update the wave speed model using nonlinear conjugate gradient method to obtain $c_{\ell+1}$\;
    Update global iteration number $\ell\gets \ell+1$\;
  }
 }
 \KwOut{Approximate wave speed $c$, from which the scattering potential $f$ can be computed via Equations \ref{eq:f} and \ref{eq:refractive:omega-c}.}
 \caption{Iterative reconstruction of the wave speed model 
          following the minimization of the misfit functional. 
          At each iteration, the total field solution to 
          \autoref{eq:fwi:helmholtz-full} is computed and the gradient
          of the misfit functional is used to update the wave speed
          model.
          The algorithm stops when the prescribed number of 
          iterations is performed for all of the frequencies 
          of interest.}
 \label{algo:FWI}
\end{algorithm}

\begin{remark} \label{rk:fwi}
  In the computational experiments, the reconstruction
  with FWI assumes the availability of the total fields 
  which are solutions to \autoref{eq:fwi:helmholtz-full}, 
  and we invert with respect to the (frequency 
  independent) wave speed $c$ defined in 
  \autoref{eq:fwi:c_and_f}. 
  We could instead use the representation with 
  relation $k^2=k_0^2 + f$, and invert with respect
  to the perturbation $f$, imposing the (known) smooth 
  background $c_0$.
  Inverting with respect to $c$ rather than $f$ is 
  mainly motivated by consistency with existing
  literature in FWI \cite{Virieux2009}, in which the 
  background model ($c_0$) is usually unknown.
  Nonetheless, reformulating the minimization with respect 
  to $f$ and imposing $c_0$ could improve the efficiency of
  FWI, as advocated by the data-space reflectivity inversion 
  of \cite{Clement2001,Faucher2020basins}.

\end{remark}

\subsection{Reconstruction based on the Born and Rytov approximations} 
\label{sec:numerics_born}

In this section, we present numerical methods for the computation of the Born and Rytov approximations from \autoref{eq:uborn-rhs-uinc} and \autoref{eq:Rytov}, respectively, as well as the reconstruction of the scattering potential.
We concentrate on the case of full rotations of the object using incident waves with different wave numbers $k_0\in K$, see \autoref{sec:rot_k0}.
The tomographic reconstruction is based on the Fourier diffraction theorem, \autoref{thm:fdt}, 
and the nonuniform discrete Fourier transform.
Nonuniform Fourier methods have also been applied in computerized tomography \cite{PoSt01IMA}, magnetic resonance imaging \cite{KnKuPo}, spherical tomography \cite{HiQu15,HiQu16} or surface wave tomography \cite{HiPoQu18}.

In the following, we describe the discretization steps we apply.
For $N \in 2 \mathbb N$, let 
$$\mathcal I_N \coloneqq \left\{-\frac{N}{2} + j: j=0,\ldots,N-1\right\}.$$
We sample the scattering potential $f$ on the uniform grid $\mathcal R_N\coloneqq \frac{2\rs}{N} \mathcal I_N^2$ in the square $[-\rs,\rs]^2$ for some $\rs>0$.
We assume that we are given measurements of the Born approximation 
$$
u^{\alpha}_{k_0}(x_1,\rM),\qquad x_1\in [-\lM, \lM],
$$
for $\alpha\in A\subset[0,2\pi]$ and $k_0\in K$, cf.\ \autoref{eq:u_a_k}.
We want to reconstruct the scattering potential $f=f_1$, recall \autoref{eq:f-multifreq}, utilizing \autoref{eq:fdt_final}.
We adapt the reconstruction approach of \cite{KirQueRitSchSet21}, which is written for the 3D case.
First, we need to approximate the partial Fourier transform
\begin{equation} \label{eq:Fuborn}
\mathcal F_{1} u^{\alpha}_{k_0}(k_1,\rM)
= \frac{1}{\sqrt{2\pi}} \int_{\R} u^{\alpha}_{k_0}(x_1,\rM)\, \e^{-\i x_1k_1} \dd x_1
,\qquad 
k_1\in[-k_0,k_0].
\end{equation}
The \emph{discrete Fourier transform} (DFT) of $u(\cdot,\rM)$ on $m$ equispaced points $x_1\in({2\lM}/{m}) \mathcal I_m$ can be defined by
\begin{equation} \label{eq:dft}
\mathbf F_{1,m} u (k_1,\rM)
:=
\frac{1}{\sqrt{2\pi}}
\frac{2\lM}{m}  
\sum_{x_1\in \frac{2\lM}{m} \mathcal I_m} u(x_1,\rM)\, \e^{-\i x_1k_1} 
,\qquad k_1 \in \frac{\pi}{\lM} \mathcal I_m,
\end{equation}
which gives an approximation of \autoref{eq:Fuborn}.
Then \autoref{eq:fdt_final} yields
\begin{equation}\label{eq:fdt2}
  k_0^2 \ktran f (R_\alpha(k_1,\kappa-k_0)^\top)
  =
  -\i \sqrt{\frac{2}{\pi}} \kappa \e^{-\i\kappa \rM}
  \ktran_{1} u^{\alpha}_{k_0}(k_1,\rM)
\end{equation}
for $\abs{k_1}\le k_0$.
Considering that we sample the angle $\alpha$ on the equispaced, discrete grid $A=({2\pi}/{n_A}) \{0,1,\dots,n_A-1\}$
and some finite set $K\subset(0,\infty)$,
\autoref{eq:dft} provides an approximation of $\ktran f$ on the non-uniform grid 
\begin{multline*}
\mathcal{Y}_{m,n_A} := \Big\{ R_{\alpha}(k_1,\kappa-k_0)^\top : 
\\ k_1 \in \frac{\pi}{\lM} \mathcal I_m, \, \abs{k_1}\le k_0, \, \alpha \in \frac{2\pi}{n_A} \{0,1,\dots,n_A-1\},\, k_0\in K \Big\}
\end{multline*}
in k-space,
from which we want to reconstruct the scattering potential $f$.

Let $M$ be the cardinality of $\mathcal{Y}_{m,n_A}$. The two-dimensional \emph{nonuniform discrete Fourier transform} (NDFT) is the linear operator 
$\mathbf F_N\colon \R^{N^2} \rightarrow \R^M$ 
defined for the vector 
$\mathbf f_N 
\coloneqq
\left( f ( \bx ) \right)_{\bx \in \mathcal R_N}
$ 
elementwise by
\begin{equation}\label{eq:ndft}
  \mathbf F_N \mathbf f_N (\by) 
  \coloneqq \frac{1}{2\pi} \frac{(2\rs)^2}{N^2} \sum_{\bx \in \mathcal R_N} f(\bx) \e^{-\i \bx \cdot \by} 
  , \quad \by \in \mathcal Y_{m,n_A},
\end{equation}
see \cite[Section 7.1]{PlPoStTa18}.
It provides an approximation of the Fourier transform 
\begin{equation}\label{eq:ndft-approx}
  \ktran f(\by) \approx \mathbf F_N \mathbf f_N(\by), \qquad \by\in\mathcal Y_{m,n_A}.
\end{equation}
Solving an equation $\mathbf F_N \mathbf f_N (\by) = \bf b$ for $\mathbf f_N$ amounts to applying an \emph{inverse NDFT}, which usually utilizes an iterative method such as the conjugate gradient method on the normal equations (CGNE), see \cite{kupo04} and \cite[Section 7.6]{PlPoStTa18}.
One should be aware that the notation regarding conjugate gradient algorithms varies in the literature:
the algorithm called CGNE in \cite{Han95}, is known as CGNR in \cite{kupo04}.
Conversely, the algorithm CGME in \cite{Han95} is known as CGNE in \cite{kupo04}.

In conclusion, our method for computing $f$ given the Born approximation $\uborn$ is 
summarized in \autoref{algo:infft}.

\begin{algorithm}[ht!]
  \KwIn{
  Measurement data 
  $$u^\alpha_{k_0}(x_1,\rM),\quad x_1\in \frac{2\lM}{m} \mathcal I_m,\ \alpha\in A=\frac{2\pi}{n_A} \{0,\dots,n_A-1\},\; k_0\in K.$$}
  
  \For{ $k_0\in K$}{
    \For{$\alpha\in A$}{
    Compute $-\i \sqrt{\frac{2}{\pi}} \kappa \e^{-\i\kappa \rM}
    \mathbf F_{1,m} u^\alpha_{k_0}(k_1,\rM),$ $k_1\in \frac{\pi}{\lM}\mathcal I_m,$ with a DFT in \autoref{eq:dft}\;
  }
    }
  Solve \autoref{eq:fdt2} with \autoref{eq:ndft-approx} for $\mathbf f_N$ using the conjugate gradient method\;
  \KwOut{Approximate scattering potential $\mathbf f_N\approx (f(\bx))_{\bx\in\mathcal R_N}$.}

  \caption{Iterative reconstruction of the scattering potential $f$ based on the Born approximation using an inverse NDFT. 
  }
  \label{algo:infft}
\end{algorithm}

The Rytov approximation $\urytov$, see \autoref{eq:Rytov}, is closely related to the Born approximation,
but it has a different physical interpretation.
Assuming that the measurements arise from the Rytov approximation, we apply \autoref{eq:Born-Rytov} to obtain $\uborn$ from which we can proceed to recover $f$ as shown above.
We note that the actual implementation of \autoref{eq:Born-Rytov} requires a phase unwrapping because the complex logarithm is unique only up to adding $2\pi\i$, cf.\ \cite{ODTbrain}.
In particular, we use in the two-dimensional case
\begin{equation} \label{eq:unwrap}
  \uborn
  = \ui \left( \i \operatorname{unwrap\left(\arg\left(\frac{\urytov}{\ui}+1\right)\right)} + \ln\abs{\frac{\urytov}{\ui}+1}  \right),
\end{equation}
where $\arg$ denotes the principle argument of a complex number and $\operatorname{unwrap}$ denotes a standard unwrapping algorithm.
For the reconstruction with the Rytov approximation, we can use \autoref{algo:infft} as well, but we have to preprocess the data $u$ by \autoref{eq:unwrap}.

\section{Numerical experiments}
\label{sec:numerical-experiments}
\newcommand{\fdiskaa}{\mathbf 1^{\mathrm{disk}}_{2}  }
\newcommand{\fdiskba}{\mathbf 1^{\mathrm{disk}}_{4.5}}
\newcommand{\fdiskab}{5\cdot\mathbf 1^{\mathrm{disk}}_{2}  }
\newcommand{\fdiskbb}{5\cdot\mathbf 1^{\mathrm{disk}}_{4.5}}
\newcommand{\fdiskac}{-5\cdot\mathbf 1^{\mathrm{disk}}_{2}  }
\newcommand{\fdiskbc}{-5\cdot\mathbf 1^{\mathrm{disk}}_{4.5}}

In this section, we carry out numerical experiments 
comparing the reconstruction obtained with FWI 
(\autoref{sec:fwi}) and 
Born and Rytov approximations (\autoref{sec:numerics_born}), 
using single and multi-frequency data-sets. 
We consider different media with varying shapes and 
amplitude for the embedded objects.
Our experiments use synthetic data with added noise:
Firstly, synthetic simulations are carried out for the known 
wave speeds using the software \texttt{hawen} \cite{Hawen2020}.
The discretization relies on a fine mesh 
(usually a few hundred thousands cells in the discretized domain)
and polynomials of order $5$ to ensure accuracy.
Then, white Gaussian noise is incorporated in the synthetic
data, with a signal-to-noise ratio of \SI{50}{\dB}.
The reconstruction with FWI also relies on software \texttt{hawen}, 
but uses different discretization setups to foster the computational
time: the discretization mesh is coarser (usually a few tens of thousand
cells) and the polynomial order varies with the cells, depending on the
(local to the cell) wavelength, in order to remain as small as 
possible, as detailed in \cite{Faucher2020adjoint}. The computational
cost of FWI is further discussed in \autoref{subsection:computational-costs}.

We perform a full rotation of the object for a single or for multiple 
frequencies $\omega$ and thus wave numbers $k_0$, cf.\ \autoref{sec:rot_k0}.
For different frequencies, the scattering potential is scaled according to \autoref{eq:f-multifreq}.
We always reconstruct the rescaled scattering potential $f_1$, which we will simply denote by $f$ in the following.
In all numerical experiments, we rely on forward data generated with the forward model of line sources in \autoref{sec:line-src}.

We compare the reconstruction quality based on 
the \emph{peak signal-to-noise ratio} (PSNR) of the 
reconstruction $\mathbf g$ with respect to the ground 
truth $\mathbf f$ determined by 
\begin{equation*}
\operatorname{PSNR}(\mathbf f,\mathbf g)
\coloneqq 10 \log_{10} \frac{\max_{\bx\in\mathcal R_N} \abs{\mathbf f(\bx)}^2}{N^{-2} \sum_{\bx\in\mathcal R_N} \abs{\mathbf f(\bx) - \mathbf g(\bx)}^2},
\end{equation*}
where higher values indicate a better reconstruction quality.

\subsection{Reconstruction of circular contrast with various amplitudes and sizes}

\begin{figure}[ht!] \centering
  \pgfmathsetmacro{\xminloc}{-20}\pgfmathsetmacro{\xmaxloc}{ 20}
  \pgfmathsetmacro{\zminloc}{-20}\pgfmathsetmacro{\zmaxloc}{ 20}
  \renewcommand{\modelfile}{images/fwi/models/domain100x100_f_R4.50_-5-5}
  \pgfmathsetmacro{\cmin} {-5} \pgfmathsetmacro{\cmax} {5}
  \begin{subfigure}[t]{.31\textwidth} \centering
    \begin{tikzpicture}

\pgfmathsetmacro{\xmingb} {-50} \pgfmathsetmacro{\xmaxgb}{50}
\pgfmathsetmacro{\zmingb} {-50} \pgfmathsetmacro{\zmaxgb}{50}

\begin{axis}[
  width=\modelwidth, height=\modelheight,
  axis on top, separate axis lines,
  xmin=\xminloc, xmax=\xmaxloc, 
  ymin=\zminloc, ymax=\zmaxloc, 
  x label style={xshift=-0.0cm, yshift= 0.00cm}, 
  y label style={xshift= 0.0cm, yshift=-0.30cm},
  colormap/jet,colorbar,colorbar style={title={\scriptsize $f$},title style={yshift=-2mm, xshift=0mm},
  width=.25cm, xshift=-0.5em},
  point meta min=\cmin,point meta max=\cmax,
  label style={font=\scriptsize},
  tick label style={font=\scriptsize},
  legend style={font=\scriptsize\selectfont},
]
\addplot [forget plot] graphics [xmin=\xmingb,xmax=\xmaxgb,ymin=\zmingb,ymax=\zmaxgb] {{\modelfile}.png};
\end{axis}
\end{tikzpicture}%
\caption{Perturbation $f$ for radius \num{4.5} and amplitude~\num{5}: model $f=\fdiskbb$.}
  \end{subfigure}\hfill
  \renewcommand{\modelfile}{images/fwi/models/domain100x100_f_R4.50_-5-5_min}
  \pgfmathsetmacro{\cmin} {-5} \pgfmathsetmacro{\cmax} {5}
  \begin{subfigure}[t]{.31\textwidth} \centering
    \begin{tikzpicture}

\pgfmathsetmacro{\xmingb} {-50} \pgfmathsetmacro{\xmaxgb}{50}
\pgfmathsetmacro{\zmingb} {-50} \pgfmathsetmacro{\zmaxgb}{50}

\begin{axis}[
  width=\modelwidth, height=\modelheight,
  axis on top, separate axis lines,
  xmin=\xminloc, xmax=\xmaxloc, 
  ymin=\zminloc, ymax=\zmaxloc, 
  x label style={xshift=-0.0cm, yshift= 0.00cm}, 
  y label style={xshift= 0.0cm, yshift=-0.30cm},
  colormap/jet,colorbar,colorbar style={title={\scriptsize $f$},title style={yshift=-2mm, xshift=0mm},
  width=.25cm, xshift=-0.5em},
  point meta min=\cmin,point meta max=\cmax,
  label style={font=\scriptsize},
  tick label style={font=\scriptsize},
  legend style={font=\scriptsize\selectfont},
]
\addplot [forget plot] graphics [xmin=\xmingb,xmax=\xmaxgb,ymin=\zmingb,ymax=\zmaxgb] {{\modelfile}.png};
\end{axis}
\end{tikzpicture}%
    \caption{Perturbation $f$ for radius \num{4.5} and amplitude~\num{-5}: model $f=\fdiskbc$.}
  \end{subfigure}\hfill
  \renewcommand{\modelfile}{images/fwi/models/domain100x100_f_R2.00_-1-1}
  \pgfmathsetmacro{\cmin} {-1} \pgfmathsetmacro{\cmax} {1}
  \begin{subfigure}[t]{.31\textwidth} \centering
    \begin{tikzpicture}

\pgfmathsetmacro{\xmingb} {-50} \pgfmathsetmacro{\xmaxgb}{50}
\pgfmathsetmacro{\zmingb} {-50} \pgfmathsetmacro{\zmaxgb}{50}

\begin{axis}[
  width=\modelwidth, height=\modelheight,
  axis on top, separate axis lines,
  xmin=\xminloc, xmax=\xmaxloc, 
  ymin=\zminloc, ymax=\zmaxloc, 
  x label style={xshift=-0.0cm, yshift= 0.00cm}, 
  y label style={xshift= 0.0cm, yshift=-0.30cm},
  colormap/jet,colorbar,colorbar style={title={\scriptsize $f$},title style={yshift=-2mm, xshift=0mm},
  width=.25cm, xshift=-0.5em},
  point meta min=\cmin,point meta max=\cmax,
  label style={font=\scriptsize},
  tick label style={font=\scriptsize},
  legend style={font=\scriptsize\selectfont},
]
\addplot [forget plot] graphics [xmin=\xmingb,xmax=\xmaxgb,ymin=\zmingb,ymax=\zmaxgb] {{\modelfile}.png};
\end{axis}
\end{tikzpicture}%
    \caption{Perturbation $f$ for radius \num{2} and amplitude~\num{1}: model $f=\fdiskaa$.}
  \end{subfigure}
  \caption{Different perturbation models $f$ 
           used for the computational experiments, 
           given for frequency $\omega/(2\pi)=\num{1}$, 
           with the relation to the wave speed given in 
           \autoref{eq:fwi:c_and_f}.
           Both the size and contrast 
           vary: we consider two radii (\num{4.5} and \num{2}) 
           and three contrasts ($\num{1}$, $\num{5}$ and $\num{-5}$
           with corresponding wave speeds $c=\num{0.9876}$, $c=\num{0.9421}$ and $c=\num{1.0701}$,
            respectively), for a total of six configurations.
           The computations are carried out on the domain $[-50,50]\times[-50,50]$,
           i.e., a slightly larger setup than \autoref{fig:fwi:modeling-2d_Uline-rhs},
           and we only picture the area near the origin for clearer visualization.
           }
  \label{fig:fwi:true-models}
\end{figure}
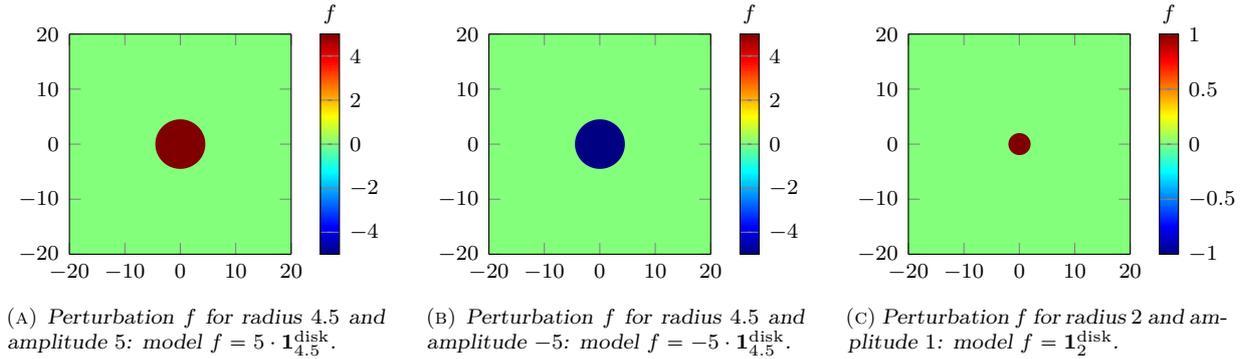

For the initial reconstruction experiments, 
we consider a circular object in a homogeneous 
background, namely the scattering potential $f$ of \autoref{eq:f0}.
We investigate different sizes and contrasts for 
the object, as shown in \autoref{fig:fwi:true-models}.
The data are generated for $n_A=\num{40}$ angles of incidence
equally partitioned between \SI{0}{\degree} to \SI{351}{\degree},
every \SI{9}{\degree},
and the measurement line is sampled on the 200 point uniform grid $ {10}^{-1}\, \mathcal I_{200} \subset [-\lM,\lM]$ 
with $\lM=10$.
Let us note that in this 
context of a circularly symmetric object, 
the data of each angle are similar and correspond to that of
\autoref{fig:fwi:modeling-2d_Uline-rhs} for $f=\mathbf 1_{4.5}^{\mathrm{disk}}$.

\subsubsection{Reconstruction using FWI with single-frequency data-sets}

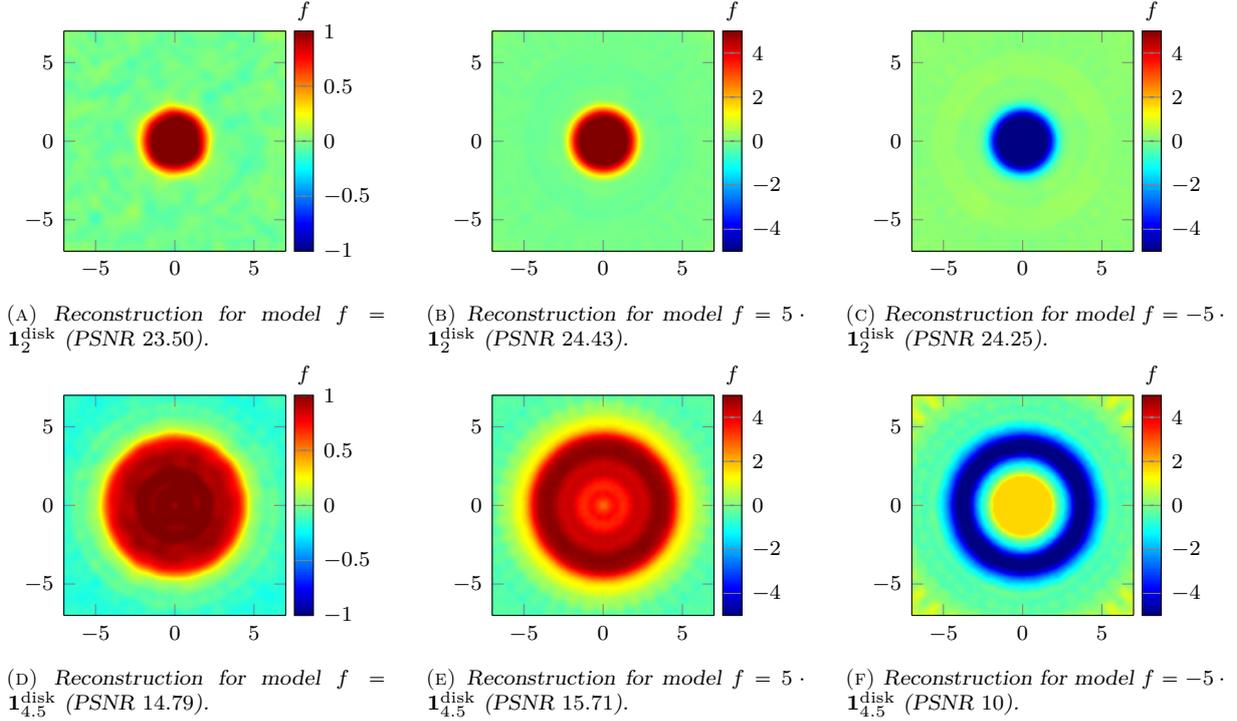
\begin{figure}[ht!] \centering
  \pgfmathsetmacro{\xminloc}{-7}\pgfmathsetmacro{\xmaxloc}{ 7}
  \pgfmathsetmacro{\zminloc}{-7}\pgfmathsetmacro{\zmaxloc}{ 7}

  \renewcommand{\modelfile}{images/fwi/reconstruction/fwi_R2.00-c0.98_unifreq1Hz-f}
  \pgfmathsetmacro{\cmin} {-1} \pgfmathsetmacro{\cmax} {1}
  \begin{subfigure}[t]{.31\textwidth} \centering
    \begin{tikzpicture}

\pgfmathsetmacro{\xmingb} {-50} \pgfmathsetmacro{\xmaxgb}{50}
\pgfmathsetmacro{\zmingb} {-50} \pgfmathsetmacro{\zmaxgb}{50}

\begin{axis}[
  width=\modelwidth, height=\modelheight,
  axis on top, separate axis lines,
  xmin=\xminloc, xmax=\xmaxloc, 
  ymin=\zminloc, ymax=\zmaxloc, 
  x label style={xshift=-0.0cm, yshift= 0.00cm}, 
  y label style={xshift= 0.0cm, yshift=-0.30cm},
  colormap/jet,colorbar,colorbar style={title={\scriptsize $f$},title style={yshift=-2mm, xshift=0mm},
  width=.25cm, xshift=-0.5em},
  point meta min=\cmin,point meta max=\cmax,
  label style={font=\scriptsize},
  tick label style={font=\scriptsize},
  legend style={font=\scriptsize\selectfont},
]
\addplot [forget plot] graphics [xmin=\xmingb,xmax=\xmaxgb,ymin=\zmingb,ymax=\zmaxgb] {{\modelfile}.png};
\end{axis}
\end{tikzpicture}%
    \caption{Reconstruction for model $f=\fdiskaa$ (PSNR \num{23.50}).}
    \label{fig:fwi:reconstruction_unifreqa}
  \end{subfigure}\hfill
  \renewcommand{\modelfile}{images/fwi/reconstruction/fwi_R2.00-c0.94_unifreq1Hz-f}
  \pgfmathsetmacro{\cmin} {-5} \pgfmathsetmacro{\cmax} {5}  
  \begin{subfigure}[t]{.31\textwidth} \centering
    \begin{tikzpicture}

\pgfmathsetmacro{\xmingb} {-50} \pgfmathsetmacro{\xmaxgb}{50}
\pgfmathsetmacro{\zmingb} {-50} \pgfmathsetmacro{\zmaxgb}{50}

\begin{axis}[
  width=\modelwidth, height=\modelheight,
  axis on top, separate axis lines,
  xmin=\xminloc, xmax=\xmaxloc, 
  ymin=\zminloc, ymax=\zmaxloc, 
  x label style={xshift=-0.0cm, yshift= 0.00cm}, 
  y label style={xshift= 0.0cm, yshift=-0.30cm},
  colormap/jet,colorbar,colorbar style={title={\scriptsize $f$},title style={yshift=-2mm, xshift=0mm},
  width=.25cm, xshift=-0.5em},
  point meta min=\cmin,point meta max=\cmax,
  label style={font=\scriptsize},
  tick label style={font=\scriptsize},
  legend style={font=\scriptsize\selectfont},
]
\addplot [forget plot] graphics [xmin=\xmingb,xmax=\xmaxgb,ymin=\zmingb,ymax=\zmaxgb] {{\modelfile}.png};
\end{axis}
\end{tikzpicture}%
    \caption{Reconstruction for model $f=\fdiskab$ 
             (PSNR \num{24.43}).}
    \label{fig:fwi:reconstruction_unifreqb}
  \end{subfigure}\hfill
  \renewcommand{\modelfile}{images/fwi/reconstruction/fwi_R2.00-c1.07_unifreq1Hz-f}
  \pgfmathsetmacro{\cmin} {-5} \pgfmathsetmacro{\cmax} {5}
  \begin{subfigure}[t]{.31\textwidth} \centering
    \begin{tikzpicture}

\pgfmathsetmacro{\xmingb} {-50} \pgfmathsetmacro{\xmaxgb}{50}
\pgfmathsetmacro{\zmingb} {-50} \pgfmathsetmacro{\zmaxgb}{50}

\begin{axis}[
  width=\modelwidth, height=\modelheight,
  axis on top, separate axis lines,
  xmin=\xminloc, xmax=\xmaxloc, 
  ymin=\zminloc, ymax=\zmaxloc, 
  x label style={xshift=-0.0cm, yshift= 0.00cm}, 
  y label style={xshift= 0.0cm, yshift=-0.30cm},
  colormap/jet,colorbar,colorbar style={title={\scriptsize $f$},title style={yshift=-2mm, xshift=0mm},
  width=.25cm, xshift=-0.5em},
  point meta min=\cmin,point meta max=\cmax,
  label style={font=\scriptsize},
  tick label style={font=\scriptsize},
  legend style={font=\scriptsize\selectfont},
]
\addplot [forget plot] graphics [xmin=\xmingb,xmax=\xmaxgb,ymin=\zmingb,ymax=\zmaxgb] {{\modelfile}.png};
\end{axis}
\end{tikzpicture}%
    \caption{Reconstruction for model $f=\fdiskac$ 
             (PSNR \num{24.25}).}
    \label{fig:fwi:reconstruction_unifreqc}
  \end{subfigure}\\

  \renewcommand{\modelfile}{images/fwi/reconstruction/fwi_R4.50-c0.98_unifreq1Hz-f}
  \pgfmathsetmacro{\cmin} {-1} \pgfmathsetmacro{\cmax} {1}
  \begin{subfigure}[t]{.31\textwidth} \centering
    \begin{tikzpicture}

\pgfmathsetmacro{\xmingb} {-50} \pgfmathsetmacro{\xmaxgb}{50}
\pgfmathsetmacro{\zmingb} {-50} \pgfmathsetmacro{\zmaxgb}{50}

\begin{axis}[
  width=\modelwidth, height=\modelheight,
  axis on top, separate axis lines,
  xmin=\xminloc, xmax=\xmaxloc, 
  ymin=\zminloc, ymax=\zmaxloc, 
  x label style={xshift=-0.0cm, yshift= 0.00cm}, 
  y label style={xshift= 0.0cm, yshift=-0.30cm},
  colormap/jet,colorbar,colorbar style={title={\scriptsize $f$},title style={yshift=-2mm, xshift=0mm},
  width=.25cm, xshift=-0.5em},
  point meta min=\cmin,point meta max=\cmax,
  label style={font=\scriptsize},
  tick label style={font=\scriptsize},
  legend style={font=\scriptsize\selectfont},
]
\addplot [forget plot] graphics [xmin=\xmingb,xmax=\xmaxgb,ymin=\zmingb,ymax=\zmaxgb] {{\modelfile}.png};
\end{axis}
\end{tikzpicture}%
    \caption{Reconstruction for model $f=\fdiskba$ 
             (PSNR \num{14.79}).}
    \label{fig:fwi:reconstruction_unifreqd}
  \end{subfigure}\hfill
  \renewcommand{\modelfile}{images/fwi/reconstruction/fwi_R4.50-c0.94_unifreq1Hz-f}
  \pgfmathsetmacro{\cmin} {-5} \pgfmathsetmacro{\cmax} {5}  
  \begin{subfigure}[t]{.31\textwidth} \centering
    \begin{tikzpicture}

\pgfmathsetmacro{\xmingb} {-50} \pgfmathsetmacro{\xmaxgb}{50}
\pgfmathsetmacro{\zmingb} {-50} \pgfmathsetmacro{\zmaxgb}{50}

\begin{axis}[
  width=\modelwidth, height=\modelheight,
  axis on top, separate axis lines,
  xmin=\xminloc, xmax=\xmaxloc, 
  ymin=\zminloc, ymax=\zmaxloc, 
  x label style={xshift=-0.0cm, yshift= 0.00cm}, 
  y label style={xshift= 0.0cm, yshift=-0.30cm},
  colormap/jet,colorbar,colorbar style={title={\scriptsize $f$},title style={yshift=-2mm, xshift=0mm},
  width=.25cm, xshift=-0.5em},
  point meta min=\cmin,point meta max=\cmax,
  label style={font=\scriptsize},
  tick label style={font=\scriptsize},
  legend style={font=\scriptsize\selectfont},
]
\addplot [forget plot] graphics [xmin=\xmingb,xmax=\xmaxgb,ymin=\zmingb,ymax=\zmaxgb] {{\modelfile}.png};
\end{axis}
\end{tikzpicture}%
    \caption{Reconstruction for model $f=\fdiskbb$ 
             (PSNR \num{15.71}).}
    \label{fig:fwi:reconstruction_unifreqe}
  \end{subfigure}\hfill
  \renewcommand{\modelfile}{images/fwi/reconstruction/fwi_R4.50-c1.07_unifreq1Hz-f}
  \pgfmathsetmacro{\cmin} {-5} \pgfmathsetmacro{\cmax} {5}
  \begin{subfigure}[t]{.31\textwidth} \centering
    \begin{tikzpicture}

\pgfmathsetmacro{\xmingb} {-50} \pgfmathsetmacro{\xmaxgb}{50}
\pgfmathsetmacro{\zmingb} {-50} \pgfmathsetmacro{\zmaxgb}{50}

\begin{axis}[
  width=\modelwidth, height=\modelheight,
  axis on top, separate axis lines,
  xmin=\xminloc, xmax=\xmaxloc, 
  ymin=\zminloc, ymax=\zmaxloc, 
  x label style={xshift=-0.0cm, yshift= 0.00cm}, 
  y label style={xshift= 0.0cm, yshift=-0.30cm},
  colormap/jet,colorbar,colorbar style={title={\scriptsize $f$},title style={yshift=-2mm, xshift=0mm},
  width=.25cm, xshift=-0.5em},
  point meta min=\cmin,point meta max=\cmax,
  label style={font=\scriptsize},
  tick label style={font=\scriptsize},
  legend style={font=\scriptsize\selectfont},
]
\addplot [forget plot] graphics [xmin=\xmingb,xmax=\xmaxgb,ymin=\zmingb,ymax=\zmaxgb] {{\modelfile}.png};
\end{axis}
\end{tikzpicture}%
    \caption{Reconstruction for model $f=\fdiskbc$ 
             (PSNR \num{10}).}
    \label{fig:fwi:reconstruction_unifreqf}
  \end{subfigure}  
  \caption{Reconstruction using iterative minimization using data
           of frequency $\omega/(2\pi)=\num{1}$ only. 
           In each cases, \num{50} iterations are performed and the initial model consists 
           in a constant background where $k_0=2\pi$. The data consist of $n_A=\num{40}$ different
           angles of incidence from \SI{0}{\degree} to \SI{351}{\degree}.}
  \label{fig:fwi:reconstruction_unifreq}
\end{figure}

We first only use data at frequency $\omega/(2\pi)=\num{1}$,
that is, wave number $k_0=2\pi$ 
for the reconstruction of the different perturbations illustrated 
in \autoref{fig:fwi:true-models}. 
With the background $k_0=2\pi$ (i.e., wave speed of \num{1}), 
it means that we 
only rely on waves with wavelength \num{1} 
when propagating in the (homogeneous) background.
Then all measurements in $\bx$ are in multiples of the wavelength.
In the case of a single frequency, only the inner loop remains 
in \autoref{algo:FWI}, and we perform \num{50} iterations.
In \autoref{fig:fwi:reconstruction_unifreq}, we picture the 
reconstruction obtained for the six different perturbations $f$.
We observe that the reconstructions of the smaller object 
of radius \num{2} (Figures \ref{fig:fwi:reconstruction_unifreqa},
 \ref{fig:fwi:reconstruction_unifreqb} and
 \ref{fig:fwi:reconstruction_unifreqc}) 
are more accurate, both in terms of the circular shape, 
and in terms of amplitude.
In the case of the larger object, 
the mild amplitude (\autoref{fig:fwi:reconstruction_unifreqd}) 
is accurately recovered, while the stronger contrasts 
(Figures \ref{fig:fwi:reconstruction_unifreqe} and
 \ref{fig:fwi:reconstruction_unifreqf})
are only partially retrieved.
Here, the outer part of the disk appears, 
but the amplitude is incorrect with a ring effect 
and incorrect values in the inner area. 
Therefore, the reconstruction using single-frequency data
is limited and its success depends on two factors: the size
of the object and its contrast.

\subsubsection{Reconstruction using FWI with multiple frequency data-sets}

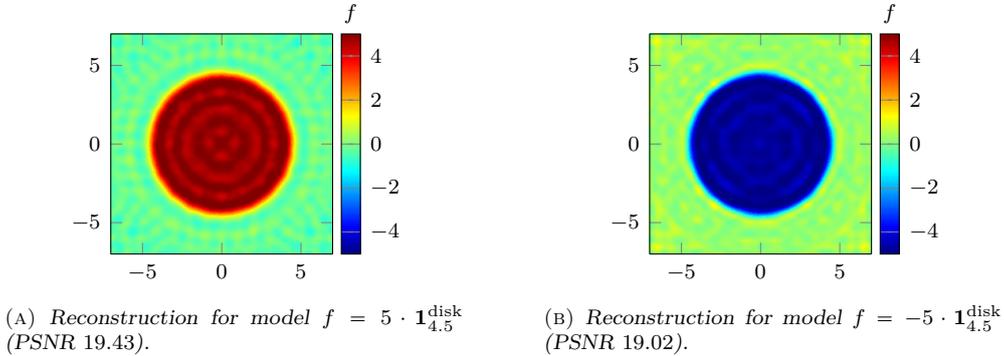
\begin{figure}[ht!] \centering
  \pgfmathsetmacro{\xminloc}{-7}\pgfmathsetmacro{\xmaxloc}{ 7}
  \pgfmathsetmacro{\zminloc}{-7}\pgfmathsetmacro{\zmaxloc}{ 7}

  \renewcommand{\modelfile}{images/fwi/reconstruction/fwi_R4.50-c0.94_multifreq-start0.2Hz-f}
  \pgfmathsetmacro{\cmin} {-5} \pgfmathsetmacro{\cmax} {5}
  \begin{subfigure}[t]{.375\textwidth} \centering
    \begin{tikzpicture}

\pgfmathsetmacro{\xmingb} {-50} \pgfmathsetmacro{\xmaxgb}{50}
\pgfmathsetmacro{\zmingb} {-50} \pgfmathsetmacro{\zmaxgb}{50}

\begin{axis}[
  width=\modelwidth, height=\modelheight,
  axis on top, separate axis lines,
  xmin=\xminloc, xmax=\xmaxloc, 
  ymin=\zminloc, ymax=\zmaxloc, 
  x label style={xshift=-0.0cm, yshift= 0.00cm}, 
  y label style={xshift= 0.0cm, yshift=-0.30cm},
  colormap/jet,colorbar,colorbar style={title={\scriptsize $f$},title style={yshift=-2mm, xshift=0mm},
  width=.25cm, xshift=-0.5em},
  point meta min=\cmin,point meta max=\cmax,
  label style={font=\scriptsize},
  tick label style={font=\scriptsize},
  legend style={font=\scriptsize\selectfont},
]
\addplot [forget plot] graphics [xmin=\xmingb,xmax=\xmaxgb,ymin=\zmingb,ymax=\zmaxgb] {{\modelfile}.png};
\end{axis}
\end{tikzpicture}%
    \caption{Reconstruction for model $f=\fdiskbb$ 
             (PSNR \num{19.43}).}
  \end{subfigure}\hspace*{1.50em}
  \renewcommand{\modelfile}{images/fwi/reconstruction/fwi_R4.50-c1.07_multifreq-start0.2Hz-f}
  \pgfmathsetmacro{\cmin} {-5} \pgfmathsetmacro{\cmax} {5}  
  \begin{subfigure}[t]{.375\textwidth} \centering
    \begin{tikzpicture}

\pgfmathsetmacro{\xmingb} {-50} \pgfmathsetmacro{\xmaxgb}{50}
\pgfmathsetmacro{\zmingb} {-50} \pgfmathsetmacro{\zmaxgb}{50}

\begin{axis}[
  width=\modelwidth, height=\modelheight,
  axis on top, separate axis lines,
  xmin=\xminloc, xmax=\xmaxloc, 
  ymin=\zminloc, ymax=\zmaxloc, 
  x label style={xshift=-0.0cm, yshift= 0.00cm}, 
  y label style={xshift= 0.0cm, yshift=-0.30cm},
  colormap/jet,colorbar,colorbar style={title={\scriptsize $f$},title style={yshift=-2mm, xshift=0mm},
  width=.25cm, xshift=-0.5em},
  point meta min=\cmin,point meta max=\cmax,
  label style={font=\scriptsize},
  tick label style={font=\scriptsize},
  legend style={font=\scriptsize\selectfont},
]
\addplot [forget plot] graphics [xmin=\xmingb,xmax=\xmaxgb,ymin=\zmingb,ymax=\zmaxgb] {{\modelfile}.png};
\end{axis}
\end{tikzpicture}%
    \caption{Reconstruction for model $f=\fdiskbc$ 
             (PSNR \num{19.02}).}
  \end{subfigure}\hfill

  \caption{Reconstruction using multi-frequency data from 
           $\omega/(2\pi)=\num{0.2}$ to $\omega/(2\pi)=\num{1}$. 
           The initial model consists in a constant wave speed $c_0=1$.
           The data consist of $n_A=\num{40}$ different angles of incidence from 
           \SI{0}{\degree} to \SI{351}{\degree}.}
  \label{fig:fwi:reconstruction_multifreq}
\end{figure}

The difficulty of recovering a large object with a 
strong contrast can be mitigated by the use of multi-frequency
data-sets, allowing a multi-scale reconstruction \cite{Bunks1995,Faucher2020basins}.
We carry out the iterative reconstruction using 
increasing frequencies, starting with $\omega/(2\pi)=\num{0.2}$ and
up to $\omega/(2\pi)=\num{1}$.
Following \cite{Faucher2020basins}, we use a sequential
progression, that is, every frequency is inverted 
separately.
The reconstructions for the object of radius \num{4.5}
and contrast $f=\pm\num{5}$ are pictured 
in \autoref{fig:fwi:reconstruction_multifreq}.
Contrary to the case of a single frequency
(see~\autoref{fig:fwi:reconstruction_unifreqd}), the 
reconstruction is now accurate and stable: the 
amplitude is accurately retrieved and 
the circular shape is clear, avoiding the circular 
artifacts observed in Figures
\ref{fig:fwi:reconstruction_unifreqe}
and \ref{fig:fwi:reconstruction_unifreqf}.

\begin{remark}
  It is possible to recover the model with a single frequency, 
  that needs to be carefully chosen depending on the size of 
  the object and the amplitude of the contrast. 
  We have seen in \autoref{fig:fwi:reconstruction_unifreq} that
  the frequency $\omega/(2\pi)=\num{1}$ is sufficient for the 
  object of radius \num{2} but for the radius 
  \num{4.5}, we need a lower frequency (i.e., larger wavelength) 
  to uncover the larger object. 
  We illustrate in \autoref{fig:fwi:reconstruction_unifreq0.7Hz}
  the reconstruction using data at only $\omega/(2\pi)={0.7}$, where
  we see that the shape and contrast are retrieved accurately.
  Nonetheless, it is hard to predict this frequency a-priori and
  we believe it remains more natural to use multiple frequencies 
  (when available in the data), to ensure the robustness of the 
  algorithm.

  \begin{figure}[ht!] \centering
  \pgfmathsetmacro{\xminloc}{-7}\pgfmathsetmacro{\xmaxloc}{ 7}
  \pgfmathsetmacro{\zminloc}{-7}\pgfmathsetmacro{\zmaxloc}{ 7}
  \renewcommand{\modelfile}{images/fwi/reconstruction/fwi_R4.50-c0.94_unifreq0.7Hz-f}
  \pgfmathsetmacro{\cmin} {-5} \pgfmathsetmacro{\cmax} {5}
  \begin{subfigure}[t]{.375\textwidth} \centering
    \begin{tikzpicture}

\pgfmathsetmacro{\xmingb} {-50} \pgfmathsetmacro{\xmaxgb}{50}
\pgfmathsetmacro{\zmingb} {-50} \pgfmathsetmacro{\zmaxgb}{50}

\begin{axis}[
  width=\modelwidth, height=\modelheight,
  axis on top, separate axis lines,
  xmin=\xminloc, xmax=\xmaxloc, 
  ymin=\zminloc, ymax=\zmaxloc, 
  x label style={xshift=-0.0cm, yshift= 0.00cm}, 
  y label style={xshift= 0.0cm, yshift=-0.30cm},
  colormap/jet,colorbar,colorbar style={title={\scriptsize $f$},title style={yshift=-2mm, xshift=0mm},
  width=.25cm, xshift=-0.5em},
  point meta min=\cmin,point meta max=\cmax,
  label style={font=\scriptsize},
  tick label style={font=\scriptsize},
  legend style={font=\scriptsize\selectfont},
]
\addplot [forget plot] graphics [xmin=\xmingb,xmax=\xmaxgb,ymin=\zmingb,ymax=\zmaxgb] {{\modelfile}.png};
\end{axis}
\end{tikzpicture}%
    \caption{Reconstruction for model $f=\fdiskbb$ 
             (PSNR \num{15.33}).}
  \end{subfigure}\hspace*{1.50em}
  \renewcommand{\modelfile}{images/fwi/reconstruction/fwi_R4.50-c1.07_unifreq0.7Hz-f}
  \pgfmathsetmacro{\cmin} {-5} \pgfmathsetmacro{\cmax} {5}  
  \begin{subfigure}[t]{.375\textwidth} \centering
    \begin{tikzpicture}

\pgfmathsetmacro{\xmingb} {-50} \pgfmathsetmacro{\xmaxgb}{50}
\pgfmathsetmacro{\zmingb} {-50} \pgfmathsetmacro{\zmaxgb}{50}

\begin{axis}[
  width=\modelwidth, height=\modelheight,
  axis on top, separate axis lines,
  xmin=\xminloc, xmax=\xmaxloc, 
  ymin=\zminloc, ymax=\zmaxloc, 
  x label style={xshift=-0.0cm, yshift= 0.00cm}, 
  y label style={xshift= 0.0cm, yshift=-0.30cm},
  colormap/jet,colorbar,colorbar style={title={\scriptsize $f$},title style={yshift=-2mm, xshift=0mm},
  width=.25cm, xshift=-0.5em},
  point meta min=\cmin,point meta max=\cmax,
  label style={font=\scriptsize},
  tick label style={font=\scriptsize},
  legend style={font=\scriptsize\selectfont},
]
\addplot [forget plot] graphics [xmin=\xmingb,xmax=\xmaxgb,ymin=\zmingb,ymax=\zmaxgb] {{\modelfile}.png};
\end{axis}
\end{tikzpicture}%
    \caption{Reconstruction for model $f=\fdiskbc$ 
             (PSNR \num{15.03}).}
  \end{subfigure}\hfill

  \caption{Reconstruction using frequency $\omega/(2\pi)=\num{0.7}$. 
           The initial model consists in a constant wave speed $c_0=1$.
           The data consist of $n_A=\num{40}$ different angles of incidence from 
           \SI{0}{\degree} to \SI{351}{\degree}.}
  \label{fig:fwi:reconstruction_unifreq0.7Hz}
  \end{figure}
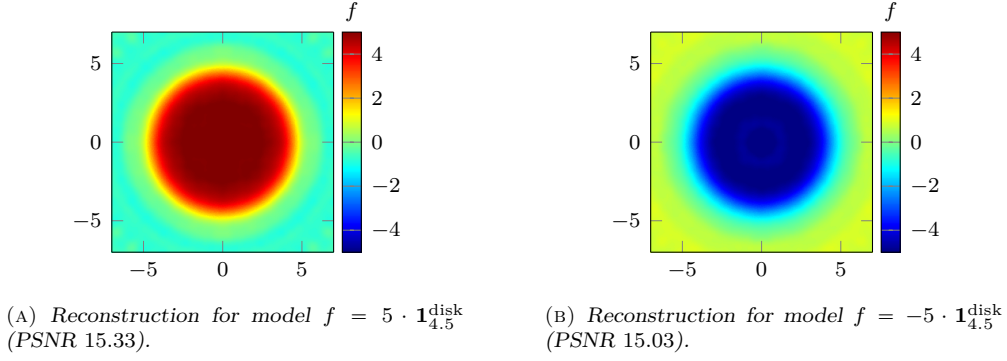
\end{remark}

\subsubsection{Reconstruction using Born and Rytov approximations}

For the reconstruction with \autoref{algo:infft}, which relies on the Born or Rytov approximation,
we use the same data as in the above experiment.
We use a grid with $N=240$ and $\rs = {N}/({8\sqrt2})\approx10$. 
The numerical results indicate that $\rs$ should not be smaller than $\lM$.
Since we have $k_1^2\le k_0^2$ and the distance between two grid points of $k_1$ is ${\pi}/{\lM}$, only around ${2k_0 \lM}/{\pi}\approx40$ of them contribute to the data of the inverse NDFT.

In the following reconstructions, we use a fixed number of 20 iteration steps in the conjugate gradient method.
Initially, we use the frequency $\omega/(2\pi)=1$ of the incident wave, therefore $k_0=2\pi$.
Reconstructions of the circular model $f=\mathbf 1^{\mathrm{disk}}_a$ 
are shown in \autoref{fig:f0-full-waveform-data}.
We note that all reconstructions are reasonably good, 
where the Rytov reconstruction looks slightly better inside the object.

\setlength{\modelwidth}{5.00cm}
\setlength{\modelheight}{\modelwidth}
\setlength{\plotwidth}{10cm}
\setlength{\plotheight}{3cm}

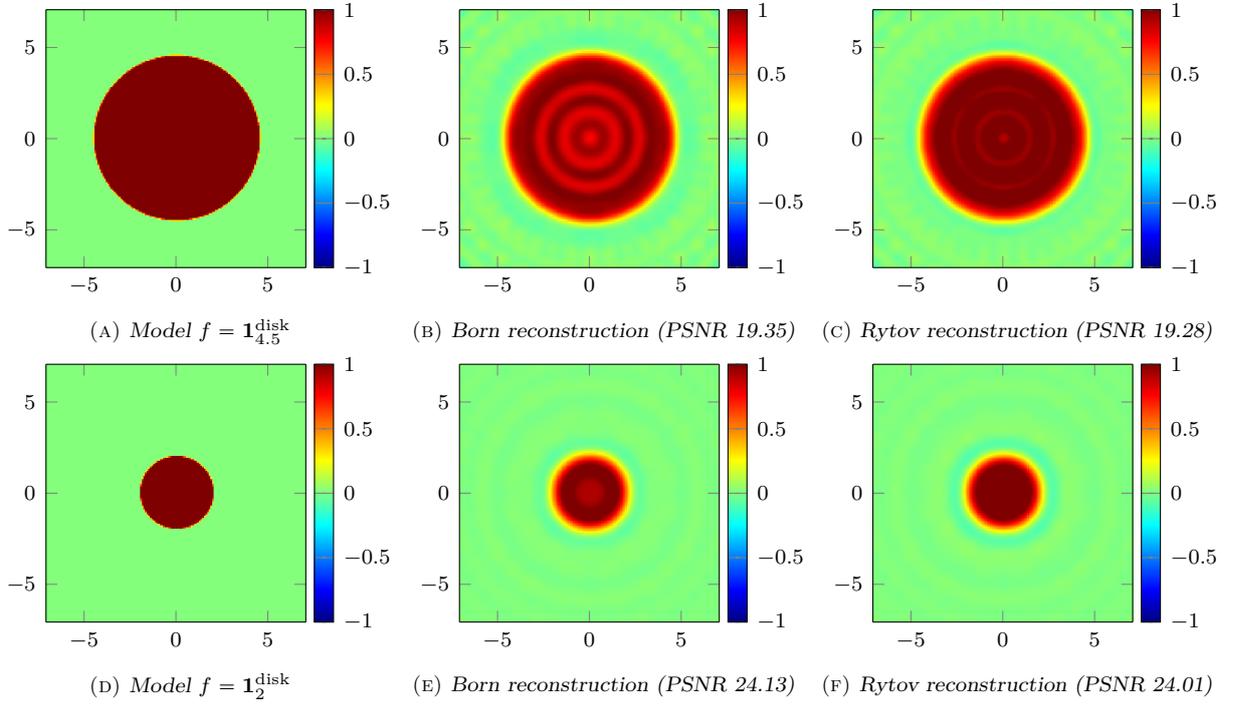
\begin{figure}[ht!]
  \pgfmathsetmacro{\xminloc}{-7.0711}\pgfmathsetmacro{\xmaxloc}{ 7.0711}
  \pgfmathsetmacro{\zminloc}{-7.0711}\pgfmathsetmacro{\zmaxloc}{ 7.0711}
  \pgfmathsetmacro{\xminpic}{-7.0711}\pgfmathsetmacro{\xmaxpic}{ 7.0711}
  \pgfmathsetmacro{\zminpic}{-7.0711}\pgfmathsetmacro{\zmaxpic}{ 7.0711}
  \pgfmathsetmacro{\cmin} {-1} \pgfmathsetmacro{\cmax} {1}
  \begin{subfigure}[t]{.32\textwidth} \centering
    \renewcommand{\modelfile}{images/disk/disk-1000mHz-r4-f1}
    \begin{tikzpicture}
  \begin{axis}[
    width=\modelwidth, height=\modelheight,
    axis on top, separate axis lines,
    xmin=\xminpic, xmax=\xmaxpic, 
    ymin=\zminpic, ymax=\zmaxpic, 
    x label style={xshift=-0.0cm, yshift= 0.00cm}, 
    y label style={xshift= 0.0cm, yshift=-0.30cm},
    colormap/jet,colorbar,colorbar style={
      width=.25cm, xshift=-.5em},
    point meta min=\cmin,point meta max=\cmax,
    label style={font=\scriptsize},
    tick label style={font=\scriptsize},
    legend style={font=\scriptsize\selectfont},
    ]
    \addplot [forget plot] graphics [xmin=\xminloc,xmax=\xmaxloc,ymin=\zminloc,ymax=\zmaxloc] {{\modelfile}.png};
  \end{axis}
\end{tikzpicture}%
    \caption{Model $f=\mathbf 1^{\mathrm{disk}}_{4.5}$}
  \end{subfigure}\hfill
  \begin{subfigure}[t]{.32\textwidth} \centering
    \renewcommand{\modelfile}{images/disk/disk-1000mHz-r4-f1-rec}
    \begin{tikzpicture}
  \begin{axis}[
    width=\modelwidth, height=\modelheight,
    axis on top, separate axis lines,
    xmin=\xminpic, xmax=\xmaxpic, 
    ymin=\zminpic, ymax=\zmaxpic, 
    x label style={xshift=-0.0cm, yshift= 0.00cm}, 
    y label style={xshift= 0.0cm, yshift=-0.30cm},
    colormap/jet,colorbar,colorbar style={
      width=.25cm, xshift=-.5em},
    point meta min=\cmin,point meta max=\cmax,
    label style={font=\scriptsize},
    tick label style={font=\scriptsize},
    legend style={font=\scriptsize\selectfont},
    ]
    \addplot [forget plot] graphics [xmin=\xminloc,xmax=\xmaxloc,ymin=\zminloc,ymax=\zmaxloc] {{\modelfile}.png};
  \end{axis}
\end{tikzpicture}%
    \caption{Born reconstruction (PSNR 19.35)}
  \end{subfigure}\hfill
  \begin{subfigure}[t]{.32\textwidth} \centering
    \renewcommand{\modelfile}{images/disk/disk-1000mHz-r4-f1-rec-Rytov}
    \begin{tikzpicture}
  \begin{axis}[
    width=\modelwidth, height=\modelheight,
    axis on top, separate axis lines,
    xmin=\xminpic, xmax=\xmaxpic, 
    ymin=\zminpic, ymax=\zmaxpic, 
    x label style={xshift=-0.0cm, yshift= 0.00cm}, 
    y label style={xshift= 0.0cm, yshift=-0.30cm},
    colormap/jet,colorbar,colorbar style={
      width=.25cm, xshift=-.5em},
    point meta min=\cmin,point meta max=\cmax,
    label style={font=\scriptsize},
    tick label style={font=\scriptsize},
    legend style={font=\scriptsize\selectfont},
    ]
    \addplot [forget plot] graphics [xmin=\xminloc,xmax=\xmaxloc,ymin=\zminloc,ymax=\zmaxloc] {{\modelfile}.png};
  \end{axis}
\end{tikzpicture}%
    \caption{Rytov reconstruction (PSNR 19.28)}
  \end{subfigure}
  \begin{subfigure}[t]{.32\textwidth} \centering
    \renewcommand{\modelfile}{images/disk/disk-1000mHz-r2-f1}
    \begin{tikzpicture}
  \begin{axis}[
    width=\modelwidth, height=\modelheight,
    axis on top, separate axis lines,
    xmin=\xminpic, xmax=\xmaxpic, 
    ymin=\zminpic, ymax=\zmaxpic, 
    x label style={xshift=-0.0cm, yshift= 0.00cm}, 
    y label style={xshift= 0.0cm, yshift=-0.30cm},
    colormap/jet,colorbar,colorbar style={
      width=.25cm, xshift=-.5em},
    point meta min=\cmin,point meta max=\cmax,
    label style={font=\scriptsize},
    tick label style={font=\scriptsize},
    legend style={font=\scriptsize\selectfont},
    ]
    \addplot [forget plot] graphics [xmin=\xminloc,xmax=\xmaxloc,ymin=\zminloc,ymax=\zmaxloc] {{\modelfile}.png};
  \end{axis}
\end{tikzpicture}%
  \caption{Model $f=\mathbf 1^{\mathrm{disk}}_{2}$}
  \end{subfigure}\hfill
  \begin{subfigure}[t]{.32\textwidth} \centering
    \renewcommand{\modelfile}{images/disk/disk-1000mHz-r2-f1-rec}
    \begin{tikzpicture}
  \begin{axis}[
    width=\modelwidth, height=\modelheight,
    axis on top, separate axis lines,
    xmin=\xminpic, xmax=\xmaxpic, 
    ymin=\zminpic, ymax=\zmaxpic, 
    x label style={xshift=-0.0cm, yshift= 0.00cm}, 
    y label style={xshift= 0.0cm, yshift=-0.30cm},
    colormap/jet,colorbar,colorbar style={
      width=.25cm, xshift=-.5em},
    point meta min=\cmin,point meta max=\cmax,
    label style={font=\scriptsize},
    tick label style={font=\scriptsize},
    legend style={font=\scriptsize\selectfont},
    ]
    \addplot [forget plot] graphics [xmin=\xminloc,xmax=\xmaxloc,ymin=\zminloc,ymax=\zmaxloc] {{\modelfile}.png};
  \end{axis}
\end{tikzpicture}%
  \caption{Born reconstruction (PSNR 24.13)}
  \end{subfigure}\hfill
  \begin{subfigure}[t]{.32\textwidth} \centering
    \renewcommand{\modelfile}{images/disk/disk-1000mHz-r2-f1-rec-Rytov}
    \begin{tikzpicture}
  \begin{axis}[
    width=\modelwidth, height=\modelheight,
    axis on top, separate axis lines,
    xmin=\xminpic, xmax=\xmaxpic, 
    ymin=\zminpic, ymax=\zmaxpic, 
    x label style={xshift=-0.0cm, yshift= 0.00cm}, 
    y label style={xshift= 0.0cm, yshift=-0.30cm},
    colormap/jet,colorbar,colorbar style={
      width=.25cm, xshift=-.5em},
    point meta min=\cmin,point meta max=\cmax,
    label style={font=\scriptsize},
    tick label style={font=\scriptsize},
    legend style={font=\scriptsize\selectfont},
    ]
    \addplot [forget plot] graphics [xmin=\xminloc,xmax=\xmaxloc,ymin=\zminloc,ymax=\zmaxloc] {{\modelfile}.png};
  \end{axis}
\end{tikzpicture}%
  \caption{Rytov reconstruction (PSNR 24.01)}
  \end{subfigure}
  \caption{Reconstructions with with the Born and Rytov approximation, where the data $u(\cdot,\rM)$ is generated with the line source model.
  The incident field has the frequency $\omega/(2\pi)=1$.
  Visible is only the cut out center, where we compute the PSNR.
    \label{fig:f0-full-waveform-data}}
\end{figure}

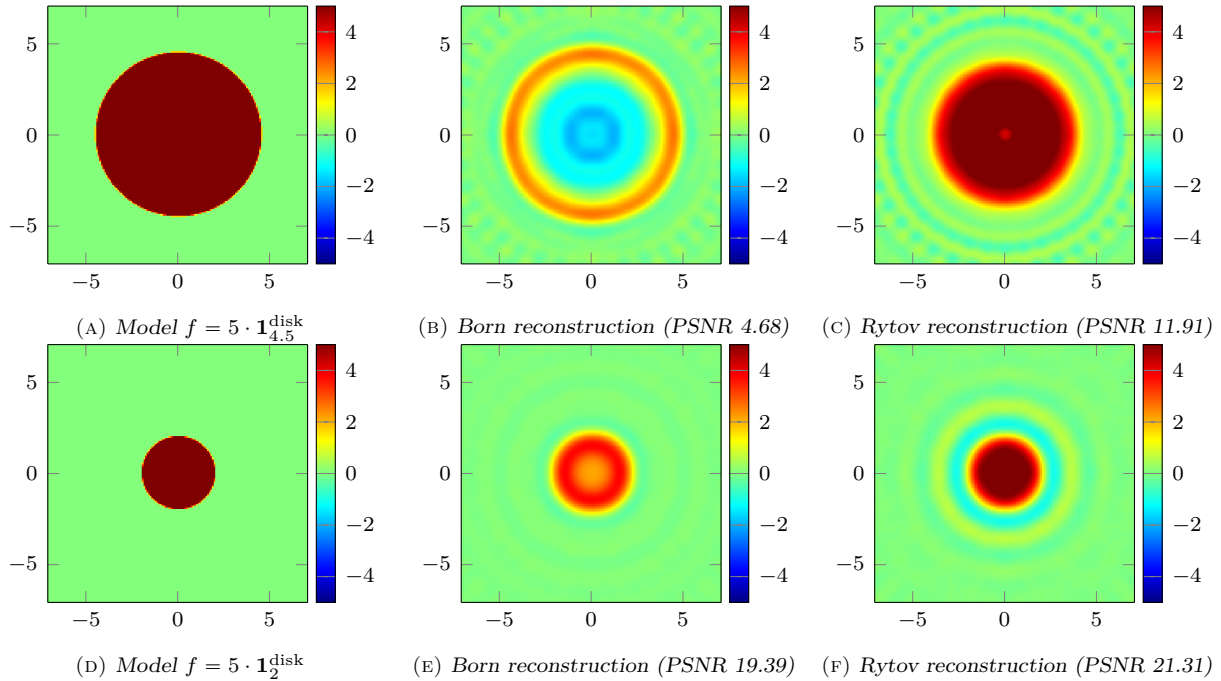
\begin{figure}[ht!]
\pgfmathsetmacro{\xminloc}{-7.0711}\pgfmathsetmacro{\xmaxloc}{ 7.0711}
\pgfmathsetmacro{\zminloc}{-7.0711}\pgfmathsetmacro{\zmaxloc}{ 7.0711}
\pgfmathsetmacro{\xminpic}{-7.0711}\pgfmathsetmacro{\xmaxpic}{ 7.0711}
\pgfmathsetmacro{\zminpic}{-7.0711}\pgfmathsetmacro{\zmaxpic}{ 7.0711}
\pgfmathsetmacro{\cmin} {-5} \pgfmathsetmacro{\cmax} {5}
\begin{subfigure}[t]{.32\textwidth} \centering
  \renewcommand{\modelfile}{images/disk/disk-1000mHz-r4-f5}
  \begin{tikzpicture}
  \begin{axis}[
    width=\modelwidth, height=\modelheight,
    axis on top, separate axis lines,
    xmin=\xminpic, xmax=\xmaxpic, 
    ymin=\zminpic, ymax=\zmaxpic, 
    x label style={xshift=-0.0cm, yshift= 0.00cm}, 
    y label style={xshift= 0.0cm, yshift=-0.30cm},
    colormap/jet,colorbar,colorbar style={
      width=.25cm, xshift=-.5em},
    point meta min=\cmin,point meta max=\cmax,
    label style={font=\scriptsize},
    tick label style={font=\scriptsize},
    legend style={font=\scriptsize\selectfont},
    ]
    \addplot [forget plot] graphics [xmin=\xminloc,xmax=\xmaxloc,ymin=\zminloc,ymax=\zmaxloc] {{\modelfile}.png};
  \end{axis}
\end{tikzpicture}%
  \caption{Model $f=5\cdot\mathbf 1^{\mathrm{disk}}_{4.5}$}
\end{subfigure}\hfill
\begin{subfigure}[t]{.32\textwidth} \centering
  \renewcommand{\modelfile}{images/disk/disk-1000mHz-r4-f5-rec}
  \begin{tikzpicture}
  \begin{axis}[
    width=\modelwidth, height=\modelheight,
    axis on top, separate axis lines,
    xmin=\xminpic, xmax=\xmaxpic, 
    ymin=\zminpic, ymax=\zmaxpic, 
    x label style={xshift=-0.0cm, yshift= 0.00cm}, 
    y label style={xshift= 0.0cm, yshift=-0.30cm},
    colormap/jet,colorbar,colorbar style={
      width=.25cm, xshift=-.5em},
    point meta min=\cmin,point meta max=\cmax,
    label style={font=\scriptsize},
    tick label style={font=\scriptsize},
    legend style={font=\scriptsize\selectfont},
    ]
    \addplot [forget plot] graphics [xmin=\xminloc,xmax=\xmaxloc,ymin=\zminloc,ymax=\zmaxloc] {{\modelfile}.png};
  \end{axis}
\end{tikzpicture}%
  \caption{Born reconstruction (PSNR 4.68)}
\end{subfigure}\hfill
\begin{subfigure}[t]{.32\textwidth} \centering
  \renewcommand{\modelfile}{images/disk/disk-1000mHz-r4-f5-rec-Rytov}
  \begin{tikzpicture}
  \begin{axis}[
    width=\modelwidth, height=\modelheight,
    axis on top, separate axis lines,
    xmin=\xminpic, xmax=\xmaxpic, 
    ymin=\zminpic, ymax=\zmaxpic, 
    x label style={xshift=-0.0cm, yshift= 0.00cm}, 
    y label style={xshift= 0.0cm, yshift=-0.30cm},
    colormap/jet,colorbar,colorbar style={
      width=.25cm, xshift=-.5em},
    point meta min=\cmin,point meta max=\cmax,
    label style={font=\scriptsize},
    tick label style={font=\scriptsize},
    legend style={font=\scriptsize\selectfont},
    ]
    \addplot [forget plot] graphics [xmin=\xminloc,xmax=\xmaxloc,ymin=\zminloc,ymax=\zmaxloc] {{\modelfile}.png};
  \end{axis}
\end{tikzpicture}%
  \caption{Rytov reconstruction (PSNR 11.91)}
\end{subfigure}
\begin{subfigure}[t]{.32\textwidth} \centering
  \renewcommand{\modelfile}{images/disk/disk-1000mHz-r2-f5}
  \begin{tikzpicture}
  \begin{axis}[
    width=\modelwidth, height=\modelheight,
    axis on top, separate axis lines,
    xmin=\xminpic, xmax=\xmaxpic, 
    ymin=\zminpic, ymax=\zmaxpic, 
    x label style={xshift=-0.0cm, yshift= 0.00cm}, 
    y label style={xshift= 0.0cm, yshift=-0.30cm},
    colormap/jet,colorbar,colorbar style={
      width=.25cm, xshift=-.5em},
    point meta min=\cmin,point meta max=\cmax,
    label style={font=\scriptsize},
    tick label style={font=\scriptsize},
    legend style={font=\scriptsize\selectfont},
    ]
    \addplot [forget plot] graphics [xmin=\xminloc,xmax=\xmaxloc,ymin=\zminloc,ymax=\zmaxloc] {{\modelfile}.png};
  \end{axis}
\end{tikzpicture}%
  \caption{Model $f=5\cdot\mathbf 1^{\mathrm{disk}}_{2}$}
\end{subfigure}\hfill
\begin{subfigure}[t]{.32\textwidth} \centering
  \renewcommand{\modelfile}{images/disk/disk-1000mHz-r2-f5-rec}
  \begin{tikzpicture}
  \begin{axis}[
    width=\modelwidth, height=\modelheight,
    axis on top, separate axis lines,
    xmin=\xminpic, xmax=\xmaxpic, 
    ymin=\zminpic, ymax=\zmaxpic, 
    x label style={xshift=-0.0cm, yshift= 0.00cm}, 
    y label style={xshift= 0.0cm, yshift=-0.30cm},
    colormap/jet,colorbar,colorbar style={
      width=.25cm, xshift=-.5em},
    point meta min=\cmin,point meta max=\cmax,
    label style={font=\scriptsize},
    tick label style={font=\scriptsize},
    legend style={font=\scriptsize\selectfont},
    ]
    \addplot [forget plot] graphics [xmin=\xminloc,xmax=\xmaxloc,ymin=\zminloc,ymax=\zmaxloc] {{\modelfile}.png};
  \end{axis}
\end{tikzpicture}%
  \caption{Born reconstruction (PSNR 19.39)}
\end{subfigure}\hfill
\begin{subfigure}[t]{.32\textwidth} \centering
  \renewcommand{\modelfile}{images/disk/disk-1000mHz-r2-f5-rec-Rytov}
  \begin{tikzpicture}
  \begin{axis}[
    width=\modelwidth, height=\modelheight,
    axis on top, separate axis lines,
    xmin=\xminpic, xmax=\xmaxpic, 
    ymin=\zminpic, ymax=\zmaxpic, 
    x label style={xshift=-0.0cm, yshift= 0.00cm}, 
    y label style={xshift= 0.0cm, yshift=-0.30cm},
    colormap/jet,colorbar,colorbar style={
      width=.25cm, xshift=-.5em},
    point meta min=\cmin,point meta max=\cmax,
    label style={font=\scriptsize},
    tick label style={font=\scriptsize},
    legend style={font=\scriptsize\selectfont},
    ]
    \addplot [forget plot] graphics [xmin=\xminloc,xmax=\xmaxloc,ymin=\zminloc,ymax=\zmaxloc] {{\modelfile}.png};
  \end{axis}
\end{tikzpicture}%
  \caption{Rytov reconstruction (PSNR 21.31)}
\end{subfigure}
\caption{Same setting as in \autoref{fig:f0-full-waveform-data}, but with a higher amplitude of $5$.
\label{fig:f5-full-waveform-data}}
\end{figure}

For a higher amplitude of the model function $f$, the limitations of the linear models become apparent.
In \autoref{fig:f5-full-waveform-data}, we see that the Born reconstruction of the larger object fails,
and for the smaller object only the Rytov approximation yields a good reconstruction,
which is consistent with \autoref{rem:Born-Rytov}.
With the Born approximation, we recognize the object's shape but not its amplitude, which is consistent with the observations in \cite{MulSchGuc15}.
However, as we see in \autoref{fig:fwi:reconstruction_unifreq}, even the FWI reconstruction makes a considerable error in the object's interior, and we cannot expect the linear models to be better.

We see that the FWI and the Born/Rytov reconstructions contain different kinds of artifacts.
Therefore, a comparison of the visual image quality perception does not necessarily yield the same conclusions as for the computed PSNR values.
Furthermore, the size of the object has a considerable effect on the PSNR,
e.g., the images in \autoref{fig:f5-full-waveform-data} (c) and (f) show a comparable visual quality, but the latter's PSNR is considerably better because of the lower error in the background farther away from the object,
see also \cite{HuyGha2010} for a study on the PSNR.

\begin{figure}[ht!]
  \pgfmathsetmacro{\xminloc}{-7.0711}\pgfmathsetmacro{\xmaxloc}{ 7.0711}
  \pgfmathsetmacro{\zminloc}{-7.0711}\pgfmathsetmacro{\zmaxloc}{ 7.0711}
  \pgfmathsetmacro{\xminpic}{-7.0711}\pgfmathsetmacro{\xmaxpic}{ 7.0711}
  \pgfmathsetmacro{\zminpic}{-7.0711}\pgfmathsetmacro{\zmaxpic}{ 7.0711}
  \pgfmathsetmacro{\cmin} {-1} \pgfmathsetmacro{\cmax} {1}
  \centering
  \begin{subfigure}[t]{.375\textwidth}\centering
    \renewcommand{\modelfile}{images/disk/disk-700mHz-r4-f1-rec}
    \begin{tikzpicture}
  \begin{axis}[
    width=\modelwidth, height=\modelheight,
    axis on top, separate axis lines,
    xmin=\xminpic, xmax=\xmaxpic, 
    ymin=\zminpic, ymax=\zmaxpic, 
    x label style={xshift=-0.0cm, yshift= 0.00cm}, 
    y label style={xshift= 0.0cm, yshift=-0.30cm},
    colormap/jet,colorbar,colorbar style={
      width=.25cm, xshift=-.5em},
    point meta min=\cmin,point meta max=\cmax,
    label style={font=\scriptsize},
    tick label style={font=\scriptsize},
    legend style={font=\scriptsize\selectfont},
    ]
    \addplot [forget plot] graphics [xmin=\xminloc,xmax=\xmaxloc,ymin=\zminloc,ymax=\zmaxloc] {{\modelfile}.png};
  \end{axis}
\end{tikzpicture}%
    \caption{Born reconstruction (PSNR 18.05)}
  \end{subfigure}\qquad
  \begin{subfigure}[t]{.375\textwidth}\centering
    \renewcommand{\modelfile}{images/disk/disk-700mHz-r4-f1-rec-Rytov}
    \begin{tikzpicture}
  \begin{axis}[
    width=\modelwidth, height=\modelheight,
    axis on top, separate axis lines,
    xmin=\xminpic, xmax=\xmaxpic, 
    ymin=\zminpic, ymax=\zmaxpic, 
    x label style={xshift=-0.0cm, yshift= 0.00cm}, 
    y label style={xshift= 0.0cm, yshift=-0.30cm},
    colormap/jet,colorbar,colorbar style={
      width=.25cm, xshift=-.5em},
    point meta min=\cmin,point meta max=\cmax,
    label style={font=\scriptsize},
    tick label style={font=\scriptsize},
    legend style={font=\scriptsize\selectfont},
    ]
    \addplot [forget plot] graphics [xmin=\xminloc,xmax=\xmaxloc,ymin=\zminloc,ymax=\zmaxloc] {{\modelfile}.png};
  \end{axis}
\end{tikzpicture}%
    \caption{Rytov reconstruction (PSNR 16.52)}
  \end{subfigure}
  \caption{
    Reconstructions of $\mathbf 1^{\mathrm{disk}}_{4.5}$, where the incident field has the frequency $\omega/(2\pi)=0.7$ instead of $1$.
    The rest of the setting is from \autoref{fig:f0-full-waveform-data}.
    \label{fig:f0-full-waveform-data-07}}
\end{figure}
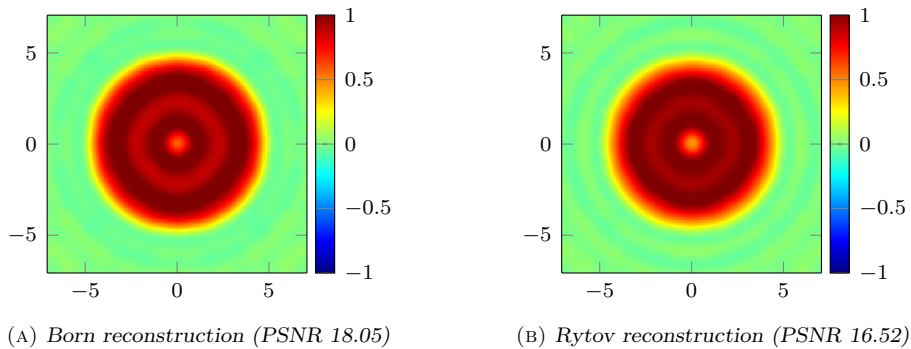

In \autoref{fig:f0-full-waveform-data-07}, we use the same setup as before, but with the 
frequency $\omega/(2\pi)=0.7$ instead of $\omega/(2\pi)=1$ and thus the wave number $k_0=\omega$.
Apparently, the reconstruction becomes worse with lower frequency, because it provides a smaller k-space coverage.

\subsection{Reconstruction of embedded shapes: Phantom 1}

We consider a more challenging scenario with 
shapes embedded in the background medium.
In \autoref{fig:fwi:heart-model_f05_data_1}, we picture
the perturbation $f$ consisting of a disk and heart 
included in an ellipse, 
with $f$ varying from \num{0} to \num{0.5}.
The computational domain corresponds to 
$[-20,20]\times[-20,20]$, with line-sources positioned at
a distance $R=\num{10}$ and receivers in $\rM=\num{6}$ to 
capture the data.
The data are generated using $n_A=\num{100}$ incidence angles $\alpha$ equispaced on $[0,2\pi]$,
following the steps described in \autoref{subsection:modeling:line-src}.
This is illustrated in \autoref{fig:fwi:heart-model_f05_data}.

\setlength{\modelwidth}{4.50cm}
\setlength{\modelheight}{\modelwidth}
\setlength{\plotwidth}{10cm}
\setlength{\plotheight}{3cm}
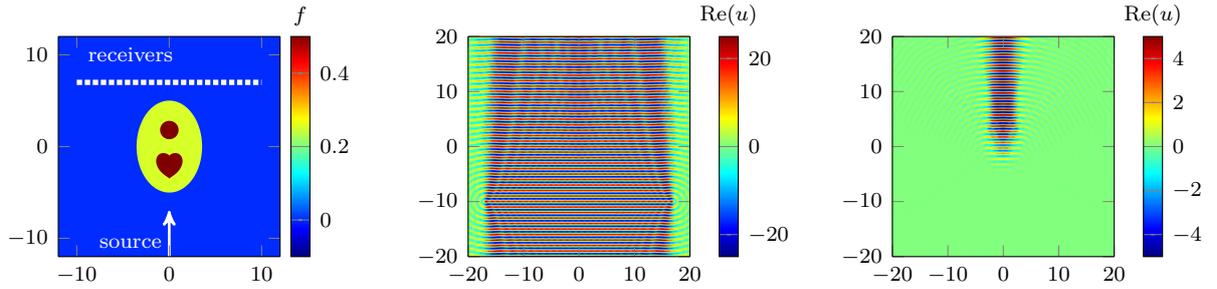
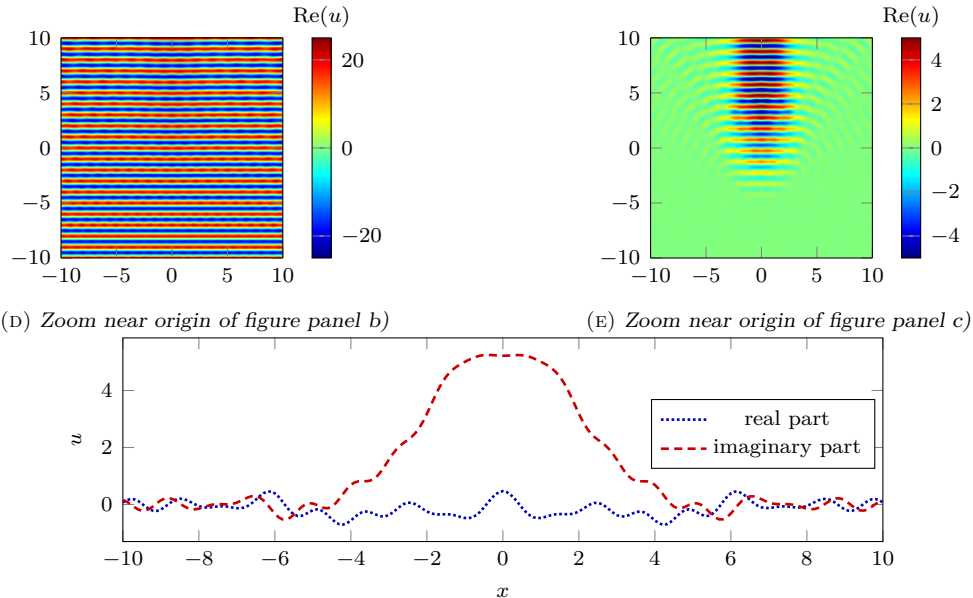
\begin{figure}[ht!]\centering
  \pgfmathsetmacro{\xminloc}{-12}\pgfmathsetmacro{\xmaxloc}{12}
  \pgfmathsetmacro{\zminloc}{-12}\pgfmathsetmacro{\zmaxloc}{12}
  \renewcommand{\modelfile}{images/fwi/model_heart/models/cp1_f05_f_scale-01-05}
  \pgfmathsetmacro{\cmin} {-0.1} \pgfmathsetmacro{\cmax} {0.5}
  \begin{subfigure}[t]{.305\textwidth}\centering
    \begin{tikzpicture}

\pgfmathsetmacro{\xmingb} {-20} \pgfmathsetmacro{\xmaxgb}{20}
\pgfmathsetmacro{\zmingb} {-20} \pgfmathsetmacro{\zmaxgb}{20}

\begin{axis}[
  width=\modelwidth, height=\modelheight,
  axis on top, separate axis lines,
  xmin=\xminloc, xmax=\xmaxloc, 
  ymin=\zminloc, ymax=\zmaxloc, 
  x label style={xshift=-0.0cm, yshift= 0.00cm}, 
  y label style={xshift= 0.0cm, yshift=-0.30cm},
  colormap/jet,colorbar,colorbar style={title={\scriptsize $f$},
  title style={yshift=-2mm, xshift=0mm},
  width=.25cm, xshift=-0.5em},
  point meta min=\cmin,point meta max=\cmax,
  label style={font=\scriptsize},
  tick label style={font=\scriptsize},
  legend style={font=\scriptsize\selectfont},
]
\addplot [forget plot] graphics [xmin=\xmingb,xmax=\xmaxgb,ymin=\zmingb,ymax=\zmaxgb] {{\modelfile}.png};

\draw[line width=2, white, densely dotted] (-10,7) -- (10,7); 
\draw[line width=1, white,->] (0,-15) -- (0,-7); 

\node[xshift=5.0em,yshift= 3.0em,white ] (cs) {\scriptsize source}; 
\node[xshift=5.0em,yshift= 9.5em,white] (cs) {\scriptsize receivers}; 

\end{axis}
\end{tikzpicture}%
    \caption{Perturbation model at frequency $\omega/(2\pi)=1$, 
             the wave speed is equal to \num{1} in the background.
             The positions of the source and the receivers recording
             transmission data are pictured in white.
             }
    \label{fig:fwi:heart-model_f05_data_1}
  \end{subfigure}\hfill
  \pgfmathsetmacro{\xminloc}{-20}\pgfmathsetmacro{\xmaxloc}{20}
  \pgfmathsetmacro{\zminloc}{-20}\pgfmathsetmacro{\zmaxloc}{20}
  \pgfmathsetmacro{\cmin} {-25} \pgfmathsetmacro{\cmax} { 25}
  \renewcommand{\modelfile}{images/fwi/model_heart/modeling/wave2d_domain40x40_src34_cp-true_scale25}
  \begin{subfigure}[t]{.32\textwidth}\centering
    \begin{tikzpicture}

\pgfmathsetmacro{\xmingb} {-20} \pgfmathsetmacro{\xmaxgb}{20}
\pgfmathsetmacro{\zmingb} {-20} \pgfmathsetmacro{\zmaxgb}{20}

\begin{axis}[
  width=\modelwidth, height=\modelheight,
  axis on top, separate axis lines,
  xmin=\xminloc, xmax=\xmaxloc, 
  ymin=\zminloc, ymax=\zmaxloc, 
  x label style={xshift=-0.0cm, yshift= 0.00cm}, 
  y label style={xshift= 0.0cm, yshift=-0.30cm},
  colormap/jet,colorbar,colorbar style={title={\scriptsize $\mathrm{Re}(u)$},
  title style={yshift=-2mm, xshift=0mm},
  width=.25cm, xshift=-0.5em},
  point meta min=\cmin,point meta max=\cmax,
  label style={font=\scriptsize},
  tick label style={font=\scriptsize},
  legend style={font=\scriptsize\selectfont},
]
\addplot [forget plot] graphics [xmin=\xmingb,xmax=\xmaxgb,ymin=\zmingb,ymax=\zmaxgb] {{\modelfile}.png};
\end{axis}
\end{tikzpicture}%
    \caption{Real part of the global solution to \autoref{eq:fwi:helmholtz-full}
             at frequency $\omega/(2\pi)=\num{1}$, the source is discretized by 
             $N_{\mathrm{sim}}=\num{1361}$ simultaneous excitations at fixed height $x_2=-10$.}
  \end{subfigure} \hfill
  \pgfmathsetmacro{\cmin} {-5} \pgfmathsetmacro{\cmax} { 5}
  \renewcommand{\modelfile}{images/fwi/model_heart/modeling/wave2d_domain40x40_src34_scat_scale5}
  \begin{subfigure}[t]{.31\textwidth}\centering
    \begin{tikzpicture}

\pgfmathsetmacro{\xmingb} {-20} \pgfmathsetmacro{\xmaxgb}{20}
\pgfmathsetmacro{\zmingb} {-20} \pgfmathsetmacro{\zmaxgb}{20}

\begin{axis}[
  width=\modelwidth, height=\modelheight,
  axis on top, separate axis lines,
  xmin=\xminloc, xmax=\xmaxloc, 
  ymin=\zminloc, ymax=\zmaxloc, 
  x label style={xshift=-0.0cm, yshift= 0.00cm}, 
  y label style={xshift= 0.0cm, yshift=-0.30cm},
  colormap/jet,colorbar,colorbar style={title={\scriptsize $\mathrm{Re}(u)$},
  title style={yshift=-2mm, xshift=0mm},
  width=.25cm, xshift=-0.5em},
  point meta min=\cmin,point meta max=\cmax,
  label style={font=\scriptsize},
  tick label style={font=\scriptsize},
  legend style={font=\scriptsize\selectfont},
]
\addplot [forget plot] graphics [xmin=\xmingb,xmax=\xmaxgb,ymin=\zmingb,ymax=\zmaxgb] {{\modelfile}.png};
\end{axis}
\end{tikzpicture}%
    \caption{Real part of the scattering solution at frequency $\omega/(2\pi)=\num{1}$.}
  \end{subfigure} \\

  \pgfmathsetmacro{\xminloc}{-10}\pgfmathsetmacro{\xmaxloc}{10}
  \pgfmathsetmacro{\zminloc}{-10}\pgfmathsetmacro{\zmaxloc}{10}
  \pgfmathsetmacro{\cmin} {-25} \pgfmathsetmacro{\cmax} { 25}
  \renewcommand{\modelfile}{images/fwi/model_heart/modeling/wave2d_domain40x40_src34_cp-true_scale25}
  \begin{subfigure}[t]{.40\textwidth}\centering
    \begin{tikzpicture}

\pgfmathsetmacro{\xmingb} {-20} \pgfmathsetmacro{\xmaxgb}{20}
\pgfmathsetmacro{\zmingb} {-20} \pgfmathsetmacro{\zmaxgb}{20}

\begin{axis}[
  width=\modelwidth, height=\modelheight,
  axis on top, separate axis lines,
  xmin=\xminloc, xmax=\xmaxloc, 
  ymin=\zminloc, ymax=\zmaxloc, 
  x label style={xshift=-0.0cm, yshift= 0.00cm}, 
  y label style={xshift= 0.0cm, yshift=-0.30cm},
  colormap/jet,colorbar,colorbar style={title={\scriptsize $\mathrm{Re}(u)$},
  title style={yshift=-2mm, xshift=0mm},
  width=.25cm, xshift=-0.5em},
  point meta min=\cmin,point meta max=\cmax,
  label style={font=\scriptsize},
  tick label style={font=\scriptsize},
  legend style={font=\scriptsize\selectfont},
]
\addplot [forget plot] graphics [xmin=\xmingb,xmax=\xmaxgb,ymin=\zmingb,ymax=\zmaxgb] {{\modelfile}.png};
\end{axis}
\end{tikzpicture}%
    \caption{Zoom near origin of figure panel b)
            }
  \end{subfigure} \qquad
  \pgfmathsetmacro{\cmin} {-5} \pgfmathsetmacro{\cmax} { 5}
  \renewcommand{\modelfile}{images/fwi/model_heart/modeling/wave2d_domain40x40_src34_scat_scale5}
  \begin{subfigure}[t]{.40\textwidth}\centering
    \begin{tikzpicture}

\pgfmathsetmacro{\xmingb} {-20} \pgfmathsetmacro{\xmaxgb}{20}
\pgfmathsetmacro{\zmingb} {-20} \pgfmathsetmacro{\zmaxgb}{20}

\begin{axis}[
  width=\modelwidth, height=\modelheight,
  axis on top, separate axis lines,
  xmin=\xminloc, xmax=\xmaxloc, 
  ymin=\zminloc, ymax=\zmaxloc, 
  x label style={xshift=-0.0cm, yshift= 0.00cm}, 
  y label style={xshift= 0.0cm, yshift=-0.30cm},
  colormap/jet,colorbar,colorbar style={title={\scriptsize $\mathrm{Re}(u)$},
  title style={yshift=-2mm, xshift=0mm},
  width=.25cm, xshift=-0.5em},
  point meta min=\cmin,point meta max=\cmax,
  label style={font=\scriptsize},
  tick label style={font=\scriptsize},
  legend style={font=\scriptsize\selectfont},
]
\addplot [forget plot] graphics [xmin=\xmingb,xmax=\xmaxgb,ymin=\zmingb,ymax=\zmaxgb] {{\modelfile}.png};
\end{axis}
\end{tikzpicture}%
    \caption{Zoom near origin of figure panel c)}
  \end{subfigure}

  \renewcommand{\datafile}{images/fwi/model_heart/modeling/data_scattscaled_domain40x40_src34.txt}
  \begin{subfigure}[t]{\linewidth} \centering
\begin{tikzpicture}
\begin{axis}[
             enlargelimits=false, 
             ylabel={\scriptsize $u$},
             xlabel={$x$},
             enlarge y limits=true, 
             enlarge x limits=false,
             yminorticks=true,
             height=.9*\plotheight,width=\plotwidth,
             scale only axis,
             ylabel style = {yshift =0mm, xshift=0mm},
             xlabel style = {yshift =0mm, xshift=0mm},
             legend pos={north east}, 
             clip mode=individual,
             legend columns=1,
             label style={font=\scriptsize},
             tick label style={font=\scriptsize},
             legend style={font=\scriptsize\selectfont},
             legend style={at={(0.99,0.70)}},
             ]  

     \pgfmathsetmacro{\scale}{1}

     \addplot[color=blue!75!black,line width=1,densely dotted]
              table[x expr = \thisrow{xrcv}, 
                    y expr = \scale*\thisrow{ureal}, 
                   ]
              {\datafile}; \addlegendentry{real part}    

     \addplot[color=red!80!black,line width=1,densely dashed] 
              table[x expr = \thisrow{xrcv}, 
                    y expr = \scale*\thisrow{uimag}]
              {\datafile}; \addlegendentry{imaginary part}    
%
%
%
              
\end{axis}
\end{tikzpicture}
    \caption{Scattered solution measured at the \num{201} 
             receivers positioned at fixed height $x_2=\num{6}$.}
  \end{subfigure}
  \caption{Illustration of the acquisition setup and  
           generated data. The computations are carried
           out on the domain $[-20,20]\times[-20,20]$.
           While FWI uses the total field, 
           the reconstruction based upon Born and Rytov 
           approximations use the scattered solutions,
           obtained after removing a reference solution 
           corresponding to a propagation in an homogeneous 
           medium, cf.~\autoref{sec:data}.}
  \label{fig:fwi:heart-model_f05_data}
\end{figure}

\subsubsection{Reconstruction using FWI}
\label{sec:fwi:heart-shapes}

We carry out the reconstruction following \autoref{algo:FWI}, and 
the results are pictured in \autoref{fig:fwi:heart-model_f05_fwi},
where we compare the use of single and multi-frequency data.
In this example, we see that with relatively low frequency data 
(i.e., relatively large wavelength), such as 
for frequency $\omega/(2\pi)=\num{0.7}$ and $\omega/(2\pi)=\num{1}$, the reconstruction is smooth, 
see \autoref{fig:fwi:heart-model_f05_fwi_0.7} and \autoref{fig:fwi:heart-model_f05_fwi_1},
and one needs to use smaller wavelengths to obtain a better reconstruction, 
see \autoref{fig:fwi:heart-model_f05_fwi_1.2} and \autoref{fig:fwi:heart-model_f05_fwi_1.4}.
In \autoref{fig:fwi:heart-model_f05_fwi_multi}, we see that multi-frequency data 
gives the best reconstruction, it is also the most robust as one does not need
to anticipate the appropriate wavelength before carrying out the reconstruction.
Here both the shapes and contrast in amplitude are accurately obtained.
We notice some oscillatory noise in the reconstructed models, which 
could certainly be reduced by incorporating a regularization criterion
in the minimization \cite{Faucher2020EV}.

\setlength{\modelwidth}{4.50cm}
\setlength{\modelheight}{\modelwidth}
\setlength{\plotwidth}{10cm}
\setlength{\plotheight}{3cm}
\begin{figure}[ht!] \centering
  \pgfmathsetmacro{\xminloc}{-5.50}\pgfmathsetmacro{\xmaxloc}{5.50}
  \pgfmathsetmacro{\zminloc}{-5.50}\pgfmathsetmacro{\zmaxloc}{5.50}
  \pgfmathsetmacro{\cmin} {-0.1} \pgfmathsetmacro{\cmax} {0.5}

  \renewcommand{\modelfile}{images/fwi/model_heart/inversion/fwi_cp1_f0.5_unifreq0.7hz_f_scale-01-05}
  \begin{subfigure}[t]{.31\textwidth}\centering
    \begin{tikzpicture}

\pgfmathsetmacro{\xmingb} {-20} \pgfmathsetmacro{\xmaxgb}{20}
\pgfmathsetmacro{\zmingb} {-20} \pgfmathsetmacro{\zmaxgb}{20}

\begin{axis}[
  width=\modelwidth, height=\modelheight,
  axis on top, separate axis lines,
  xmin=\xminloc, xmax=\xmaxloc, 
  ymin=\zminloc, ymax=\zmaxloc, 
  x label style={xshift=-0.0cm, yshift= 0.00cm}, 
  y label style={xshift= 0.0cm, yshift=-0.30cm},
  colormap/jet,colorbar,colorbar style={title={\scriptsize $f$},
  title style={yshift=-2mm, xshift=0mm},
  width=.25cm, xshift=-0.5em},
  point meta min=\cmin,point meta max=\cmax,
  label style={font=\scriptsize},
  tick label style={font=\scriptsize},
  legend style={font=\scriptsize\selectfont},
]
\addplot [forget plot] graphics [xmin=\xmingb,xmax=\xmaxgb,ymin=\zmingb,ymax=\zmaxgb] {{\modelfile}.png};


\end{axis}
\end{tikzpicture}%
    \caption{Using single-frequency,\\ $\omega/(2\pi)=\num{0.7}$
    (PSNR \num{22.38}).}
    \label{fig:fwi:heart-model_f05_fwi_0.7}
  \end{subfigure}\hfill
  \renewcommand{\modelfile}{images/fwi/model_heart/inversion/fwi_cp1_f0.5_unifreq1hz_f_scale-01-05}
  \begin{subfigure}[t]{.31\textwidth}\centering
    \begin{tikzpicture}

\pgfmathsetmacro{\xmingb} {-20} \pgfmathsetmacro{\xmaxgb}{20}
\pgfmathsetmacro{\zmingb} {-20} \pgfmathsetmacro{\zmaxgb}{20}

\begin{axis}[
  width=\modelwidth, height=\modelheight,
  axis on top, separate axis lines,
  xmin=\xminloc, xmax=\xmaxloc, 
  ymin=\zminloc, ymax=\zmaxloc, 
  x label style={xshift=-0.0cm, yshift= 0.00cm}, 
  y label style={xshift= 0.0cm, yshift=-0.30cm},
  colormap/jet,colorbar,colorbar style={title={\scriptsize $f$},
  title style={yshift=-2mm, xshift=0mm},
  width=.25cm, xshift=-0.5em},
  point meta min=\cmin,point meta max=\cmax,
  label style={font=\scriptsize},
  tick label style={font=\scriptsize},
  legend style={font=\scriptsize\selectfont},
]
\addplot [forget plot] graphics [xmin=\xmingb,xmax=\xmaxgb,ymin=\zmingb,ymax=\zmaxgb] {{\modelfile}.png};


\end{axis}
\end{tikzpicture}%
    \caption{Using single-frequency,\\ $\omega/(2\pi)=\num{1}$.
    (PSNR \num{22.91}).}
    \label{fig:fwi:heart-model_f05_fwi_1}
  \end{subfigure}\hfill
  \renewcommand{\modelfile}{images/fwi/model_heart/inversion/fwi_cp1_f0.5_unifreq1.2hz_f_scale-01-05}
  \begin{subfigure}[t]{.31\textwidth}\centering
    \begin{tikzpicture}

\pgfmathsetmacro{\xmingb} {-20} \pgfmathsetmacro{\xmaxgb}{20}
\pgfmathsetmacro{\zmingb} {-20} \pgfmathsetmacro{\zmaxgb}{20}

\begin{axis}[
  width=\modelwidth, height=\modelheight,
  axis on top, separate axis lines,
  xmin=\xminloc, xmax=\xmaxloc, 
  ymin=\zminloc, ymax=\zmaxloc, 
  x label style={xshift=-0.0cm, yshift= 0.00cm}, 
  y label style={xshift= 0.0cm, yshift=-0.30cm},
  colormap/jet,colorbar,colorbar style={title={\scriptsize $f$},
  title style={yshift=-2mm, xshift=0mm},
  width=.25cm, xshift=-0.5em},
  point meta min=\cmin,point meta max=\cmax,
  label style={font=\scriptsize},
  tick label style={font=\scriptsize},
  legend style={font=\scriptsize\selectfont},
]
\addplot [forget plot] graphics [xmin=\xmingb,xmax=\xmaxgb,ymin=\zmingb,ymax=\zmaxgb] {{\modelfile}.png};


\end{axis}
\end{tikzpicture}%
    \caption{Using single-frequency,\\ $\omega/(2\pi)=\num{1.2}$.
    (PSNR \num{23.10}).}
    \label{fig:fwi:heart-model_f05_fwi_1.2}
  \end{subfigure}

  \renewcommand{\modelfile}{images/fwi/model_heart/inversion/fwi_cp1_f0.5_unifreq1.4hz_f_scale-01-05}
  \begin{subfigure}[t]{.31\textwidth}\centering
    \begin{tikzpicture}

\pgfmathsetmacro{\xmingb} {-20} \pgfmathsetmacro{\xmaxgb}{20}
\pgfmathsetmacro{\zmingb} {-20} \pgfmathsetmacro{\zmaxgb}{20}

\begin{axis}[
  width=\modelwidth, height=\modelheight,
  axis on top, separate axis lines,
  xmin=\xminloc, xmax=\xmaxloc, 
  ymin=\zminloc, ymax=\zmaxloc, 
  x label style={xshift=-0.0cm, yshift= 0.00cm}, 
  y label style={xshift= 0.0cm, yshift=-0.30cm},
  colormap/jet,colorbar,colorbar style={title={\scriptsize $f$},
  title style={yshift=-2mm, xshift=0mm},
  width=.25cm, xshift=-0.5em},
  point meta min=\cmin,point meta max=\cmax,
  label style={font=\scriptsize},
  tick label style={font=\scriptsize},
  legend style={font=\scriptsize\selectfont},
]
\addplot [forget plot] graphics [xmin=\xmingb,xmax=\xmaxgb,ymin=\zmingb,ymax=\zmaxgb] {{\modelfile}.png};


\end{axis}
\end{tikzpicture}%
    \caption{Using single-frequency, $\omega/(2\pi)=\num{1.4}$.
    (PSNR \num{23.31}).}
    \label{fig:fwi:heart-model_f05_fwi_1.4}
  \end{subfigure}
  \renewcommand{\modelfile}{images/fwi/model_heart/inversion/fwi_cp1_f05_multifreq0.7-1.4hz_f_scale-01-05}
  \begin{subfigure}[t]{.31\textwidth}\centering
    \begin{tikzpicture}

\pgfmathsetmacro{\xmingb} {-20} \pgfmathsetmacro{\xmaxgb}{20}
\pgfmathsetmacro{\zmingb} {-20} \pgfmathsetmacro{\zmaxgb}{20}

\begin{axis}[
  width=\modelwidth, height=\modelheight,
  axis on top, separate axis lines,
  xmin=\xminloc, xmax=\xmaxloc, 
  ymin=\zminloc, ymax=\zmaxloc, 
  x label style={xshift=-0.0cm, yshift= 0.00cm}, 
  y label style={xshift= 0.0cm, yshift=-0.30cm},
  colormap/jet,colorbar,colorbar style={title={\scriptsize $f$},
  title style={yshift=-2mm, xshift=0mm},
  width=.25cm, xshift=-0.5em},
  point meta min=\cmin,point meta max=\cmax,
  label style={font=\scriptsize},
  tick label style={font=\scriptsize},
  legend style={font=\scriptsize\selectfont},
]
\addplot [forget plot] graphics [xmin=\xmingb,xmax=\xmaxgb,ymin=\zmingb,ymax=\zmaxgb] {{\modelfile}.png};


\end{axis}
\end{tikzpicture}%
    \caption{Using multi-frequency, $\omega/(2\pi)\in\{\num{0.7}, \num{1}, \num{1.2}, \num{1.4} \}$.
    (PSNR \num{27.28}).}
    \label{fig:fwi:heart-model_f05_fwi_multi}
  \end{subfigure}
  \caption{Reconstruction of the model 
           \autoref{fig:fwi:heart-model_f05_data_1} 
           with FWI and starting from a homogeneous
           background with $f=0$. The models are given
           at frequency $\omega/(2\pi)=1$ and the wave speed
           is equal to \num{1} in the background.}
  \label{fig:fwi:heart-model_f05_fwi}
\end{figure}
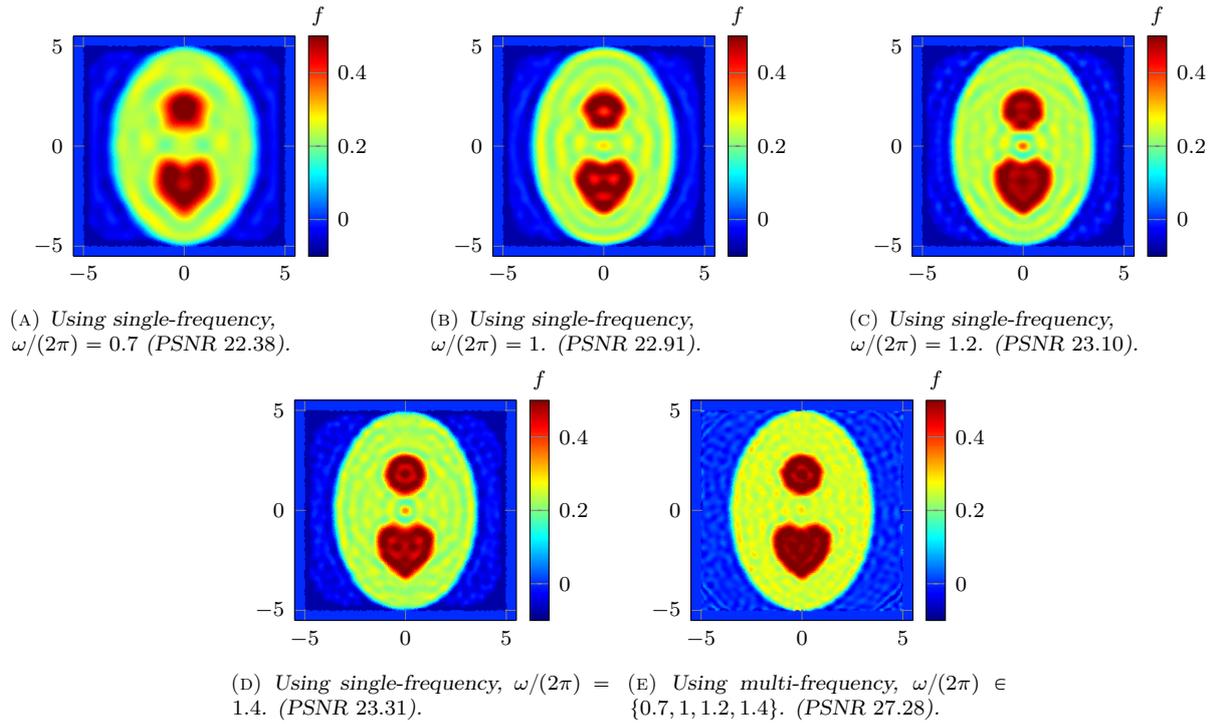

In \autoref{fig:fwi:heart-model_f2_fwi}, we conduct a 
similar computational experiment, but increasing the contrast in the 
included heart shape where $f$ has now a value of 
\num{2}, see \autoref{fig:fwi:heart-model_f2_fwi_a}.
We provide single and multi-frequency reconstructions,
and observe that large wavelengths still provide
a smooth reconstruction. The high contrast in the 
heart is well recovered.

\setlength{\modelwidth}{4.50cm}
\setlength{\modelheight}{\modelwidth}
\begin{figure}[ht!] \centering
  \pgfmathsetmacro{\xminloc}{-5.50}\pgfmathsetmacro{\xmaxloc}{5.50}
  \pgfmathsetmacro{\zminloc}{-5.50}\pgfmathsetmacro{\zmaxloc}{5.50}
  \pgfmathsetmacro{\cmin} {-0.1} \pgfmathsetmacro{\cmax}  {2}
  \renewcommand{\modelfile}{images/fwi/model_heart/models/cp1_f2_f_scale-01-02}
  \begin{subfigure}[t]{.31\textwidth}\centering
    \begin{tikzpicture}

\pgfmathsetmacro{\xmingb} {-20} \pgfmathsetmacro{\xmaxgb}{20}
\pgfmathsetmacro{\zmingb} {-20} \pgfmathsetmacro{\zmaxgb}{20}

\begin{axis}[
  width=\modelwidth, height=\modelheight,
  axis on top, separate axis lines,
  xmin=\xminloc, xmax=\xmaxloc, 
  ymin=\zminloc, ymax=\zmaxloc, 
  x label style={xshift=-0.0cm, yshift= 0.00cm}, 
  y label style={xshift= 0.0cm, yshift=-0.30cm},
  colormap/jet,colorbar,colorbar style={title={\scriptsize $f$},
  title style={yshift=-2mm, xshift=0mm},
  width=.25cm, xshift=-0.5em},
  point meta min=\cmin,point meta max=\cmax,
  label style={font=\scriptsize},
  tick label style={font=\scriptsize},
  legend style={font=\scriptsize\selectfont},
]
\addplot [forget plot] graphics [xmin=\xmingb,xmax=\xmaxgb,ymin=\zmingb,ymax=\zmaxgb] {{\modelfile}.png};


\end{axis}
\end{tikzpicture}%
    \caption{True model.}
    \label{fig:fwi:heart-model_f2_fwi_a}
  \end{subfigure}\hfill
  \renewcommand{\modelfile}{images/fwi/model_heart/inversion/fwi_cp1_f2_unifreq0.7hz_f_scale-01-2}
  \begin{subfigure}[t]{.31\textwidth}\centering
    \begin{tikzpicture}

\pgfmathsetmacro{\xmingb} {-20} \pgfmathsetmacro{\xmaxgb}{20}
\pgfmathsetmacro{\zmingb} {-20} \pgfmathsetmacro{\zmaxgb}{20}

\begin{axis}[
  width=\modelwidth, height=\modelheight,
  axis on top, separate axis lines,
  xmin=\xminloc, xmax=\xmaxloc, 
  ymin=\zminloc, ymax=\zmaxloc, 
  x label style={xshift=-0.0cm, yshift= 0.00cm}, 
  y label style={xshift= 0.0cm, yshift=-0.30cm},
  colormap/jet,colorbar,colorbar style={title={\scriptsize $f$},
  title style={yshift=-2mm, xshift=0mm},
  width=.25cm, xshift=-0.5em},
  point meta min=\cmin,point meta max=\cmax,
  label style={font=\scriptsize},
  tick label style={font=\scriptsize},
  legend style={font=\scriptsize\selectfont},
]
\addplot [forget plot] graphics [xmin=\xmingb,xmax=\xmaxgb,ymin=\zmingb,ymax=\zmaxgb] {{\modelfile}.png};


\end{axis}
\end{tikzpicture}%
    \caption{Using single-frequency,\\ $\omega/(2\pi)=\num{0.7}$
            (PSNR \num{25.04}).}
  \end{subfigure}\hfill
  \renewcommand{\modelfile}{images/fwi/model_heart/inversion/fwi_cp1_f2_unifreq1hz_f_scale-01-2}
  \begin{subfigure}[t]{.31\textwidth}\centering
    \begin{tikzpicture}

\pgfmathsetmacro{\xmingb} {-20} \pgfmathsetmacro{\xmaxgb}{20}
\pgfmathsetmacro{\zmingb} {-20} \pgfmathsetmacro{\zmaxgb}{20}

\begin{axis}[
  width=\modelwidth, height=\modelheight,
  axis on top, separate axis lines,
  xmin=\xminloc, xmax=\xmaxloc, 
  ymin=\zminloc, ymax=\zmaxloc, 
  x label style={xshift=-0.0cm, yshift= 0.00cm}, 
  y label style={xshift= 0.0cm, yshift=-0.30cm},
  colormap/jet,colorbar,colorbar style={title={\scriptsize $f$},
  title style={yshift=-2mm, xshift=0mm},
  width=.25cm, xshift=-0.5em},
  point meta min=\cmin,point meta max=\cmax,
  label style={font=\scriptsize},
  tick label style={font=\scriptsize},
  legend style={font=\scriptsize\selectfont},
]
\addplot [forget plot] graphics [xmin=\xmingb,xmax=\xmaxgb,ymin=\zmingb,ymax=\zmaxgb] {{\modelfile}.png};


\end{axis}
\end{tikzpicture}%
    \caption{Using single-frequency,\\ $\omega/(2\pi)=\num{1}$
            (PSNR \num{25.91}).}
  \end{subfigure}

  \renewcommand{\modelfile}{images/fwi/model_heart/inversion/fwi_cp1_f2_unifreq1.2hz_f_scale-01-2}
  \begin{subfigure}[t]{.31\textwidth}\centering
    \begin{tikzpicture}

\pgfmathsetmacro{\xmingb} {-20} \pgfmathsetmacro{\xmaxgb}{20}
\pgfmathsetmacro{\zmingb} {-20} \pgfmathsetmacro{\zmaxgb}{20}

\begin{axis}[
  width=\modelwidth, height=\modelheight,
  axis on top, separate axis lines,
  xmin=\xminloc, xmax=\xmaxloc, 
  ymin=\zminloc, ymax=\zmaxloc, 
  x label style={xshift=-0.0cm, yshift= 0.00cm}, 
  y label style={xshift= 0.0cm, yshift=-0.30cm},
  colormap/jet,colorbar,colorbar style={title={\scriptsize $f$},
  title style={yshift=-2mm, xshift=0mm},
  width=.25cm, xshift=-0.5em},
  point meta min=\cmin,point meta max=\cmax,
  label style={font=\scriptsize},
  tick label style={font=\scriptsize},
  legend style={font=\scriptsize\selectfont},
]
\addplot [forget plot] graphics [xmin=\xmingb,xmax=\xmaxgb,ymin=\zmingb,ymax=\zmaxgb] {{\modelfile}.png};


\end{axis}
\end{tikzpicture}%
    \caption{Using single-frequency,\\ $\omega/(2\pi)=\num{1.2}$
            (PSNR \num{26.19}).}
  \end{subfigure}\hfill
  \renewcommand{\modelfile}{images/fwi/model_heart/inversion/fwi_cp1_f2_unifreq1.4hz_f_scale-01-2}
  \begin{subfigure}[t]{.31\textwidth}\centering
    \begin{tikzpicture}

\pgfmathsetmacro{\xmingb} {-20} \pgfmathsetmacro{\xmaxgb}{20}
\pgfmathsetmacro{\zmingb} {-20} \pgfmathsetmacro{\zmaxgb}{20}

\begin{axis}[
  width=\modelwidth, height=\modelheight,
  axis on top, separate axis lines,
  xmin=\xminloc, xmax=\xmaxloc, 
  ymin=\zminloc, ymax=\zmaxloc, 
  x label style={xshift=-0.0cm, yshift= 0.00cm}, 
  y label style={xshift= 0.0cm, yshift=-0.30cm},
  colormap/jet,colorbar,colorbar style={title={\scriptsize $f$},
  title style={yshift=-2mm, xshift=0mm},
  width=.25cm, xshift=-0.5em},
  point meta min=\cmin,point meta max=\cmax,
  label style={font=\scriptsize},
  tick label style={font=\scriptsize},
  legend style={font=\scriptsize\selectfont},
]
\addplot [forget plot] graphics [xmin=\xmingb,xmax=\xmaxgb,ymin=\zmingb,ymax=\zmaxgb] {{\modelfile}.png};


\end{axis}
\end{tikzpicture}%
    \caption{Using single-frequency,\\ $\omega/(2\pi)=\num{1.4}$
            (PSNR \num{26.75}).}
  \end{subfigure}\hfill
  \renewcommand{\modelfile}{images/fwi/model_heart/inversion/fwi_cp1_f2_multifreq0.7-1.4hz_f_scale-01-2}
  \begin{subfigure}[t]{.31\textwidth}\centering
    \begin{tikzpicture}

\pgfmathsetmacro{\xmingb} {-20} \pgfmathsetmacro{\xmaxgb}{20}
\pgfmathsetmacro{\zmingb} {-20} \pgfmathsetmacro{\zmaxgb}{20}

\begin{axis}[
  width=\modelwidth, height=\modelheight,
  axis on top, separate axis lines,
  xmin=\xminloc, xmax=\xmaxloc, 
  ymin=\zminloc, ymax=\zmaxloc, 
  x label style={xshift=-0.0cm, yshift= 0.00cm}, 
  y label style={xshift= 0.0cm, yshift=-0.30cm},
  colormap/jet,colorbar,colorbar style={title={\scriptsize $f$},
  title style={yshift=-2mm, xshift=0mm},
  width=.25cm, xshift=-0.5em},
  point meta min=\cmin,point meta max=\cmax,
  label style={font=\scriptsize},
  tick label style={font=\scriptsize},
  legend style={font=\scriptsize\selectfont},
]
\addplot [forget plot] graphics [xmin=\xmingb,xmax=\xmaxgb,ymin=\zmingb,ymax=\zmaxgb] {{\modelfile}.png};


\end{axis}
\end{tikzpicture}%
    \caption{Using multi-frequency,
             $\omega/(2\pi)\in\{\num{0.7}, \num{1}, \num{1.2}, \num{1.4} \}$
            (PSNR \num{28.76}).}
  \end{subfigure}
  \caption{Reconstruction with FWI starting 
           from a homogeneous background 
           with $f=0$. The models are given
           at frequency $\omega/(2\pi)=1$ and the wave speed
           is equal to \num{1} in the background.}
  \label{fig:fwi:heart-model_f2_fwi}
\end{figure}
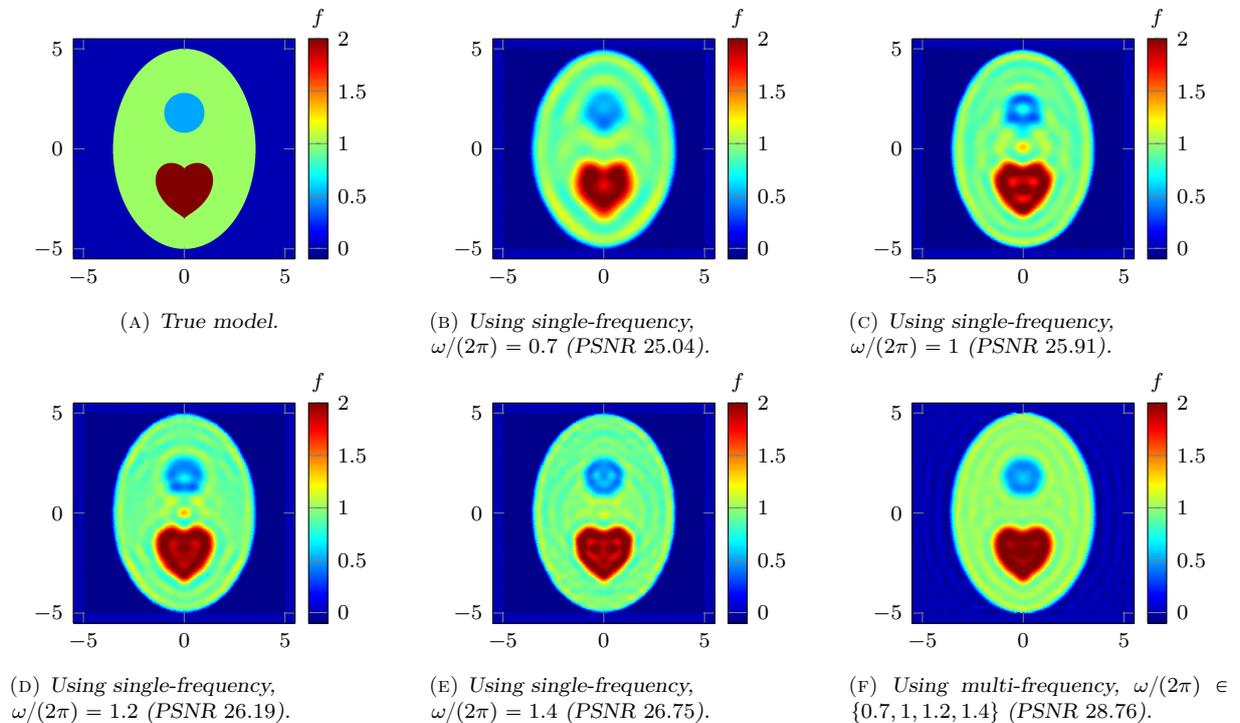

\subsubsection{Reconstruction using Born and Rytov approximations}
\label{sec:Born-heart-shapes}

We perform the reconstruction with \autoref{algo:infft} of the test model $f$ from \autoref{fig:fwi:heart-model_f05_data}.
In the following tests, we discretize $f$ on a finer grid of resolution $N=720$, which covers the radius $\rs = 15/\sqrt2$.
The PSNR is computed only on the central part of the grid that is visible in the image.
Since we know that the $f$ is real-valued, we take only the real part of the reconstruction.
For simplicity, we use a constant number of $12$ iterations in the conjugate gradient method of the inverse NDFT.

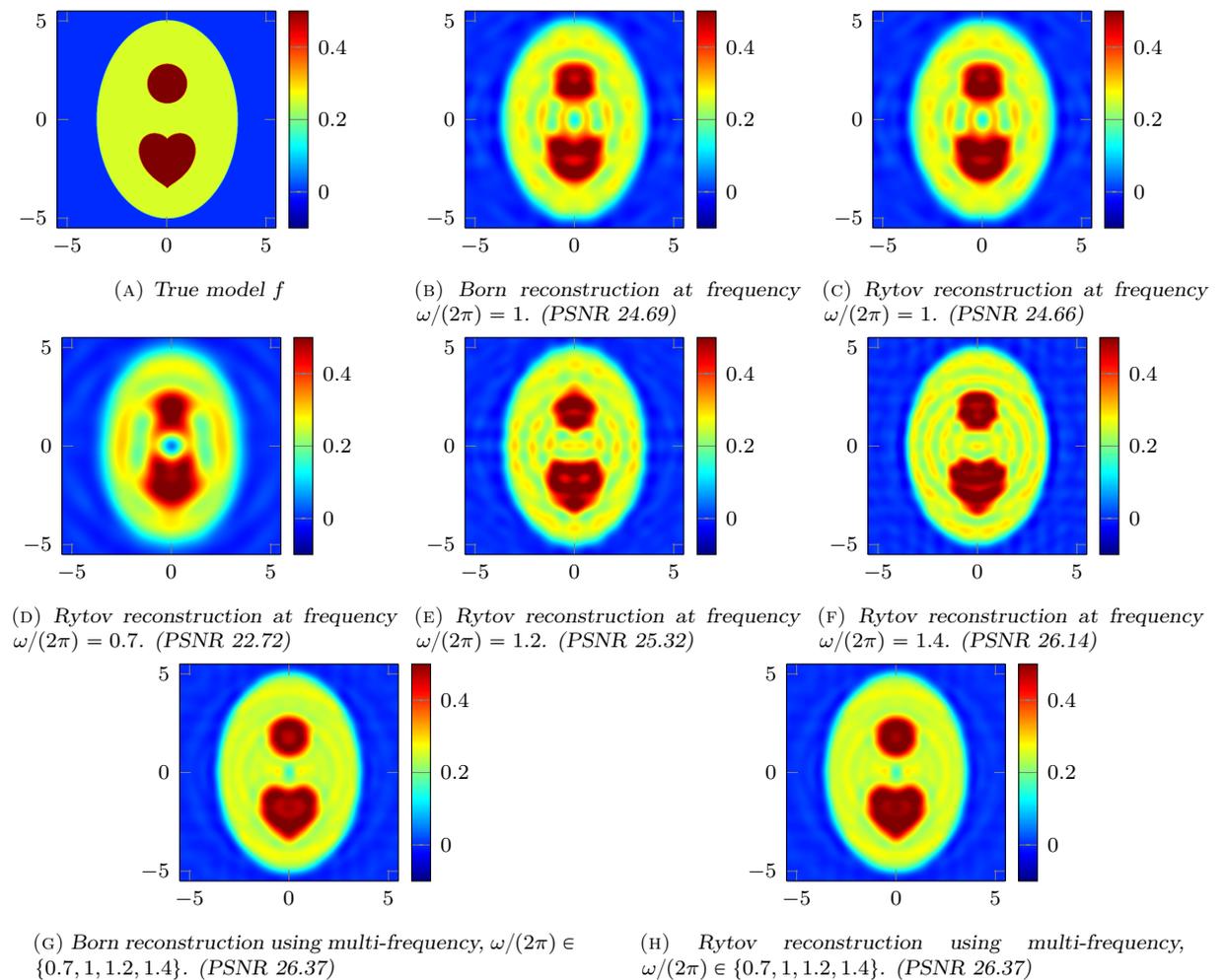
\begin{figure}[ht]\centering
  \pgfmathsetmacro{\xminloc}{-7.0711}\pgfmathsetmacro{\xmaxloc}{7.0711}
  \pgfmathsetmacro{\zminloc}{-7.0711}\pgfmathsetmacro{\zmaxloc}{7.0711}
  \pgfmathsetmacro{\xminpic}{-5.5000}\pgfmathsetmacro{\xmaxpic}{5.5000}
  \pgfmathsetmacro{\zminpic}{-5.5000}\pgfmathsetmacro{\zmaxpic}{5.5000}
  \pgfmathsetmacro{\cmin} {-0.1} \pgfmathsetmacro{\cmax} { 0.5}
  
  \renewcommand{\modelfile}{images/hearts/heartsnoise-1000mHz}
  \begin{subfigure}[t]{.32\textwidth}
    \begin{tikzpicture}
  \begin{axis}[
    width=\modelwidth, height=\modelheight,
    axis on top, separate axis lines,
    xmin=\xminpic, xmax=\xmaxpic, 
    ymin=\zminpic, ymax=\zmaxpic, 
    x label style={xshift=-0.0cm, yshift= 0.00cm}, 
    y label style={xshift= 0.0cm, yshift=-0.30cm},
    colormap/jet,colorbar,colorbar style={
      width=.25cm, xshift=-.5em},
    point meta min=\cmin,point meta max=\cmax,
    label style={font=\scriptsize},
    tick label style={font=\scriptsize},
    legend style={font=\scriptsize\selectfont},
    ]
    \addplot [forget plot] graphics [xmin=\xminloc,xmax=\xmaxloc,ymin=\zminloc,ymax=\zmaxloc] {{\modelfile}.png};
  \end{axis}
\end{tikzpicture}%
    \caption{True model $f$}
  \end{subfigure}\hfill
  \renewcommand{\modelfile}{images/hearts/heartsnoise-1000mHz-rec}
  \begin{subfigure}[t]{.32\textwidth}
  \begin{tikzpicture}
  \begin{axis}[
    width=\modelwidth, height=\modelheight,
    axis on top, separate axis lines,
    xmin=\xminpic, xmax=\xmaxpic, 
    ymin=\zminpic, ymax=\zmaxpic, 
    x label style={xshift=-0.0cm, yshift= 0.00cm}, 
    y label style={xshift= 0.0cm, yshift=-0.30cm},
    colormap/jet,colorbar,colorbar style={
      width=.25cm, xshift=-.5em},
    point meta min=\cmin,point meta max=\cmax,
    label style={font=\scriptsize},
    tick label style={font=\scriptsize},
    legend style={font=\scriptsize\selectfont},
    ]
    \addplot [forget plot] graphics [xmin=\xminloc,xmax=\xmaxloc,ymin=\zminloc,ymax=\zmaxloc] {{\modelfile}.png};
  \end{axis}
\end{tikzpicture}%
  \caption{Born reconstruction at frequency $\omega/(2\pi)={1} $. (PSNR 24.69)}
  \end{subfigure}\hfill
  \renewcommand{\modelfile}{images/hearts/heartsnoise-1000mHz-rec-Rytov}
  \begin{subfigure}[t]{.32\textwidth}
  \begin{tikzpicture}
  \begin{axis}[
    width=\modelwidth, height=\modelheight,
    axis on top, separate axis lines,
    xmin=\xminpic, xmax=\xmaxpic, 
    ymin=\zminpic, ymax=\zmaxpic, 
    x label style={xshift=-0.0cm, yshift= 0.00cm}, 
    y label style={xshift= 0.0cm, yshift=-0.30cm},
    colormap/jet,colorbar,colorbar style={
      width=.25cm, xshift=-.5em},
    point meta min=\cmin,point meta max=\cmax,
    label style={font=\scriptsize},
    tick label style={font=\scriptsize},
    legend style={font=\scriptsize\selectfont},
    ]
    \addplot [forget plot] graphics [xmin=\xminloc,xmax=\xmaxloc,ymin=\zminloc,ymax=\zmaxloc] {{\modelfile}.png};
  \end{axis}
\end{tikzpicture}%
  \caption{Rytov reconstruction at frequency $\omega/(2\pi)={1} $. (PSNR 24.66)}
  \end{subfigure}
  \renewcommand{\modelfile}{images/hearts/heartsnoise-700mHz-rec-Rytov}
  \begin{subfigure}[t]{.32\textwidth}
  \begin{tikzpicture}
  \begin{axis}[
    width=\modelwidth, height=\modelheight,
    axis on top, separate axis lines,
    xmin=\xminpic, xmax=\xmaxpic, 
    ymin=\zminpic, ymax=\zmaxpic, 
    x label style={xshift=-0.0cm, yshift= 0.00cm}, 
    y label style={xshift= 0.0cm, yshift=-0.30cm},
    colormap/jet,colorbar,colorbar style={
      width=.25cm, xshift=-.5em},
    point meta min=\cmin,point meta max=\cmax,
    label style={font=\scriptsize},
    tick label style={font=\scriptsize},
    legend style={font=\scriptsize\selectfont},
    ]
    \addplot [forget plot] graphics [xmin=\xminloc,xmax=\xmaxloc,ymin=\zminloc,ymax=\zmaxloc] {{\modelfile}.png};
  \end{axis}
\end{tikzpicture}%
  \caption{Rytov reconstruction at frequency $\omega/(2\pi)={0.7} $. (PSNR 22.72)} 
  \end{subfigure}
  \renewcommand{\modelfile}{images/hearts/heartsnoise-1200mHz-rec-Rytov}
  \begin{subfigure}[t]{.32\textwidth}
  \begin{tikzpicture}
  \begin{axis}[
    width=\modelwidth, height=\modelheight,
    axis on top, separate axis lines,
    xmin=\xminpic, xmax=\xmaxpic, 
    ymin=\zminpic, ymax=\zmaxpic, 
    x label style={xshift=-0.0cm, yshift= 0.00cm}, 
    y label style={xshift= 0.0cm, yshift=-0.30cm},
    colormap/jet,colorbar,colorbar style={
      width=.25cm, xshift=-.5em},
    point meta min=\cmin,point meta max=\cmax,
    label style={font=\scriptsize},
    tick label style={font=\scriptsize},
    legend style={font=\scriptsize\selectfont},
    ]
    \addplot [forget plot] graphics [xmin=\xminloc,xmax=\xmaxloc,ymin=\zminloc,ymax=\zmaxloc] {{\modelfile}.png};
  \end{axis}
\end{tikzpicture}%
  \caption{Rytov reconstruction at frequency $\omega/(2\pi)={1.2} $. (PSNR 25.32)} 
  \end{subfigure}
  \renewcommand{\modelfile}{images/hearts/heartsnoise-1400mHz-rec-Rytov}
  \begin{subfigure}[t]{.32\textwidth}
  \begin{tikzpicture}
  \begin{axis}[
    width=\modelwidth, height=\modelheight,
    axis on top, separate axis lines,
    xmin=\xminpic, xmax=\xmaxpic, 
    ymin=\zminpic, ymax=\zmaxpic, 
    x label style={xshift=-0.0cm, yshift= 0.00cm}, 
    y label style={xshift= 0.0cm, yshift=-0.30cm},
    colormap/jet,colorbar,colorbar style={
      width=.25cm, xshift=-.5em},
    point meta min=\cmin,point meta max=\cmax,
    label style={font=\scriptsize},
    tick label style={font=\scriptsize},
    legend style={font=\scriptsize\selectfont},
    ]
    \addplot [forget plot] graphics [xmin=\xminloc,xmax=\xmaxloc,ymin=\zminloc,ymax=\zmaxloc] {{\modelfile}.png};
  \end{axis}
\end{tikzpicture}%
  \caption{Rytov reconstruction at frequency $\omega/(2\pi)={1.4} $. (PSNR 26.14)} 
  \end{subfigure}
  
  \renewcommand{\modelfile}{images/hearts/heartsnoise-700_1000_1200_1400mHz-rec}
  \begin{subfigure}[t]{.45\textwidth}\centering
  \begin{tikzpicture}
  \begin{axis}[
    width=\modelwidth, height=\modelheight,
    axis on top, separate axis lines,
    xmin=\xminpic, xmax=\xmaxpic, 
    ymin=\zminpic, ymax=\zmaxpic, 
    x label style={xshift=-0.0cm, yshift= 0.00cm}, 
    y label style={xshift= 0.0cm, yshift=-0.30cm},
    colormap/jet,colorbar,colorbar style={
      width=.25cm, xshift=-.5em},
    point meta min=\cmin,point meta max=\cmax,
    label style={font=\scriptsize},
    tick label style={font=\scriptsize},
    legend style={font=\scriptsize\selectfont},
    ]
    \addplot [forget plot] graphics [xmin=\xminloc,xmax=\xmaxloc,ymin=\zminloc,ymax=\zmaxloc] {{\modelfile}.png};
  \end{axis}
\end{tikzpicture}%
  \caption{Born reconstruction using multi-frequency, $\omega/(2\pi)\in\{0.7,1,1.2,1.4\}$. (PSNR 26.37)} 
  \end{subfigure}\qquad
  \renewcommand{\modelfile}{images/hearts/heartsnoise-700_1000_1200_1400mHz-rec-Rytov}
  \begin{subfigure}[t]{.45\textwidth}\centering
  \begin{tikzpicture}
  \begin{axis}[
    width=\modelwidth, height=\modelheight,
    axis on top, separate axis lines,
    xmin=\xminpic, xmax=\xmaxpic, 
    ymin=\zminpic, ymax=\zmaxpic, 
    x label style={xshift=-0.0cm, yshift= 0.00cm}, 
    y label style={xshift= 0.0cm, yshift=-0.30cm},
    colormap/jet,colorbar,colorbar style={
      width=.25cm, xshift=-.5em},
    point meta min=\cmin,point meta max=\cmax,
    label style={font=\scriptsize},
    tick label style={font=\scriptsize},
    legend style={font=\scriptsize\selectfont},
    ]
    \addplot [forget plot] graphics [xmin=\xminloc,xmax=\xmaxloc,ymin=\zminloc,ymax=\zmaxloc] {{\modelfile}.png};
  \end{axis}
\end{tikzpicture}%
  \caption{Rytov reconstruction using multi-frequency, $\omega/(2\pi)\in\{0.7,1,1.2,1.4\} $. (PSNR 26.37)} 
  \end{subfigure}

  \caption{Reconstructions for different frequencies of the incident wave.
  The PSNR is computed on the visible part of the grid for the real part of the reconstruction, since we know that $f$ must be real.
  \label{fig:heart-model-rec}}
\end{figure}

The Born and Rytov reconstructions are shown in \autoref{fig:heart-model-rec},
where the data $u$ is the same as in \autoref{sec:fwi:heart-shapes}.
The reconstruction with a higher frequency of the incident wave is more accurate, since it provides a larger k-space coverage, which is the disk of radius
$\sqrt{2} k_0 = \sqrt{2} \omega$, 
see \autoref{sec:coverage}.
Moreover, the multi-frequency reconstruction is shown in \autoref{fig:heart-model-rec} (g) and (h).
Even though the multi-frequency setup covers the same disk in k-space, it still seems superior because we have more data points of the Fourier transform $\mathcal Ff$.

\begin{figure}[ht]\centering
  \pgfmathsetmacro{\xminloc}{-7.0711}\pgfmathsetmacro{\xmaxloc}{7.0711}
  \pgfmathsetmacro{\zminloc}{-7.0711}\pgfmathsetmacro{\zmaxloc}{7.0711}
  \pgfmathsetmacro{\xminpic}{-5.5000}\pgfmathsetmacro{\xmaxpic}{5.5000}
  \pgfmathsetmacro{\zminpic}{-5.5000}\pgfmathsetmacro{\zmaxpic}{5.5000}
  \pgfmathsetmacro{\cmin} {-0.1} \pgfmathsetmacro{\cmax} { 2}
  
  \renewcommand{\modelfile}{images/hearts/hearts2noise-1000mHz}
  \begin{subfigure}[t]{.32\textwidth}
    \begin{tikzpicture}
  \begin{axis}[
    width=\modelwidth, height=\modelheight,
    axis on top, separate axis lines,
    xmin=\xminpic, xmax=\xmaxpic, 
    ymin=\zminpic, ymax=\zmaxpic, 
    x label style={xshift=-0.0cm, yshift= 0.00cm}, 
    y label style={xshift= 0.0cm, yshift=-0.30cm},
    colormap/jet,colorbar,colorbar style={
      width=.25cm, xshift=-.5em},
    point meta min=\cmin,point meta max=\cmax,
    label style={font=\scriptsize},
    tick label style={font=\scriptsize},
    legend style={font=\scriptsize\selectfont},
    ]
    \addplot [forget plot] graphics [xmin=\xminloc,xmax=\xmaxloc,ymin=\zminloc,ymax=\zmaxloc] {{\modelfile}.png};
  \end{axis}
\end{tikzpicture}%
    \caption{True model $f$}
  \end{subfigure}\hfill
  \renewcommand{\modelfile}{images/hearts/hearts2noise-1000mHz-rec}
  \begin{subfigure}[t]{.32\textwidth}
    \begin{tikzpicture}
  \begin{axis}[
    width=\modelwidth, height=\modelheight,
    axis on top, separate axis lines,
    xmin=\xminpic, xmax=\xmaxpic, 
    ymin=\zminpic, ymax=\zmaxpic, 
    x label style={xshift=-0.0cm, yshift= 0.00cm}, 
    y label style={xshift= 0.0cm, yshift=-0.30cm},
    colormap/jet,colorbar,colorbar style={
      width=.25cm, xshift=-.5em},
    point meta min=\cmin,point meta max=\cmax,
    label style={font=\scriptsize},
    tick label style={font=\scriptsize},
    legend style={font=\scriptsize\selectfont},
    ]
    \addplot [forget plot] graphics [xmin=\xminloc,xmax=\xmaxloc,ymin=\zminloc,ymax=\zmaxloc] {{\modelfile}.png};
  \end{axis}
\end{tikzpicture}%
    \caption{Born reconstruction at frequency $\omega/(2\pi)={1} $. (PSNR 23.53)}
  \end{subfigure}\hfill
  \renewcommand{\modelfile}{images/hearts/hearts2noise-1000mHz-rec-Rytov}
  \begin{subfigure}[t]{.32\textwidth}
    \begin{tikzpicture}
  \begin{axis}[
    width=\modelwidth, height=\modelheight,
    axis on top, separate axis lines,
    xmin=\xminpic, xmax=\xmaxpic, 
    ymin=\zminpic, ymax=\zmaxpic, 
    x label style={xshift=-0.0cm, yshift= 0.00cm}, 
    y label style={xshift= 0.0cm, yshift=-0.30cm},
    colormap/jet,colorbar,colorbar style={
      width=.25cm, xshift=-.5em},
    point meta min=\cmin,point meta max=\cmax,
    label style={font=\scriptsize},
    tick label style={font=\scriptsize},
    legend style={font=\scriptsize\selectfont},
    ]
    \addplot [forget plot] graphics [xmin=\xminloc,xmax=\xmaxloc,ymin=\zminloc,ymax=\zmaxloc] {{\modelfile}.png};
  \end{axis}
\end{tikzpicture}%
    \caption{Rytov reconstruction at frequency $\omega/(2\pi)={1} $. (PSNR 24.47)} 
  \end{subfigure}
  \renewcommand{\modelfile}{images/hearts/hearts2noise-700mHz-rec-Rytov}
  \begin{subfigure}[t]{.32\textwidth}
  \begin{tikzpicture}
  \begin{axis}[
    width=\modelwidth, height=\modelheight,
    axis on top, separate axis lines,
    xmin=\xminpic, xmax=\xmaxpic, 
    ymin=\zminpic, ymax=\zmaxpic, 
    x label style={xshift=-0.0cm, yshift= 0.00cm}, 
    y label style={xshift= 0.0cm, yshift=-0.30cm},
    colormap/jet,colorbar,colorbar style={
      width=.25cm, xshift=-.5em},
    point meta min=\cmin,point meta max=\cmax,
    label style={font=\scriptsize},
    tick label style={font=\scriptsize},
    legend style={font=\scriptsize\selectfont},
    ]
    \addplot [forget plot] graphics [xmin=\xminloc,xmax=\xmaxloc,ymin=\zminloc,ymax=\zmaxloc] {{\modelfile}.png};
  \end{axis}
\end{tikzpicture}%
  \caption{Rytov reconstruction at frequency $\omega/(2\pi)={0.7} $. (PSNR 21.78)} 
  \end{subfigure}
  \renewcommand{\modelfile}{images/hearts/hearts2noise-1200mHz-rec-Rytov}
  \begin{subfigure}[t]{.32\textwidth}
  \begin{tikzpicture}
  \begin{axis}[
    width=\modelwidth, height=\modelheight,
    axis on top, separate axis lines,
    xmin=\xminpic, xmax=\xmaxpic, 
    ymin=\zminpic, ymax=\zmaxpic, 
    x label style={xshift=-0.0cm, yshift= 0.00cm}, 
    y label style={xshift= 0.0cm, yshift=-0.30cm},
    colormap/jet,colorbar,colorbar style={
      width=.25cm, xshift=-.5em},
    point meta min=\cmin,point meta max=\cmax,
    label style={font=\scriptsize},
    tick label style={font=\scriptsize},
    legend style={font=\scriptsize\selectfont},
    ]
    \addplot [forget plot] graphics [xmin=\xminloc,xmax=\xmaxloc,ymin=\zminloc,ymax=\zmaxloc] {{\modelfile}.png};
  \end{axis}
\end{tikzpicture}%
  \caption{Rytov reconstruction at frequency $\omega/(2\pi)={1.2} $. (PSNR 25.55)} 
  \end{subfigure}
  \renewcommand{\modelfile}{images/hearts/hearts2noise-1400mHz-rec-Rytov}
  \begin{subfigure}[t]{.32\textwidth}
    \begin{tikzpicture}
  \begin{axis}[
    width=\modelwidth, height=\modelheight,
    axis on top, separate axis lines,
    xmin=\xminpic, xmax=\xmaxpic, 
    ymin=\zminpic, ymax=\zmaxpic, 
    x label style={xshift=-0.0cm, yshift= 0.00cm}, 
    y label style={xshift= 0.0cm, yshift=-0.30cm},
    colormap/jet,colorbar,colorbar style={
      width=.25cm, xshift=-.5em},
    point meta min=\cmin,point meta max=\cmax,
    label style={font=\scriptsize},
    tick label style={font=\scriptsize},
    legend style={font=\scriptsize\selectfont},
    ]
    \addplot [forget plot] graphics [xmin=\xminloc,xmax=\xmaxloc,ymin=\zminloc,ymax=\zmaxloc] {{\modelfile}.png};
  \end{axis}
\end{tikzpicture}%
    \caption{Rytov reconstruction at frequency $\omega/(2\pi)={1.4} $. (PSNR 26.40)} 
  \end{subfigure}
  
  \renewcommand{\modelfile}{images/hearts/hearts2noise-700_1000_1200_1400mHz-rec}
  \begin{subfigure}[t]{.45\textwidth}\centering
    \begin{tikzpicture}
  \begin{axis}[
    width=\modelwidth, height=\modelheight,
    axis on top, separate axis lines,
    xmin=\xminpic, xmax=\xmaxpic, 
    ymin=\zminpic, ymax=\zmaxpic, 
    x label style={xshift=-0.0cm, yshift= 0.00cm}, 
    y label style={xshift= 0.0cm, yshift=-0.30cm},
    colormap/jet,colorbar,colorbar style={
      width=.25cm, xshift=-.5em},
    point meta min=\cmin,point meta max=\cmax,
    label style={font=\scriptsize},
    tick label style={font=\scriptsize},
    legend style={font=\scriptsize\selectfont},
    ]
    \addplot [forget plot] graphics [xmin=\xminloc,xmax=\xmaxloc,ymin=\zminloc,ymax=\zmaxloc] {{\modelfile}.png};
  \end{axis}
\end{tikzpicture}%
    \caption{Born reconstruction using multi-frequency, $\omega/(2\pi)\in\{0.7,1,1.2,1.4\} $. (PSNR 23.77)} 
  \end{subfigure}\qquad
  \renewcommand{\modelfile}{images/hearts/hearts2noise-700_1000_1200_1400mHz-rec-Rytov}
  \begin{subfigure}[t]{.45\textwidth}\centering
    \begin{tikzpicture}
  \begin{axis}[
    width=\modelwidth, height=\modelheight,
    axis on top, separate axis lines,
    xmin=\xminpic, xmax=\xmaxpic, 
    ymin=\zminpic, ymax=\zmaxpic, 
    x label style={xshift=-0.0cm, yshift= 0.00cm}, 
    y label style={xshift= 0.0cm, yshift=-0.30cm},
    colormap/jet,colorbar,colorbar style={
      width=.25cm, xshift=-.5em},
    point meta min=\cmin,point meta max=\cmax,
    label style={font=\scriptsize},
    tick label style={font=\scriptsize},
    legend style={font=\scriptsize\selectfont},
    ]
    \addplot [forget plot] graphics [xmin=\xminloc,xmax=\xmaxloc,ymin=\zminloc,ymax=\zmaxloc] {{\modelfile}.png};
  \end{axis}
\end{tikzpicture}%
    \caption{Rytov reconstruction using multi-frequency, $\omega/(2\pi)\in\{0.7,1,1.2,1.4\} $. (PSNR 25.92)} 
  \end{subfigure}
  
  \caption{Reconstructions with a higher contrast, where the rest of the setting is the same as in \autoref{fig:heart-model-rec}.
    \label{fig:heart2-model-rec}}
\end{figure}
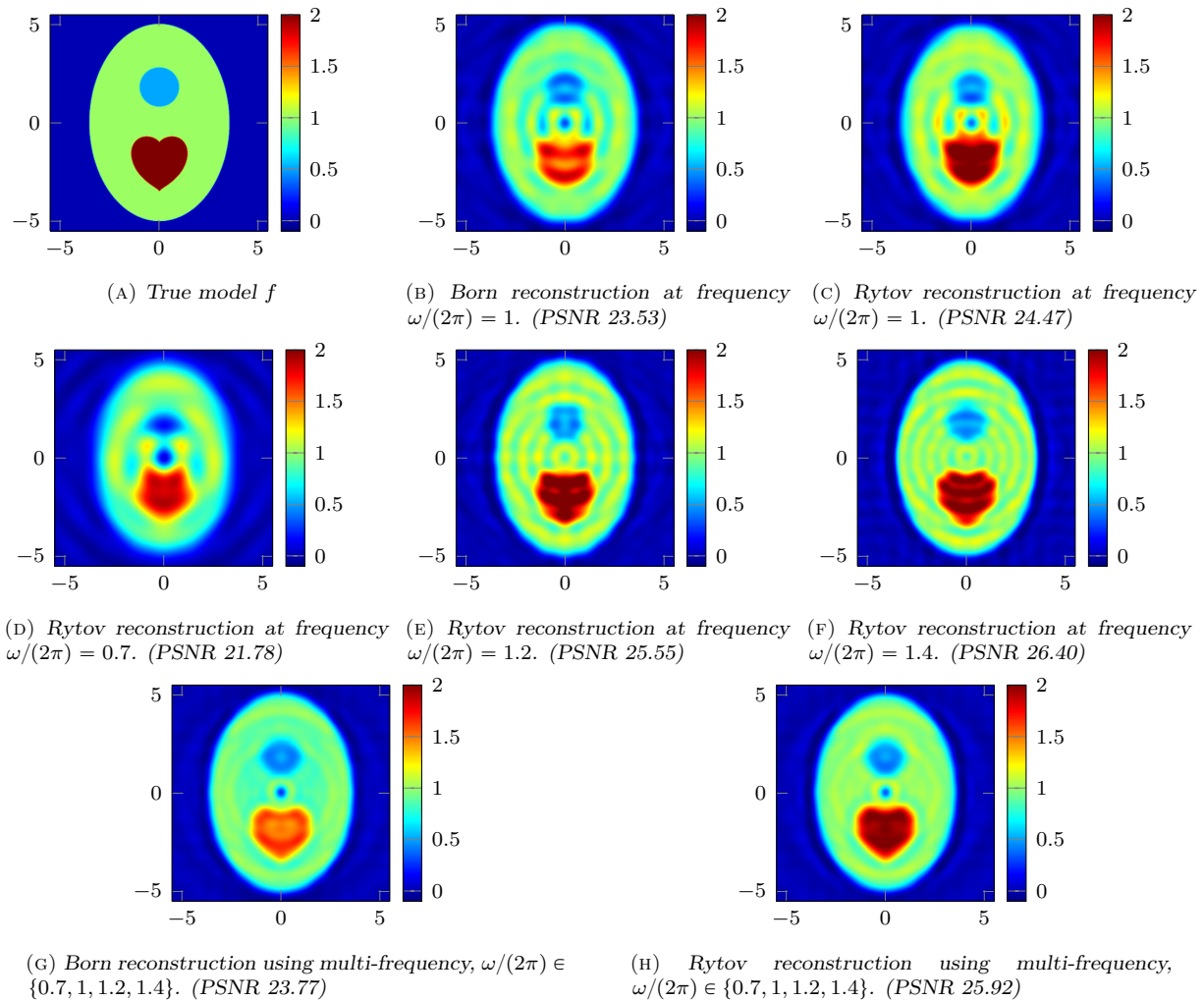

For the similar model from \autoref{fig:fwi:heart-model_f2_fwi_a} with a higher contrast,
the reconstructions with Born and Rytov approximation differ more from the FWI reconstruction because of the more severe scattering, see \autoref{fig:heart2-model-rec}.
In general, we can expect the FWI reconstruction to be better since it is a numerical approximation of the wave equation,
of which the Born or Rytov approximations are just linearizations.

\subsection{Reconstruction of embedded shapes: Phantom 2}

We now consider the case with combinations of smaller convex and 
non-convex shapes included in the background medium.

\subsubsection{Reconstruction using FWI}

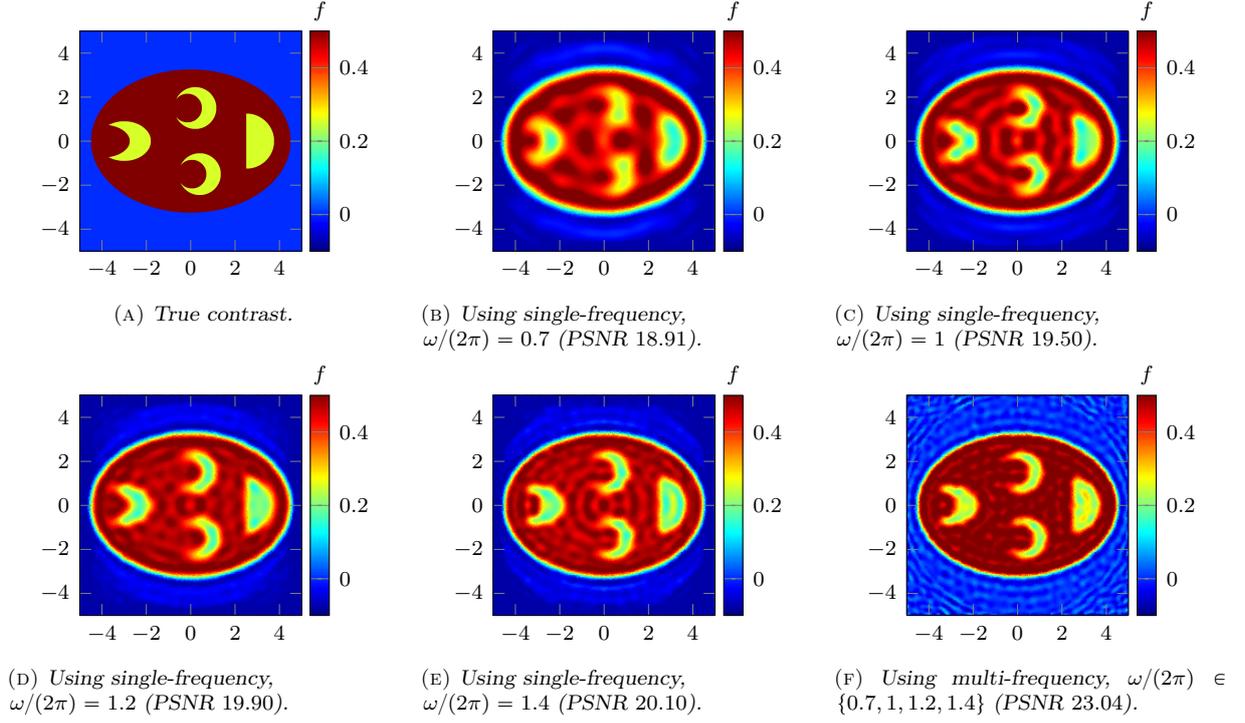
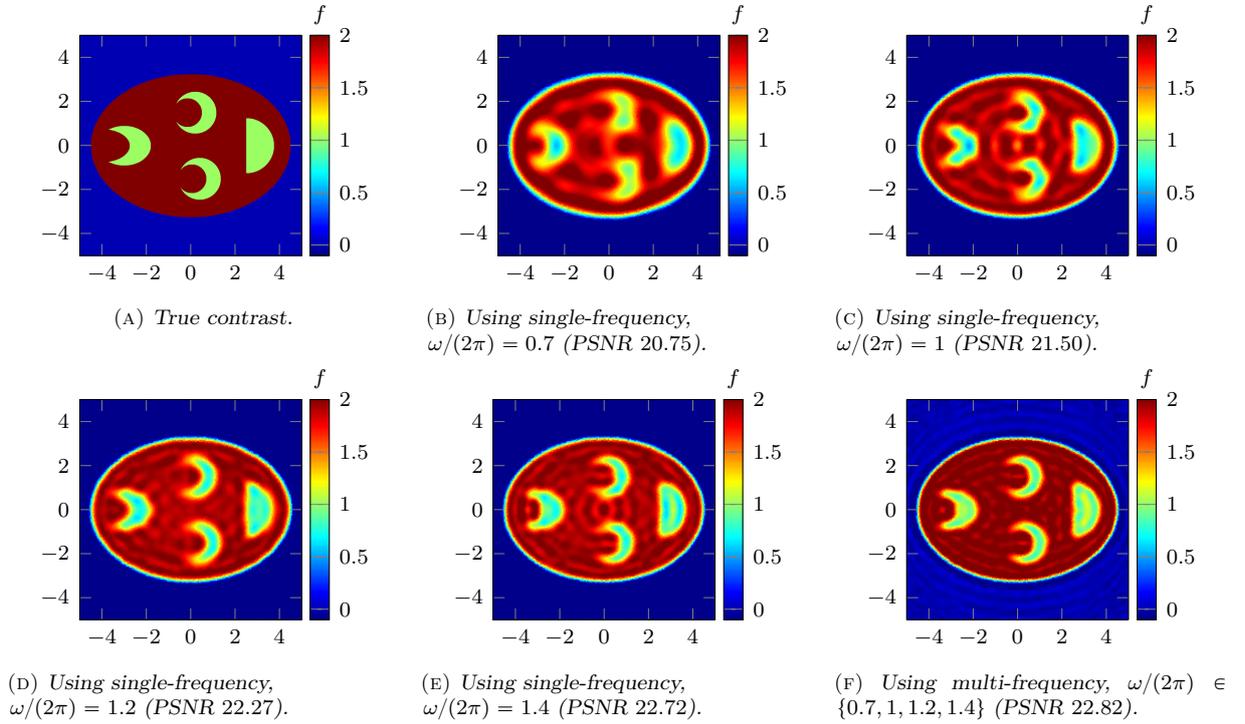
\begin{figure}[ht!] \centering
  \pgfmathsetmacro{\xminloc}{-5}\pgfmathsetmacro{\xmaxloc}{5}
  \pgfmathsetmacro{\zminloc}{-5}\pgfmathsetmacro{\zmaxloc}{5}
  \pgfmathsetmacro{\cmin} {-0.1} \pgfmathsetmacro{\cmax}  {0.5}
  \renewcommand{\modelfile}{images/fwi/model_shapes/models/cp1_f05_f_scale-01-05}
  \begin{subfigure}[t]{.32\textwidth}\centering
    \begin{tikzpicture}

\pgfmathsetmacro{\xmingb} {-20} \pgfmathsetmacro{\xmaxgb}{20}
\pgfmathsetmacro{\zmingb} {-20} \pgfmathsetmacro{\zmaxgb}{20}

\begin{axis}[
  width=\modelwidth, height=\modelheight,
  axis on top, separate axis lines,
  xmin=\xminloc, xmax=\xmaxloc, 
  ymin=\zminloc, ymax=\zmaxloc, 
  x label style={xshift=-0.0cm, yshift= 0.00cm}, 
  y label style={xshift= 0.0cm, yshift=-0.30cm},
  colormap/jet,colorbar,colorbar style={title={\scriptsize $f$},
  title style={yshift=-2mm, xshift=0mm},
  width=.25cm, xshift=-0.5em},
  point meta min=\cmin,point meta max=\cmax,
  label style={font=\scriptsize},
  tick label style={font=\scriptsize},
  legend style={font=\scriptsize\selectfont},
]
\addplot [forget plot] graphics [xmin=\xmingb,xmax=\xmaxgb,ymin=\zmingb,ymax=\zmaxgb] {{\modelfile}.png};


\end{axis}
\end{tikzpicture}%
    \caption{True contrast.}
  \end{subfigure}\hfill
  \renewcommand{\modelfile}{images/fwi/model_shapes/inversion/fwi_cp1_f0.5_unifreq0.7hz_f_scale-01-05}
  \begin{subfigure}[t]{.32\textwidth}\centering
    \begin{tikzpicture}

\pgfmathsetmacro{\xmingb} {-20} \pgfmathsetmacro{\xmaxgb}{20}
\pgfmathsetmacro{\zmingb} {-20} \pgfmathsetmacro{\zmaxgb}{20}

\begin{axis}[
  width=\modelwidth, height=\modelheight,
  axis on top, separate axis lines,
  xmin=\xminloc, xmax=\xmaxloc, 
  ymin=\zminloc, ymax=\zmaxloc, 
  x label style={xshift=-0.0cm, yshift= 0.00cm}, 
  y label style={xshift= 0.0cm, yshift=-0.30cm},
  colormap/jet,colorbar,colorbar style={title={\scriptsize $f$},
  title style={yshift=-2mm, xshift=0mm},
  width=.25cm, xshift=-0.5em},
  point meta min=\cmin,point meta max=\cmax,
  label style={font=\scriptsize},
  tick label style={font=\scriptsize},
  legend style={font=\scriptsize\selectfont},
]
\addplot [forget plot] graphics [xmin=\xmingb,xmax=\xmaxgb,ymin=\zmingb,ymax=\zmaxgb] {{\modelfile}.png};


\end{axis}
\end{tikzpicture}%
    \caption{Using single-frequency,\\ $\omega/(2\pi)=\num{0.7}$
            (PSNR \num{18.91}).}
  \end{subfigure}\hfill
  \renewcommand{\modelfile}{images/fwi/model_shapes/inversion/fwi_cp1_f0.5_unifreq1hz_f_scale-01-05}
  \begin{subfigure}[t]{.32\textwidth}\centering
    \begin{tikzpicture}

\pgfmathsetmacro{\xmingb} {-20} \pgfmathsetmacro{\xmaxgb}{20}
\pgfmathsetmacro{\zmingb} {-20} \pgfmathsetmacro{\zmaxgb}{20}

\begin{axis}[
  width=\modelwidth, height=\modelheight,
  axis on top, separate axis lines,
  xmin=\xminloc, xmax=\xmaxloc, 
  ymin=\zminloc, ymax=\zmaxloc, 
  x label style={xshift=-0.0cm, yshift= 0.00cm}, 
  y label style={xshift= 0.0cm, yshift=-0.30cm},
  colormap/jet,colorbar,colorbar style={title={\scriptsize $f$},
  title style={yshift=-2mm, xshift=0mm},
  width=.25cm, xshift=-0.5em},
  point meta min=\cmin,point meta max=\cmax,
  label style={font=\scriptsize},
  tick label style={font=\scriptsize},
  legend style={font=\scriptsize\selectfont},
]
\addplot [forget plot] graphics [xmin=\xmingb,xmax=\xmaxgb,ymin=\zmingb,ymax=\zmaxgb] {{\modelfile}.png};


\end{axis}
\end{tikzpicture}%
    \caption{Using single-frequency,\\ $\omega/(2\pi)=\num{1}$
            (PSNR \num{19.50}).}
  \end{subfigure}
  
  \renewcommand{\modelfile}{images/fwi/model_shapes/inversion/fwi_cp1_f0.5_unifreq1.2hz_f_scale-01-05}
  \begin{subfigure}[t]{.32\textwidth}\centering
    \begin{tikzpicture}

\pgfmathsetmacro{\xmingb} {-20} \pgfmathsetmacro{\xmaxgb}{20}
\pgfmathsetmacro{\zmingb} {-20} \pgfmathsetmacro{\zmaxgb}{20}

\begin{axis}[
  width=\modelwidth, height=\modelheight,
  axis on top, separate axis lines,
  xmin=\xminloc, xmax=\xmaxloc, 
  ymin=\zminloc, ymax=\zmaxloc, 
  x label style={xshift=-0.0cm, yshift= 0.00cm}, 
  y label style={xshift= 0.0cm, yshift=-0.30cm},
  colormap/jet,colorbar,colorbar style={title={\scriptsize $f$},
  title style={yshift=-2mm, xshift=0mm},
  width=.25cm, xshift=-0.5em},
  point meta min=\cmin,point meta max=\cmax,
  label style={font=\scriptsize},
  tick label style={font=\scriptsize},
  legend style={font=\scriptsize\selectfont},
]
\addplot [forget plot] graphics [xmin=\xmingb,xmax=\xmaxgb,ymin=\zmingb,ymax=\zmaxgb] {{\modelfile}.png};


\end{axis}
\end{tikzpicture}%
    \caption{Using single-frequency,\\ $\omega/(2\pi)=\num{1.2}$
            (PSNR \num{19.90}).}
  \end{subfigure}\hfill
  \renewcommand{\modelfile}{images/fwi/model_shapes/inversion/fwi_cp1_f0.5_unifreq1.4hz_f_scale-01-05}
  \begin{subfigure}[t]{.32\textwidth}\centering
    \begin{tikzpicture}

\pgfmathsetmacro{\xmingb} {-20} \pgfmathsetmacro{\xmaxgb}{20}
\pgfmathsetmacro{\zmingb} {-20} \pgfmathsetmacro{\zmaxgb}{20}

\begin{axis}[
  width=\modelwidth, height=\modelheight,
  axis on top, separate axis lines,
  xmin=\xminloc, xmax=\xmaxloc, 
  ymin=\zminloc, ymax=\zmaxloc, 
  x label style={xshift=-0.0cm, yshift= 0.00cm}, 
  y label style={xshift= 0.0cm, yshift=-0.30cm},
  colormap/jet,colorbar,colorbar style={title={\scriptsize $f$},
  title style={yshift=-2mm, xshift=0mm},
  width=.25cm, xshift=-0.5em},
  point meta min=\cmin,point meta max=\cmax,
  label style={font=\scriptsize},
  tick label style={font=\scriptsize},
  legend style={font=\scriptsize\selectfont},
]
\addplot [forget plot] graphics [xmin=\xmingb,xmax=\xmaxgb,ymin=\zmingb,ymax=\zmaxgb] {{\modelfile}.png};


\end{axis}
\end{tikzpicture}%
    \caption{Using single-frequency,\\ $\omega/(2\pi)=\num{1.4}$
            (PSNR \num{20.10}).}
  \end{subfigure}\hfill
  \renewcommand{\modelfile}{images/fwi/model_shapes/inversion/fwi_cp1_f0.5_multifreq0.7-1.4hz_f_scale-01-05}
  \begin{subfigure}[t]{.32\textwidth}\centering
    \begin{tikzpicture}

\pgfmathsetmacro{\xmingb} {-20} \pgfmathsetmacro{\xmaxgb}{20}
\pgfmathsetmacro{\zmingb} {-20} \pgfmathsetmacro{\zmaxgb}{20}

\begin{axis}[
  width=\modelwidth, height=\modelheight,
  axis on top, separate axis lines,
  xmin=\xminloc, xmax=\xmaxloc, 
  ymin=\zminloc, ymax=\zmaxloc, 
  x label style={xshift=-0.0cm, yshift= 0.00cm}, 
  y label style={xshift= 0.0cm, yshift=-0.30cm},
  colormap/jet,colorbar,colorbar style={title={\scriptsize $f$},
  title style={yshift=-2mm, xshift=0mm},
  width=.25cm, xshift=-0.5em},
  point meta min=\cmin,point meta max=\cmax,
  label style={font=\scriptsize},
  tick label style={font=\scriptsize},
  legend style={font=\scriptsize\selectfont},
]
\addplot [forget plot] graphics [xmin=\xmingb,xmax=\xmaxgb,ymin=\zmingb,ymax=\zmaxgb] {{\modelfile}.png};


\end{axis}
\end{tikzpicture}%
    \caption{Using multi-frequency, $\omega/(2\pi)\in\{\num{0.7}, \num{1}, \num{1.2}, \num{1.4} \}$
            (PSNR \num{23.04}).}
  \end{subfigure}
  \caption{Reconstruction with FWI starting 
           from a homogeneous background 
           with $f=0$. The models are given
           at frequency $\omega/(2\pi)=1$ and the wave speed
           is equal to \num{1} in the background.}
\label{fig:fwi:shape-model_f05_fwi}
\end{figure}

In \autoref{fig:fwi:shape-model_f05_fwi} and 
\autoref{fig:fwi:shape-model_f2_fwi}, we show 
the model perturbation, which consist in small
objects buried in the background. 
FWI is carried out with single and multiple 
frequencies, while we investigate a mild 
contrast in \autoref{fig:fwi:shape-model_f05_fwi} 
(where $f$ is at most \num{0.5}) and a 
stronger contrast in \autoref{fig:fwi:shape-model_f05_fwi} 
(where $f$ is at most \num{2}).
We see that the model can be recovered using a single frequency, which has to be
selected depending on the contrast. As an alternative, multi-frequency data appears
to be a robust candidate and always provides a good reconstruction, for both the 
object's shape and amplitude. 
The reconstruction quality with high contrast in \autoref{fig:fwi:shape-model_f2_fwi} seems to be of a similar level as with low contrast.

\begin{figure}[ht!] \centering
  \pgfmathsetmacro{\xminloc}{-5}\pgfmathsetmacro{\xmaxloc}{5}
  \pgfmathsetmacro{\zminloc}{-5}\pgfmathsetmacro{\zmaxloc}{5}
  \pgfmathsetmacro{\cmin} {-0.1} \pgfmathsetmacro{\cmax}  {2}
  \renewcommand{\modelfile}{images/fwi/model_shapes/models/cp1_f2_f_scale-01-2}
  \begin{subfigure}[t]{.32\textwidth}\centering
    \begin{tikzpicture}

\pgfmathsetmacro{\xmingb} {-20} \pgfmathsetmacro{\xmaxgb}{20}
\pgfmathsetmacro{\zmingb} {-20} \pgfmathsetmacro{\zmaxgb}{20}

\begin{axis}[
  width=\modelwidth, height=\modelheight,
  axis on top, separate axis lines,
  xmin=\xminloc, xmax=\xmaxloc, 
  ymin=\zminloc, ymax=\zmaxloc, 
  x label style={xshift=-0.0cm, yshift= 0.00cm}, 
  y label style={xshift= 0.0cm, yshift=-0.30cm},
  colormap/jet,colorbar,colorbar style={title={\scriptsize $f$},
  title style={yshift=-2mm, xshift=0mm},
  width=.25cm, xshift=-0.5em},
  point meta min=\cmin,point meta max=\cmax,
  label style={font=\scriptsize},
  tick label style={font=\scriptsize},
  legend style={font=\scriptsize\selectfont},
]
\addplot [forget plot] graphics [xmin=\xmingb,xmax=\xmaxgb,ymin=\zmingb,ymax=\zmaxgb] {{\modelfile}.png};


\end{axis}
\end{tikzpicture}%
    \caption{True contrast.}
  \end{subfigure}\hfill
  \renewcommand{\modelfile}{images/fwi/model_shapes/inversion/fwi_cp1_f2_unifreq0.7hz_f_scale-01-2}
  \begin{subfigure}[t]{.32\textwidth}\centering
    \begin{tikzpicture}

\pgfmathsetmacro{\xmingb} {-20} \pgfmathsetmacro{\xmaxgb}{20}
\pgfmathsetmacro{\zmingb} {-20} \pgfmathsetmacro{\zmaxgb}{20}

\begin{axis}[
  width=\modelwidth, height=\modelheight,
  axis on top, separate axis lines,
  xmin=\xminloc, xmax=\xmaxloc, 
  ymin=\zminloc, ymax=\zmaxloc, 
  x label style={xshift=-0.0cm, yshift= 0.00cm}, 
  y label style={xshift= 0.0cm, yshift=-0.30cm},
  colormap/jet,colorbar,colorbar style={title={\scriptsize $f$},
  title style={yshift=-2mm, xshift=0mm},
  width=.25cm, xshift=-0.5em},
  point meta min=\cmin,point meta max=\cmax,
  label style={font=\scriptsize},
  tick label style={font=\scriptsize},
  legend style={font=\scriptsize\selectfont},
]
\addplot [forget plot] graphics [xmin=\xmingb,xmax=\xmaxgb,ymin=\zmingb,ymax=\zmaxgb] {{\modelfile}.png};


\end{axis}
\end{tikzpicture}%
    \caption{Using single-frequency,\\ $\omega/(2\pi)=\num{0.7}$
            (PSNR \num{20.75}).}
  \end{subfigure}
  \renewcommand{\modelfile}{images/fwi/model_shapes/inversion/fwi_cp1_f2_unifreq1hz_f_scale-01-2}
  \begin{subfigure}[t]{.32\textwidth}\centering
    \begin{tikzpicture}

\pgfmathsetmacro{\xmingb} {-20} \pgfmathsetmacro{\xmaxgb}{20}
\pgfmathsetmacro{\zmingb} {-20} \pgfmathsetmacro{\zmaxgb}{20}

\begin{axis}[
  width=\modelwidth, height=\modelheight,
  axis on top, separate axis lines,
  xmin=\xminloc, xmax=\xmaxloc, 
  ymin=\zminloc, ymax=\zmaxloc, 
  x label style={xshift=-0.0cm, yshift= 0.00cm}, 
  y label style={xshift= 0.0cm, yshift=-0.30cm},
  colormap/jet,colorbar,colorbar style={title={\scriptsize $f$},
  title style={yshift=-2mm, xshift=0mm},
  width=.25cm, xshift=-0.5em},
  point meta min=\cmin,point meta max=\cmax,
  label style={font=\scriptsize},
  tick label style={font=\scriptsize},
  legend style={font=\scriptsize\selectfont},
]
\addplot [forget plot] graphics [xmin=\xmingb,xmax=\xmaxgb,ymin=\zmingb,ymax=\zmaxgb] {{\modelfile}.png};


\end{axis}
\end{tikzpicture}%
    \caption{Using single-frequency,\\ $\omega/(2\pi)=\num{1}$
            (PSNR \num{21.50}).}
  \end{subfigure}

  \renewcommand{\modelfile}{images/fwi/model_shapes/inversion/fwi_cp1_f2_unifreq1.2hz_f_scale-01-2}
  \begin{subfigure}[t]{.32\textwidth}\centering
    \begin{tikzpicture}

\pgfmathsetmacro{\xmingb} {-20} \pgfmathsetmacro{\xmaxgb}{20}
\pgfmathsetmacro{\zmingb} {-20} \pgfmathsetmacro{\zmaxgb}{20}

\begin{axis}[
  width=\modelwidth, height=\modelheight,
  axis on top, separate axis lines,
  xmin=\xminloc, xmax=\xmaxloc, 
  ymin=\zminloc, ymax=\zmaxloc, 
  x label style={xshift=-0.0cm, yshift= 0.00cm}, 
  y label style={xshift= 0.0cm, yshift=-0.30cm},
  colormap/jet,colorbar,colorbar style={title={\scriptsize $f$},
  title style={yshift=-2mm, xshift=0mm},
  width=.25cm, xshift=-0.5em},
  point meta min=\cmin,point meta max=\cmax,
  label style={font=\scriptsize},
  tick label style={font=\scriptsize},
  legend style={font=\scriptsize\selectfont},
]
\addplot [forget plot] graphics [xmin=\xmingb,xmax=\xmaxgb,ymin=\zmingb,ymax=\zmaxgb] {{\modelfile}.png};


\end{axis}
\end{tikzpicture}%
    \caption{Using single-frequency,\\ $\omega/(2\pi)=\num{1.2}$
            (PSNR \num{22.27}).}
  \end{subfigure}\hfill
  \renewcommand{\modelfile}{images/fwi/model_shapes/inversion/fwi_cp1_f2_unifreq1.4hz_f_scale-01-2}
  \begin{subfigure}[t]{.32\textwidth}\centering
    \begin{tikzpicture}

\pgfmathsetmacro{\xmingb} {-20} \pgfmathsetmacro{\xmaxgb}{20}
\pgfmathsetmacro{\zmingb} {-20} \pgfmathsetmacro{\zmaxgb}{20}

\begin{axis}[
  width=\modelwidth, height=\modelheight,
  axis on top, separate axis lines,
  xmin=\xminloc, xmax=\xmaxloc, 
  ymin=\zminloc, ymax=\zmaxloc, 
  x label style={xshift=-0.0cm, yshift= 0.00cm}, 
  y label style={xshift= 0.0cm, yshift=-0.30cm},
  colormap/jet,colorbar,colorbar style={title={\scriptsize $f$},
  title style={yshift=-2mm, xshift=0mm},
  width=.25cm, xshift=-0.5em},
  point meta min=\cmin,point meta max=\cmax,
  label style={font=\scriptsize},
  tick label style={font=\scriptsize},
  legend style={font=\scriptsize\selectfont},
]
\addplot [forget plot] graphics [xmin=\xmingb,xmax=\xmaxgb,ymin=\zmingb,ymax=\zmaxgb] {{\modelfile}.png};


\end{axis}
\end{tikzpicture}%
    \caption{Using single-frequency,\\ $\omega/(2\pi)=\num{1.4}$
            (PSNR \num{22.72}).}
  \end{subfigure}\hfill
  \renewcommand{\modelfile}{images/fwi/model_shapes/inversion/fwi_cp1_f2_multifreq0.7-1.4hz_f_scale-01-2}
  \begin{subfigure}[t]{.32\textwidth}\centering
    \begin{tikzpicture}

\pgfmathsetmacro{\xmingb} {-20} \pgfmathsetmacro{\xmaxgb}{20}
\pgfmathsetmacro{\zmingb} {-20} \pgfmathsetmacro{\zmaxgb}{20}

\begin{axis}[
  width=\modelwidth, height=\modelheight,
  axis on top, separate axis lines,
  xmin=\xminloc, xmax=\xmaxloc, 
  ymin=\zminloc, ymax=\zmaxloc, 
  x label style={xshift=-0.0cm, yshift= 0.00cm}, 
  y label style={xshift= 0.0cm, yshift=-0.30cm},
  colormap/jet,colorbar,colorbar style={title={\scriptsize $f$},
  title style={yshift=-2mm, xshift=0mm},
  width=.25cm, xshift=-0.5em},
  point meta min=\cmin,point meta max=\cmax,
  label style={font=\scriptsize},
  tick label style={font=\scriptsize},
  legend style={font=\scriptsize\selectfont},
]
\addplot [forget plot] graphics [xmin=\xmingb,xmax=\xmaxgb,ymin=\zmingb,ymax=\zmaxgb] {{\modelfile}.png};


\end{axis}
\end{tikzpicture}%
    \caption{Using multi-frequency, 
             $\omega/(2\pi)\in\{\num{0.7}, \num{1}, \num{1.2}, \num{1.4} \}$
            (PSNR \num{22.82}).}
  \end{subfigure}
  \caption{Reconstruction with FWI starting 
           from a homogeneous background 
           with $f=0$. The models are given
           at frequency $\omega/(2\pi)=1$ and the wave speed
           is equal to \num{1} in the background.}
  \label{fig:fwi:shape-model_f2_fwi}
\end{figure}

\subsubsection{Reconstruction using Born and Rytov approximations}

\begin{figure}[ht]\centering
  \pgfmathsetmacro{\xminloc}{-7.0711}\pgfmathsetmacro{\xmaxloc}{7.0711}
  \pgfmathsetmacro{\zminloc}{-7.0711}\pgfmathsetmacro{\zmaxloc}{7.0711}
  \pgfmathsetmacro{\xminpic}{-5.0000}\pgfmathsetmacro{\xmaxpic}{5.0000}
  \pgfmathsetmacro{\zminpic}{-5.0000}\pgfmathsetmacro{\zmaxpic}{5.0000}
  \pgfmathsetmacro{\cmin} {-0.1} \pgfmathsetmacro{\cmax} { 0.5}
  
  \renewcommand{\modelfile}{images/shapes/shapesnoise-1000mHz}
  \begin{subfigure}[t]{.32\textwidth}
    \begin{tikzpicture}
  \begin{axis}[
    width=\modelwidth, height=\modelheight,
    axis on top, separate axis lines,
    xmin=\xminpic, xmax=\xmaxpic, 
    ymin=\zminpic, ymax=\zmaxpic, 
    x label style={xshift=-0.0cm, yshift= 0.00cm}, 
    y label style={xshift= 0.0cm, yshift=-0.30cm},
    colormap/jet,colorbar,colorbar style={
      width=.25cm, xshift=-.5em},
    point meta min=\cmin,point meta max=\cmax,
    label style={font=\scriptsize},
    tick label style={font=\scriptsize},
    legend style={font=\scriptsize\selectfont},
    ]
    \addplot [forget plot] graphics [xmin=\xminloc,xmax=\xmaxloc,ymin=\zminloc,ymax=\zmaxloc] {{\modelfile}.png};
  \end{axis}
\end{tikzpicture}%
    \caption{True model $f$}
  \end{subfigure}\hfill
  \renewcommand{\modelfile}{images/shapes/shapesnoise-1000mHz-rec}
  \begin{subfigure}[t]{.32\textwidth}
    \begin{tikzpicture}
  \begin{axis}[
    width=\modelwidth, height=\modelheight,
    axis on top, separate axis lines,
    xmin=\xminpic, xmax=\xmaxpic, 
    ymin=\zminpic, ymax=\zmaxpic, 
    x label style={xshift=-0.0cm, yshift= 0.00cm}, 
    y label style={xshift= 0.0cm, yshift=-0.30cm},
    colormap/jet,colorbar,colorbar style={
      width=.25cm, xshift=-.5em},
    point meta min=\cmin,point meta max=\cmax,
    label style={font=\scriptsize},
    tick label style={font=\scriptsize},
    legend style={font=\scriptsize\selectfont},
    ]
    \addplot [forget plot] graphics [xmin=\xminloc,xmax=\xmaxloc,ymin=\zminloc,ymax=\zmaxloc] {{\modelfile}.png};
  \end{axis}
\end{tikzpicture}%
    \caption{Born reconstruction at frequency $\omega/(2\pi)={1} $. (PSNR 20.14)}
  \end{subfigure}\hfill
  \renewcommand{\modelfile}{images/shapes/shapesnoise-1000mHz-rec-Rytov}
  \begin{subfigure}[t]{.32\textwidth}
    \begin{tikzpicture}
  \begin{axis}[
    width=\modelwidth, height=\modelheight,
    axis on top, separate axis lines,
    xmin=\xminpic, xmax=\xmaxpic, 
    ymin=\zminpic, ymax=\zmaxpic, 
    x label style={xshift=-0.0cm, yshift= 0.00cm}, 
    y label style={xshift= 0.0cm, yshift=-0.30cm},
    colormap/jet,colorbar,colorbar style={
      width=.25cm, xshift=-.5em},
    point meta min=\cmin,point meta max=\cmax,
    label style={font=\scriptsize},
    tick label style={font=\scriptsize},
    legend style={font=\scriptsize\selectfont},
    ]
    \addplot [forget plot] graphics [xmin=\xminloc,xmax=\xmaxloc,ymin=\zminloc,ymax=\zmaxloc] {{\modelfile}.png};
  \end{axis}
\end{tikzpicture}%
    \caption{Rytov reconstruction at frequency $\omega/(2\pi)={1} $. (PSNR 20.10)}
  \end{subfigure}

  \renewcommand{\modelfile}{images/shapes/shapesnoise-700mHz-rec-Rytov}
  \begin{subfigure}[t]{.32\textwidth}
    \begin{tikzpicture}
  \begin{axis}[
    width=\modelwidth, height=\modelheight,
    axis on top, separate axis lines,
    xmin=\xminpic, xmax=\xmaxpic, 
    ymin=\zminpic, ymax=\zmaxpic, 
    x label style={xshift=-0.0cm, yshift= 0.00cm}, 
    y label style={xshift= 0.0cm, yshift=-0.30cm},
    colormap/jet,colorbar,colorbar style={
      width=.25cm, xshift=-.5em},
    point meta min=\cmin,point meta max=\cmax,
    label style={font=\scriptsize},
    tick label style={font=\scriptsize},
    legend style={font=\scriptsize\selectfont},
    ]
    \addplot [forget plot] graphics [xmin=\xminloc,xmax=\xmaxloc,ymin=\zminloc,ymax=\zmaxloc] {{\modelfile}.png};
  \end{axis}
\end{tikzpicture}%
    \caption{Rytov reconstruction at frequency $\omega/(2\pi)={0.7} $. (PSNR 18.00)} 
  \end{subfigure}
  \renewcommand{\modelfile}{images/shapes/shapesnoise-1200mHz-rec-Rytov}
  \begin{subfigure}[t]{.32\textwidth}
    \begin{tikzpicture}
  \begin{axis}[
    width=\modelwidth, height=\modelheight,
    axis on top, separate axis lines,
    xmin=\xminpic, xmax=\xmaxpic, 
    ymin=\zminpic, ymax=\zmaxpic, 
    x label style={xshift=-0.0cm, yshift= 0.00cm}, 
    y label style={xshift= 0.0cm, yshift=-0.30cm},
    colormap/jet,colorbar,colorbar style={
      width=.25cm, xshift=-.5em},
    point meta min=\cmin,point meta max=\cmax,
    label style={font=\scriptsize},
    tick label style={font=\scriptsize},
    legend style={font=\scriptsize\selectfont},
    ]
    \addplot [forget plot] graphics [xmin=\xminloc,xmax=\xmaxloc,ymin=\zminloc,ymax=\zmaxloc] {{\modelfile}.png};
  \end{axis}
\end{tikzpicture}%
    \caption{Rytov reconstruction at frequency $\omega/(2\pi)={1.2} $. (PSNR 21.29)} 
  \end{subfigure}
  \renewcommand{\modelfile}{images/shapes/shapesnoise-1400mHz-rec-Rytov}
  \begin{subfigure}[t]{.32\textwidth}
    \begin{tikzpicture}
  \begin{axis}[
    width=\modelwidth, height=\modelheight,
    axis on top, separate axis lines,
    xmin=\xminpic, xmax=\xmaxpic, 
    ymin=\zminpic, ymax=\zmaxpic, 
    x label style={xshift=-0.0cm, yshift= 0.00cm}, 
    y label style={xshift= 0.0cm, yshift=-0.30cm},
    colormap/jet,colorbar,colorbar style={
      width=.25cm, xshift=-.5em},
    point meta min=\cmin,point meta max=\cmax,
    label style={font=\scriptsize},
    tick label style={font=\scriptsize},
    legend style={font=\scriptsize\selectfont},
    ]
    \addplot [forget plot] graphics [xmin=\xminloc,xmax=\xmaxloc,ymin=\zminloc,ymax=\zmaxloc] {{\modelfile}.png};
  \end{axis}
\end{tikzpicture}%
    \caption{Rytov reconstruction at frequency $\omega/(2\pi)={1.4} $. (PSNR 22.17)} 
  \end{subfigure}

  \renewcommand{\modelfile}{images/shapes/shapesnoise-700_1000_1200_1400mHz-rec}
  \begin{subfigure}[t]{.375\textwidth}\centering
  \begin{tikzpicture}
  \begin{axis}[
    width=\modelwidth, height=\modelheight,
    axis on top, separate axis lines,
    xmin=\xminpic, xmax=\xmaxpic, 
    ymin=\zminpic, ymax=\zmaxpic, 
    x label style={xshift=-0.0cm, yshift= 0.00cm}, 
    y label style={xshift= 0.0cm, yshift=-0.30cm},
    colormap/jet,colorbar,colorbar style={
      width=.25cm, xshift=-.5em},
    point meta min=\cmin,point meta max=\cmax,
    label style={font=\scriptsize},
    tick label style={font=\scriptsize},
    legend style={font=\scriptsize\selectfont},
    ]
    \addplot [forget plot] graphics [xmin=\xminloc,xmax=\xmaxloc,ymin=\zminloc,ymax=\zmaxloc] {{\modelfile}.png};
  \end{axis}
\end{tikzpicture}%
  \caption{Born reconstruction using multi-frequency, $\omega/(2\pi)\in\{0.7,1,1.2,1.4\} $. (PSNR 21.91)} 
  \end{subfigure}\quad
  \renewcommand{\modelfile}{images/shapes/shapesnoise-700_1000_1200_1400mHz-rec-Rytov}
  \begin{subfigure}[t]{.375\textwidth}\centering
  \begin{tikzpicture}
  \begin{axis}[
    width=\modelwidth, height=\modelheight,
    axis on top, separate axis lines,
    xmin=\xminpic, xmax=\xmaxpic, 
    ymin=\zminpic, ymax=\zmaxpic, 
    x label style={xshift=-0.0cm, yshift= 0.00cm}, 
    y label style={xshift= 0.0cm, yshift=-0.30cm},
    colormap/jet,colorbar,colorbar style={
      width=.25cm, xshift=-.5em},
    point meta min=\cmin,point meta max=\cmax,
    label style={font=\scriptsize},
    tick label style={font=\scriptsize},
    legend style={font=\scriptsize\selectfont},
    ]
    \addplot [forget plot] graphics [xmin=\xminloc,xmax=\xmaxloc,ymin=\zminloc,ymax=\zmaxloc] {{\modelfile}.png};
  \end{axis}
\end{tikzpicture}%
  \caption{Rytov reconstruction using multi-frequency, $\omega/(2\pi)\in\{0.7,1,1.2,1.4\} $. (PSNR 21.93)} 
  \end{subfigure}
  
  \caption{Reconstructions of the more complicated shapes.
    The models are given at frequency $\omega/(2\pi)=1$.
    \label{fig:shapes-model-rec}}
\end{figure}
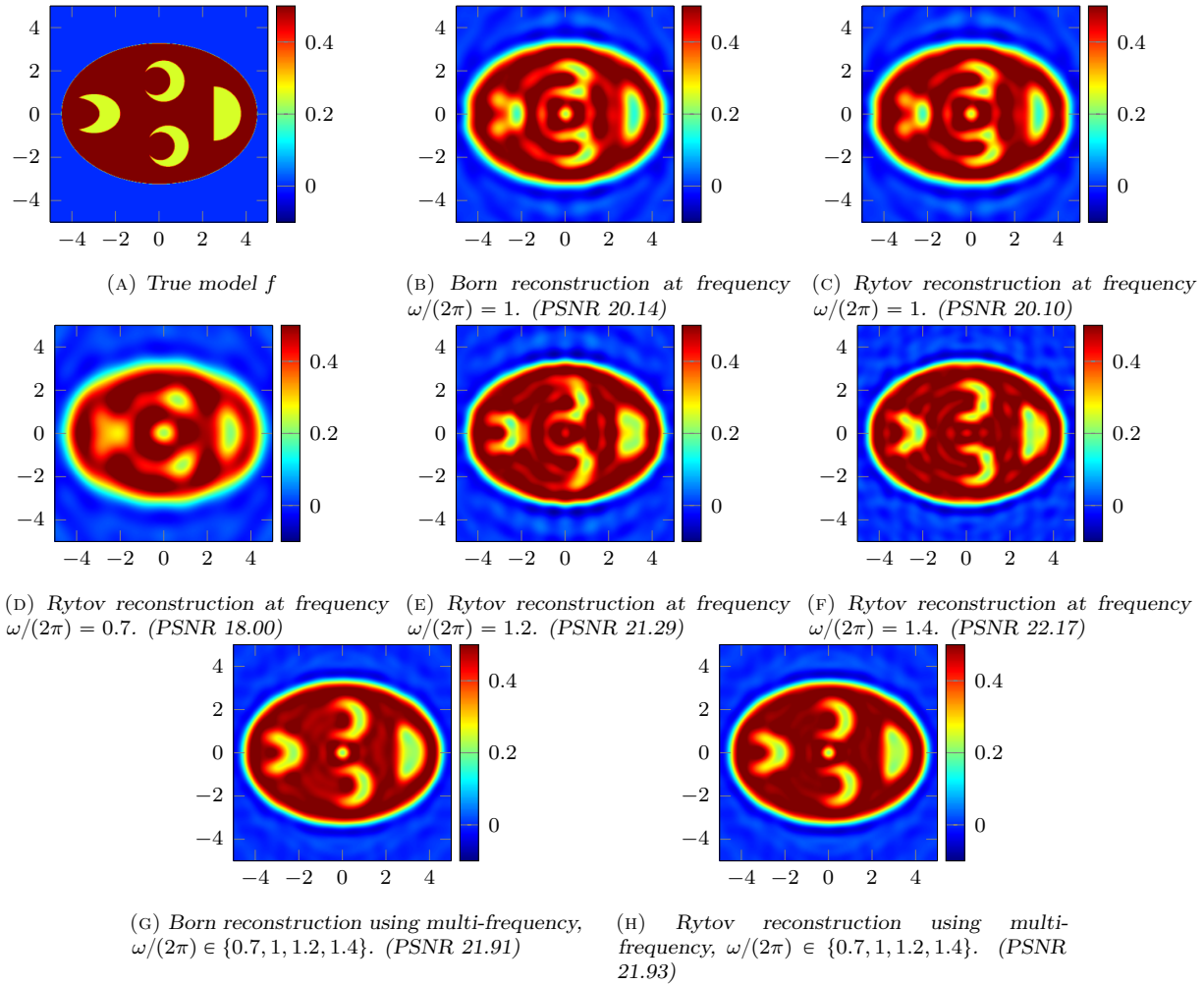

In \autoref{fig:shapes-model-rec}, we show the reconstruction 
using Born and Rytov approximations.
Here, we can clearly see that we need a higher frequency in 
order to resolve small features of the object.
However, the reconstructions are still inferior to the FWI.

\begin{figure}[ht]\centering
  \pgfmathsetmacro{\xminloc}{-7.0711}\pgfmathsetmacro{\xmaxloc}{7.0711}
  \pgfmathsetmacro{\zminloc}{-7.0711}\pgfmathsetmacro{\zmaxloc}{7.0711}
  \pgfmathsetmacro{\xminpic}{-5.0000}\pgfmathsetmacro{\xmaxpic}{5.0000}
  \pgfmathsetmacro{\zminpic}{-5.0000}\pgfmathsetmacro{\zmaxpic}{5.0000}
  \pgfmathsetmacro{\cmin} {-0.1} \pgfmathsetmacro{\cmax} { 2}
  
  \renewcommand{\modelfile}{images/shapes/shapes2noise-1000mHz}
  \begin{subfigure}[t]{.32\textwidth}
    \begin{tikzpicture}
  \begin{axis}[
    width=\modelwidth, height=\modelheight,
    axis on top, separate axis lines,
    xmin=\xminpic, xmax=\xmaxpic, 
    ymin=\zminpic, ymax=\zmaxpic, 
    x label style={xshift=-0.0cm, yshift= 0.00cm}, 
    y label style={xshift= 0.0cm, yshift=-0.30cm},
    colormap/jet,colorbar,colorbar style={
      width=.25cm, xshift=-.5em},
    point meta min=\cmin,point meta max=\cmax,
    label style={font=\scriptsize},
    tick label style={font=\scriptsize},
    legend style={font=\scriptsize\selectfont},
    ]
    \addplot [forget plot] graphics [xmin=\xminloc,xmax=\xmaxloc,ymin=\zminloc,ymax=\zmaxloc] {{\modelfile}.png};
  \end{axis}
\end{tikzpicture}%
    \caption{True model $f$}
  \end{subfigure}\hfill
  \renewcommand{\modelfile}{images/shapes/shapes2noise-1000mHz-rec}
  \begin{subfigure}[t]{.32\textwidth}
    \begin{tikzpicture}
  \begin{axis}[
    width=\modelwidth, height=\modelheight,
    axis on top, separate axis lines,
    xmin=\xminpic, xmax=\xmaxpic, 
    ymin=\zminpic, ymax=\zmaxpic, 
    x label style={xshift=-0.0cm, yshift= 0.00cm}, 
    y label style={xshift= 0.0cm, yshift=-0.30cm},
    colormap/jet,colorbar,colorbar style={
      width=.25cm, xshift=-.5em},
    point meta min=\cmin,point meta max=\cmax,
    label style={font=\scriptsize},
    tick label style={font=\scriptsize},
    legend style={font=\scriptsize\selectfont},
    ]
    \addplot [forget plot] graphics [xmin=\xminloc,xmax=\xmaxloc,ymin=\zminloc,ymax=\zmaxloc] {{\modelfile}.png};
  \end{axis}
\end{tikzpicture}%
    \caption{Born reconstruction at frequency $\omega/(2\pi)={1} $. (PSNR 17.16)}
  \end{subfigure}\hfill
  \renewcommand{\modelfile}{images/shapes/shapes2noise-1000mHz-rec-Rytov}
  \begin{subfigure}[t]{.32\textwidth}
    \begin{tikzpicture}
  \begin{axis}[
    width=\modelwidth, height=\modelheight,
    axis on top, separate axis lines,
    xmin=\xminpic, xmax=\xmaxpic, 
    ymin=\zminpic, ymax=\zmaxpic, 
    x label style={xshift=-0.0cm, yshift= 0.00cm}, 
    y label style={xshift= 0.0cm, yshift=-0.30cm},
    colormap/jet,colorbar,colorbar style={
      width=.25cm, xshift=-.5em},
    point meta min=\cmin,point meta max=\cmax,
    label style={font=\scriptsize},
    tick label style={font=\scriptsize},
    legend style={font=\scriptsize\selectfont},
    ]
    \addplot [forget plot] graphics [xmin=\xminloc,xmax=\xmaxloc,ymin=\zminloc,ymax=\zmaxloc] {{\modelfile}.png};
  \end{axis}
\end{tikzpicture}%
    \caption{Rytov reconstruction at frequency $\omega/(2\pi)={1} $. (PSNR 18.73)}
  \end{subfigure}

  \renewcommand{\modelfile}{images/shapes/shapes2noise-700mHz-rec-Rytov}
  \begin{subfigure}[t]{.32\textwidth}
    \begin{tikzpicture}
  \begin{axis}[
    width=\modelwidth, height=\modelheight,
    axis on top, separate axis lines,
    xmin=\xminpic, xmax=\xmaxpic, 
    ymin=\zminpic, ymax=\zmaxpic, 
    x label style={xshift=-0.0cm, yshift= 0.00cm}, 
    y label style={xshift= 0.0cm, yshift=-0.30cm},
    colormap/jet,colorbar,colorbar style={
      width=.25cm, xshift=-.5em},
    point meta min=\cmin,point meta max=\cmax,
    label style={font=\scriptsize},
    tick label style={font=\scriptsize},
    legend style={font=\scriptsize\selectfont},
    ]
    \addplot [forget plot] graphics [xmin=\xminloc,xmax=\xmaxloc,ymin=\zminloc,ymax=\zmaxloc] {{\modelfile}.png};
  \end{axis}
\end{tikzpicture}%
    \caption{Rytov reconstruction at frequency $\omega/(2\pi)={0.7} $. (PSNR 16.10)} 
  \end{subfigure}\hfill
  \renewcommand{\modelfile}{images/shapes/shapes2noise-1200mHz-rec-Rytov}
  \begin{subfigure}[t]{.32\textwidth}
    \begin{tikzpicture}
  \begin{axis}[
    width=\modelwidth, height=\modelheight,
    axis on top, separate axis lines,
    xmin=\xminpic, xmax=\xmaxpic, 
    ymin=\zminpic, ymax=\zmaxpic, 
    x label style={xshift=-0.0cm, yshift= 0.00cm}, 
    y label style={xshift= 0.0cm, yshift=-0.30cm},
    colormap/jet,colorbar,colorbar style={
      width=.25cm, xshift=-.5em},
    point meta min=\cmin,point meta max=\cmax,
    label style={font=\scriptsize},
    tick label style={font=\scriptsize},
    legend style={font=\scriptsize\selectfont},
    ]
    \addplot [forget plot] graphics [xmin=\xminloc,xmax=\xmaxloc,ymin=\zminloc,ymax=\zmaxloc] {{\modelfile}.png};
  \end{axis}
\end{tikzpicture}%
    \caption{Rytov reconstruction at frequency $\omega/(2\pi)={1.2} $. (PSNR 20.14)} 
  \end{subfigure}\hfill
  \renewcommand{\modelfile}{images/shapes/shapes2noise-1400mHz-rec-Rytov}
  \begin{subfigure}[t]{.32\textwidth}
    \begin{tikzpicture}
  \begin{axis}[
    width=\modelwidth, height=\modelheight,
    axis on top, separate axis lines,
    xmin=\xminpic, xmax=\xmaxpic, 
    ymin=\zminpic, ymax=\zmaxpic, 
    x label style={xshift=-0.0cm, yshift= 0.00cm}, 
    y label style={xshift= 0.0cm, yshift=-0.30cm},
    colormap/jet,colorbar,colorbar style={
      width=.25cm, xshift=-.5em},
    point meta min=\cmin,point meta max=\cmax,
    label style={font=\scriptsize},
    tick label style={font=\scriptsize},
    legend style={font=\scriptsize\selectfont},
    ]
    \addplot [forget plot] graphics [xmin=\xminloc,xmax=\xmaxloc,ymin=\zminloc,ymax=\zmaxloc] {{\modelfile}.png};
  \end{axis}
\end{tikzpicture}%
    \caption{Rytov reconstruction at frequency $\omega/(2\pi)={1.4} $. (PSNR 21.15)} 
  \end{subfigure}

  \renewcommand{\modelfile}{images/shapes/shapes2noise-700_1000_1200_1400mHz-rec}
  \begin{subfigure}[t]{.375\textwidth}\centering
  \begin{tikzpicture}
  \begin{axis}[
    width=\modelwidth, height=\modelheight,
    axis on top, separate axis lines,
    xmin=\xminpic, xmax=\xmaxpic, 
    ymin=\zminpic, ymax=\zmaxpic, 
    x label style={xshift=-0.0cm, yshift= 0.00cm}, 
    y label style={xshift= 0.0cm, yshift=-0.30cm},
    colormap/jet,colorbar,colorbar style={
      width=.25cm, xshift=-.5em},
    point meta min=\cmin,point meta max=\cmax,
    label style={font=\scriptsize},
    tick label style={font=\scriptsize},
    legend style={font=\scriptsize\selectfont},
    ]
    \addplot [forget plot] graphics [xmin=\xminloc,xmax=\xmaxloc,ymin=\zminloc,ymax=\zmaxloc] {{\modelfile}.png};
  \end{axis}
\end{tikzpicture}%
  \caption{Born reconstruction using multi-frequency, $\omega/(2\pi)\in\{0.7,1,1.2,1.4\} $. (PSNR 16.92)} 
  \end{subfigure}\quad
  \renewcommand{\modelfile}{images/shapes/shapes2noise-700_1000_1200_1400mHz-rec-Rytov}
  \begin{subfigure}[t]{.375\textwidth}\centering
  \begin{tikzpicture}
  \begin{axis}[
    width=\modelwidth, height=\modelheight,
    axis on top, separate axis lines,
    xmin=\xminpic, xmax=\xmaxpic, 
    ymin=\zminpic, ymax=\zmaxpic, 
    x label style={xshift=-0.0cm, yshift= 0.00cm}, 
    y label style={xshift= 0.0cm, yshift=-0.30cm},
    colormap/jet,colorbar,colorbar style={
      width=.25cm, xshift=-.5em},
    point meta min=\cmin,point meta max=\cmax,
    label style={font=\scriptsize},
    tick label style={font=\scriptsize},
    legend style={font=\scriptsize\selectfont},
    ]
    \addplot [forget plot] graphics [xmin=\xminloc,xmax=\xmaxloc,ymin=\zminloc,ymax=\zmaxloc] {{\modelfile}.png};
  \end{axis}
\end{tikzpicture}%
  \caption{Rytov reconstruction using multi-frequency, $\omega/(2\pi)\in\{0.7,1,1.2,1.4\} $. (PSNR 20.33)} 
  \end{subfigure}
  
  \caption{Reconstructions of the more complicated shapes with a higher contrast than in \autoref{fig:shapes-model-rec}.
    The models are given at frequency $\omega/(2\pi)=1$.
    \label{fig:shapes2-model-rec}}
\end{figure}
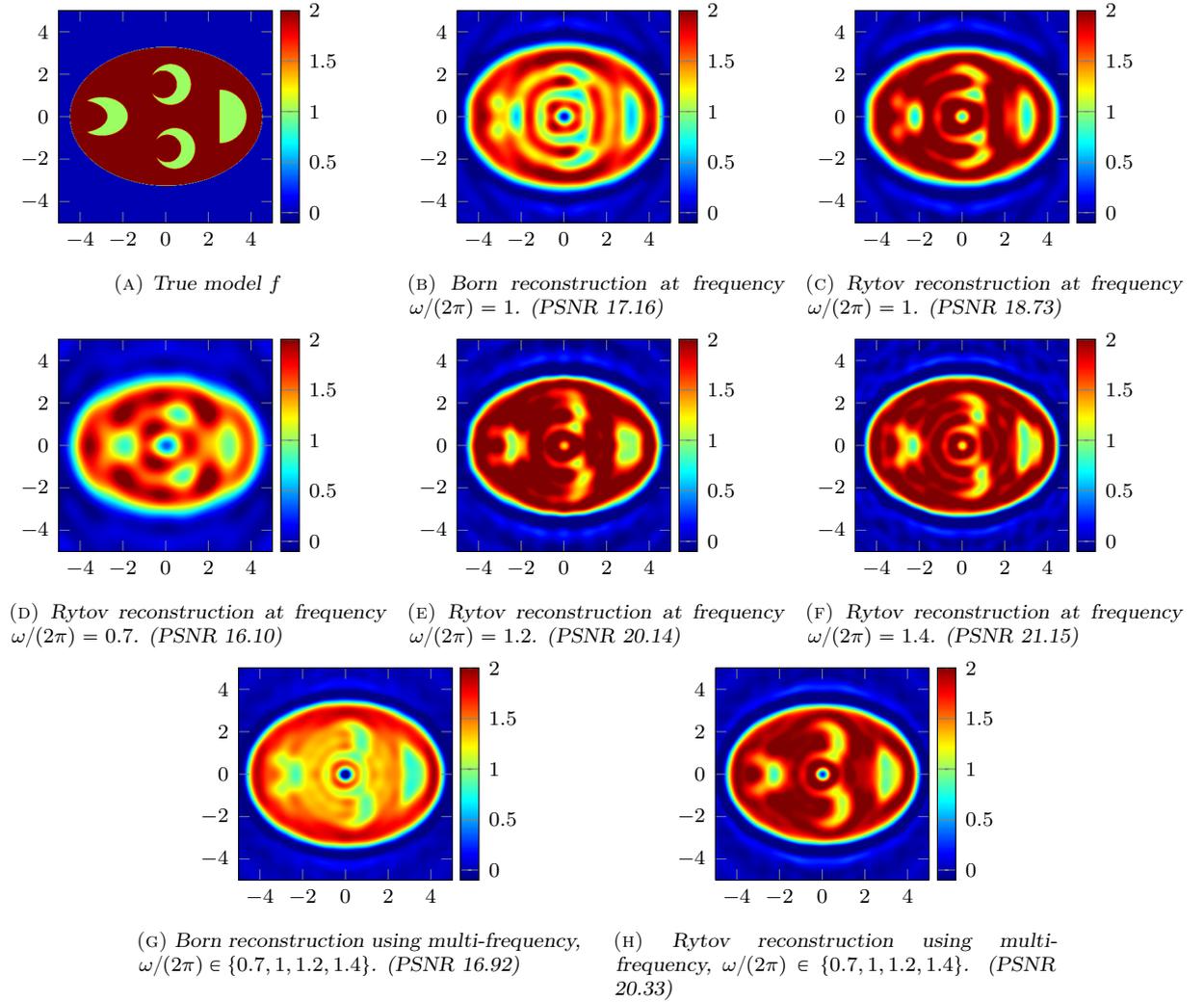

With a higher contrast, the reconstruction with the Rytov approximation is considerably better than the one with the Born approximation, see \autoref{fig:shapes2-model-rec}. 
This is consistent with \autoref{rem:Born-Rytov}.
Interestingly, the shapes reconstruction in high contrast barely profits from taking frequencies higher than 1, even though the k-space coverage is larger.
In this situation, the Rytov reconstruction is almost comparable to the one with lower contrast, but still worse than the FWI.

\subsection{Computational costs}
\label{subsection:computational-costs} 
\hspace*{\fill}

\textbf{Computational cost of FWI.}
The computational cost of FWI comes from the discretization 
and resolution of the wave problem \autoref{eq:fwi:helmholtz-full}
for each of the sources in the acquisition, coupled with the 
iterative procedure of \autoref{algo:FWI}.
In our numerical experiments, we use the software \texttt{hawen} for the 
iterative inversion, \cite{Hawen2020},~\autoref{footnote:fwi:hawen},
which relies on the Hybridizable discontinuous Galerkin 
discretization, \cite{Cockburn2009,Faucher2020adjoint}.
The number of degrees of freedom for the discretization 
depends on the number of cells in the mesh, and the 
polynomial order. 
In the inversion experiments, we use a fixed mesh for
all iterations, with about fifty thousands cells. On
the other hand, the polynomial order is selected depending 
on the wavelength on each cell. That is, each of the cells 
in the mesh is allowed to have a different order (here 
between $3$ to $7$), see \cite{Faucher2020adjoint}.
Then, when the frequency changes, while the mesh remains 
the same, the order of the polynomial evolves accordingly
to the change of wavelength. 
Once the wave equation, \autoref{eq:fwi:helmholtz-full}, is 
discretized, we obtain a linear system which size is the 
number of degrees of freedom, that must be solved for the 
different sources (i.e., the different incident angles). 
We rely on the direct solver MUMPS, \cite{Amestoy2019}, such
that once the matrix factorization is computed, the numerical 
cost of having several sources (i.e., multiple right-hand sides 
in the linear system) is drastically mitigated, motivating
the use of a direct solver instead of an iterative one.

Our numerical experiments have been carried out on the 
Vienna Scientific Cluster \texttt{VSC-4}\footnote{\url{https://vsc.ac.at/}},
using \num{48} cores. 
For the reconstructions of \autoref{fig:fwi:heart-model_f05_fwi}, 
\autoref{fig:fwi:heart-model_f2_fwi}, \autoref{fig:fwi:shape-model_f05_fwi}
and \autoref{fig:fwi:shape-model_f2_fwi}, the size of the computational
domain is \num{40}$\times$\num{40}, with about \num{350000} degrees of 
freedom.
Using single-frequency data, \num{50} iterations are performed
in \autoref{algo:FWI} and the total computational time is of 
about \SI{40}{\min}.
In the case of multiple frequencies, we have a total of \num{120}
iterations and the computational time is of about \SI{1}{\hour}\,\SI{45}{\min}.

\textbf{Computational cost of Born and Rytov approximations.}
The conjugate gradient method used in the inverse NDFT requires in each iteration step the evaluation of an NDFT \autoref{eq:ndft} and its adjoint.
We utilize the nonequispaced fast Fourier transform (NFFT) algorithm \cite{KeKuPo09}, implemented in the open-source software library \texttt{nfft} \cite{nfft3}, which can compute an NDFT in $\mathcal{O}(N^2 \log N + M)$ arithmetic operations, which is considerably less than the $\mathcal{O}(N^2M)$ operations of a straightforward implementation of \autoref{eq:ndft}.

The numerical simulations of \autoref{sec:Born-heart-shapes} have been carried out on a 4-core Intel Core i5-6500 processor.
We used 12 iterations of the conjugate gradient method and noted that the reconstruction quality hardly benefits from a higher number of iterations. 
The reconstruction of an image took about one second, with a grid size $720\times720$ of $f$ and $200\times100$ data points of $u$ for each frequency.
Therefore, the numerical computation of the Born and Rytov approximations is much faster than the FWI.

\section{Conclusion} \label{sec:conclusion}

   We study the imaging problem for diffraction tomography, where 
   wave measurements are used to quantitatively reconstruct the 
   physical properties, i.e., the refractive index.
   The forward operator that describes the wave propagation corresponds
   with the Helmholtz equation, which, under the assumption of small 
   background perturbations, can be represented via the Born and Rytov 
   approximations.
   
   Firstly, we have compared different forward models in terms of the resulting
   measured data $u$. It highlights that, even in the case of a small circular 
   object, the Born approximation is not entirely accurate to represent the 
   total wave field given by the Helmholtz equation.
   In addition, the source that initiates the phenomenon (e.g., a point source 
   located very far from the object, or simultaneous point source along a line), 
   also plays an important role as it changes the resulting signals, 
   hence leading to systematic differences depending on the choice of forward model. 
   We found that the line source model approximates the plane wave pretty well.
   
   Secondly, we have carried out the reconstruction using data from the total field $\utot$, 
   and compared the efficiency of the Full Waveform Inversion method (FWI) with that
   of the Born and Rytov approximations.
   FWI works directly with the Helmholtz problem, \autoref{eq:fwi:helmholtz-full}, 
   hence giving a robust approach that can be implemented in all configurations, 
   however at the cost of possibly intensive computations.
   On the other hand, the Born and Rytov are computationally cheap, but lack 
   accuracy when the object is too large or when the contrast is too strong.
   We have also noted that the Rytov approximation gives better results than the Born one. 
   Furthermore, for all reconstruction methods, we have shown that 
   using data of multiple frequencies allows to improve the robustness 
   of the reconstruction by providing information on multiple wavelengths.

\subsection*{Acknowledgments} 
  This work is supported by the Austrian Science Fund (FWF) 
  within SFB F68 (``Tomography across the Scales''), Projects 
  F68-06 and F68-07.
  FF is funded by the Austrian Science Fund (FWF) under the Lise 
  Meitner fellowship M 2791-N. 
  Funding by the DFG under Germany's Excellence Strategy -- The Berlin Mathematics Research
  Center MATH+ (EXC-2046/1, Projektnummer: 390685689) is gratefully acknowledged. 
  For the numerical experiments, we acknowledge the use of the 
  Vienna Scientific Cluster \texttt{VSC-4} (\url{https://vsc.ac.at/}).


\section*{References}
\renewcommand{\i}{\ii}
\printbibliography[heading=none]

\end{document}